\documentclass[11pt]{amsart}%
\usepackage{amsmath}%
\usepackage{amssymb}%
\usepackage{amscd}%
\usepackage{color}%
\usepackage{mathrsfs}%
\usepackage[all]{xy}%
\def\ol#1{\overline{#1}}
\def\wh#1{\widehat{#1}}
\def\wt#1{\widetilde{#1}}
\def\ul#1{\underline{#1}}
\theoremstyle{plain}
    \newtheorem{theorem}{Theorem}[section]
    
    \newtheorem{proposition}[theorem]{Proposition}
    
    \newtheorem{lemma}[theorem]{Lemma}
    \newtheorem{corollary}[theorem]{Corollary}
\theoremstyle{definition}
    \newtheorem{definition}[theorem]{Definition}

    \newtheorem{remark}[theorem]{Remark}
\def\Alphabet{A,B,C,D,E,F,G,H,I,J,K,L,M,N,O,P,Q,R,S,T,U,V,W,X,Y,Z}
\def\alphabet{a,b,c,d,e,f,g,h,i,j,k,l,m,n,o,p,q,r,s,t,u,v,w,x,y,z}
\def\endpiece{xxx}
\def\makeAlphabet[#1]{\expandafter\makeA#1,xxx,}
\def\makealphabet[#1]{\expandafter\makea#1,xxx,}
\def\makeA#1,{\def\temp{#1}\ifx\temp\endpiece\else%
\mkbb{#1}\mkfrak{#1}\mkbf{#1}\mkcal{#1}\mkscr{#1}\expandafter\makeA\fi}%
\def\makea#1,{\def\temp{#1}\ifx\temp\endpiece\else\mkfrak{#1}\mkbf{#1}\expandafter\makea\fi}%
\def\mkbb#1{\expandafter\def\csname bb#1\endcsname{\mathbb{#1}}}
\def\mkfrak#1{\expandafter\def\csname fr#1\endcsname{\mathfrak{#1}}}
\def\mkbf#1{\expandafter\def\csname b#1\endcsname{\mathbf{#1}}}
\def\mkcal#1{\expandafter\def\csname c#1\endcsname{\mathcal{#1}}}
\def\mkscr#1{\expandafter\def\csname s#1\endcsname{\mathscr{#1}}}
\def\makeop[#1]{\xmakeop#1,xxx,}
\def\mkop#1{\expandafter\def\csname #1\endcsname{{\mathrm{#1}}}} %
\def\xmakeop#1,{\def\temp{#1}\ifx\temp\endpiece\else\mkop{#1}\expandafter\xmakeop\fi}%
\makeAlphabet[\Alphabet]
\makealphabet[\alphabet]
\makeop[Alt,End,Ext,Hom,Sym,Tor]
\makeop[CH,Pic,div]
\makeop[Ker,Im,Coker,Coim,id,pr]
\makeop[Spec,Spf,Spm,Isom]
\makeop[Re,Im]%
\makeop[dR,Nis,crys,rig,syn,tor,Cont]
\makeop[Gal,Res,ab,res]
\makeop[pol]
\makeop[Var,an,Isoc,univ,Cusp]%
\makeop[Lie,coLie,MIC,DR,Gr,Ab,Cone,col,Frob]
\makeop[Meas,For]
\makeop[Gl,Sl]
\makeop[Eis,tr,Fil,gr,et,red]
\def\MHS{\operatorname{M}\!\operatorname{H}\!\operatorname{S}}%
\def\VMHS{\operatorname{V}\!\operatorname{M}\!\operatorname{H}\!\operatorname{S}}%
\def\sLog{\sL\!\!\operatorname{og}}
\def\cLog{\cL\!\operatorname{og}}

\def\pair#1{\left\langle #1 \right\rangle}

\def\Isocda{\Isoc^{\kern-0.5mm\dagger}}

\def\omegas{\omega^*}
\def\uomegas{\ul\omega^*}
\def\uomegav{\ul\omega^\vee}
\def\uomegasv{\ul\omega^{*\vee}}
\def\uetav{\ul\eta^\vee}
\def\ueta{\ul\eta}
\def\wuomegav{\wt{\ul\omega}^\vee}
\def\wuomegasv{\wt{\ul\omega}^{* \vee}}
\def\bsK{{\boldsymbol{K}}}

\def\theN{N}
\def\funcL{L}
\def\zerofunct{\zeta}
\def\bsKfrp{\bsK_{\!\frp}}

\def\bsD{\boldsymbol{D}}
\def\hbs{v}

\def\bN{N(\frp)}

\def\Exp#1{\exp\left[#1\right]}%
\def\pair#1{\langle#1\rangle}%
\def\markout#1{}%
\def\bnu{\boldsymbol{\nu}}

\def\bomega{\boldsymbol\omega}
\def\bomegas{\boldsymbol\omegas}

\def\bsomegav{\boldsymbol{\omega}^\vee}
\def\bsetav{\boldsymbol{\eta}^\vee}
\def\bfp{\boldsymbol{p}}

\def\pPolFunc{D^{(p)}}%
\def\pGenFunc{F^{(p)}}%

\def\AH{{\!\sA}}
\def\PolFun{D}

\def\GenFunc{F}

\begin{document}
\title[Elliptic Polylogarithm]{On the de Rham and $p$-adic realizations of the Elliptic Polylogarithm for CM elliptic curves}
\author[Bannai]{Kenichi Bannai}
\address{Department of Mathematics, Keio University, 3-14-1 Hiyoshi, Kouhoku-ku, Yokohama 223-8522, Japan}
\email{bannai@math.keio.ac.jp}
\author[Kobayashi]{Shinichi Kobayashi}
\address{Graduate school of Mathematics, Nagoya University, Furo-cho Chikusa-ku, Nagoya 464-8602, Japan}
\email{shinichi@math.nagoya-u.ac.jp}
\author[Tsuji]{Takeshi Tsuji}
\address{Department of Mathematics, University of Tokyo, Komaba, Meguro-ku,  Tokyo 153-8914, Japan}
\email{t-tsuji@ms.u-tokyo.ac.jp}
\date{September 22, 2008} 
\begin{abstract}
	In this paper, we give an explicit description of the de Rham and $p$-adic polylogarithms for elliptic curves
	using the Kronecker theta function. We prove in particular that when the elliptic curve has 
	complex multiplication and good reduction at $p$, then the specializations 
	to torsion points of the $p$-adic elliptic polylogarithm are related to $p$-adic Eisenstein-Kronecker 
	numbers, proving a $p$-adic analogue of the result of Beilinson and Levin expressing the complex 
	elliptic polylogarithm in terms of Eisenstein-Kronecker-Lerch series.   Our result is valid even if the
	elliptic curve has supersingular reduction at $p$.
\end{abstract}
\maketitle

\setcounter{tocdepth}{1}
\setcounter{section}{-1}

%
%
\section{Introduction}
%
%

%
\subsection{Introduction}
%

In the paper \cite{BL}, Beilinson and Levin constructed the elliptic polylogarithm, which is an
element in absolute Hodge or $\ell$-adic cohomology of an elliptic curve minus the identity.
This construction is a generalization to the case of elliptic curves of the construction by
Beilinson and Deligne of the polylogarithm sheaf on the projective line minus three points.
The purpose of this paper is to study the $p$-adic realization of the elliptic 
polylogarithm for an elliptic curve with complex multiplication, and to investigate its
relation to $p$-adic $L$-functions associated to the elliptic curve,
even for the case where the elliptic curve has supersingular reduction at the prime $p$.

To achieve our goal, we first describe the de Rham realization of the elliptic polylogarithm for a general
elliptic curve defined over a subfield of $\bbC$.  In particular, we explicitly describe the
connection of the elliptic polylogarithm using rational functions.  This result is of independent interest.  Similar results 
were obtained by Levin and Racinet \cite{LR} Section 5.1.3.  A different description for this connection
was also given by Besser and Solomon \cite{BS}.   

Using the de Rham realization of the elliptic polylogarithm, we then construct the 
$p$-adic realization of the elliptic polylogarithm as a filtered overconvergent $F$-isocrystal on
the elliptic curve minus the identity, when the elliptic curve has complex multiplication and good reduction at
a fixed prime $p \geq 5$.    Our main result Theorem \ref{thm: describe polylog} is an explicit description of the
$p$-adic elliptic polylogarithm sheaf in terms of overconvergent functions characterized
as the solutions of  certain $p$-adic differential equations.

Using this description, we calculate the specializations of the $p$-adic elliptic polylogarithm to torsion points of order 
prime to $p$ (more precisely, torsion points of order prime to $\frp$.  See the Overview for details), 
and prove that the specializations give the 
$p$-adic Eisenstein-Kronecker numbers, 
which are special values of the $p$-adic distribution interpolating Eisenstein-Kronecker numbers.  
This result is a generalization of the result of \cite{Ba3}, where we have dealt only with the one variable case
for an ordinary prime.   A similar result concerning the specialization in two-variables
was obtained  in \cite{BKi}, again for ordinary primes, using a very different method.  
The result of the current paper is valid even when $p$ is supersingular.

The $p$-adic Eisenstein-Kronecker numbers are related to special values of $p$-adic $L$-functions which
$p$-adically interpolate special values of Hecke $L$-functions associated to Gr\"ossencharacters of imaginary quadratic fields.
Since the $p$-adic elliptic polylogarithm is expected to be the image by the 
syntomic regulator of the motivic elliptic polylogarithm, our result may be interpreted as a $p$-adic 
analogue of Beilinson's conjecture.\medskip

In the appendix, modeling on our approach of the $p$-adic case, we calculate the real Hodge realization
of the elliptic polylogarithm by solving certain iterated differential equations as in the $p$-adic case.
The Hodge realization of the elliptic polylogarithm was first described by Beilinson-Levin \cite{BL} and Wildeshaus \cite{W}.  
We give an alternative description of the real Hodge realization in terms of multi-valued 
meromorphic functions given as the solutions of these differential equations.
Our method highlights the striking similarity between the classical and
the $p$-adic cases.

\tableofcontents

\subsection{Overview}

The detailed content of this paper is as follows.  In \S 1, 
we introduce the Kronecker theta function $\Theta(z,w)$, which is our main tool in describing the elliptic polylogarithm.
A slightly modified version of this function was previously used by Levin \cite{L} to describe the analytic aspect of the elliptic polylogarithm.  We use this function to construct rational functions $L_n$  on the elliptic curve,
which we call the connection functions.  The main result of the first section is the explicit description
of the de Rham realization of the elliptic polylogarithm in terms of the connection functions (Corollary \ref{cor: dR P}).


The main result of this paper is an explicit description of the $p$-adic elliptic polylogarithm for CM elliptic curves.
In  \cite{BK1}, we studied in detail the method to construct  from $\Theta(z,w)$
the $p$-adic distributions $p$-adically interpolating Eisenstein-Kronecker numbers
in the case of CM elliptic curves.  The relation between $p$-adic elliptic polylogarithms and $p$-adic distributions
will be proved using this result.  Let $\bsK$ be an imaginary quadratic field whose class number is one.  Let $E$ be an elliptic curve
over $\bsK$ with complex multiplication by the ring of integers $\cO_{\bsK}$ of $\bsK$, with good reduction at a prime 
$p \geq 5$.  We denote by $\psi_{E/\bsK}$ the Gr\"ossencharacter of $E$ over $\bsK$.
We fix a prime $\frp$ of $\cO_{\bsK}$ over $p$, and we let $\pi := \psi(\frp)$. Let $\Gamma$ be the
period lattice of $E$ for some invariant differential $\omega$ defined over $\cO_{\bsK}$.

In \S 2, we introduce the Eisenstein-Kronecker-Lerch series and Eisenstein-Kronecker numbers.
We fix a lattice $\Gamma$ in $\bbC$.  Let $z_0 \in \bbC \setminus \Gamma$. 
We define the Eisenstein-Kronecker numbers $e^*_{a,b}(z_0)$ for integers $a$ and $b$ by the formula
$$
	e^*_{a,b}(z_0) = \sum_{\gamma \in \Gamma\setminus \{0\}} \frac{\ol\gamma^a}{\gamma^b} \pair{\gamma, z_0},
$$
where $\pair{\gamma, z_0} := \exp((\gamma \ol z_0 - z_0 \ol\gamma)/A)$ and 
$A$ is the fundamental area of $\Gamma$ divided by $\pi$.  The above sum converges only for $b > a+2$, but
one may give it meaning for all $a$ and $b$ by analytic continuation, defining it as the
special values of Eisenstein-Kronecker-Lerch series.  
We then review the properties of the $p$-adic distribution interpolating Eisenstein-Kronecker numbers
and give the definition of $p$-adic Eisenstein-Kronecker numbers.
Let $z_0$ be a torsion point of $E(\ol\bbQ)$.  
By a theorem of Damerell, the numbers $e^*_{a,b}(z_0)/A^a$ are algebraic over $\bsK$ when $a,b \geq 0$.
We fix an embedding $i_\frp : \ol\bsK \hookrightarrow \bbC_p$ continuous for the $\frp$-adic topology on $\bsK$,
and we regard $e^*_{a,b}(z_0)/A^a$ for $a$, $b \geq 0$ as $p$-adic numbers through this embedding.
The Hodge theoretic elliptic polylogarithm is related to Eisenstein-Kronecker numbers $e^*_{a,b}(z_0)$
for $a<0$ (see Theorem \ref{thm: Beilinson and Levin}), which are complex numbers expected to be transcendental. 
Hence to consider the $p$-adic analogue, we need to define $p$-adic versions of $e^*_{a,b}(z_0)$ for $a<0$.  
In order to achieve this goal, we use $p$-adic interpolation.

Assume now that $p$ is ordinary of the form $(p) = \frp \frp^*$ in $\cO_{\bsK}$, 
and suppose that $z_0$ is non-zero of order prime to $p$.  
Then for their construction of the two-variable $p$-adic $L$-function of the CM elliptic curve (see also \cite{BK1}),
Manin-Vishik and Katz constructed a $p$-adic measure $\mu_{z_0,0}$ on 
$\bbZ_p \times \bbZ_p$ which satisfies 
$$
	\frac{1}{\Omega_\frp^{a+b}} \int_{\bbZ_p^\times \times \bbZ_p} x^a y^b d \mu_{z_0,0}(x,y) = 
	(-1)^{a+b} \left( \frac{e^*_{a,b+1}(z_0)}{A^a} - \frac{\pi^a e^*_{a,b+1}(\pi z_0)}{\ol\pi^{b+1} A^a}  \right )
$$
for any $a, b \geq 0$, where $\Omega_\frp$ is a certain $p$-adic period in $W(\ol\bbF_p)^\times$.
Using this measure, we define the $p$-adic Eisenstein-Kronecker number as follows.
\begin{definition}
	Suppose $z_0$ is a non-zero torsion point of $E(\ol\bbQ)$ of order prime to $p$.
	For any integer $a$, $b$ such that $b \geq 0$, we define the $p$-adic Eisenstein-Kronecker number
	$e_{a,b+1}^{(p)}(z_0)$ by the formula
	$$
		e^{(p)}_{a,b+1}(z_0) :=  \frac{1}{b!} \int_{\bbZ_p^\times \times \bbZ_p} x^a y^b d \mu_{z_0,0}(x,y).
	$$ 
	Note that this definition is valid even for $a < 0$. 
\end{definition}

When $p$ is supersingular, which is equivalent to the condition that $p$ remains prime in $\cO_\bsK$,
then a two-variable measure as above interpolating Eisenstein-Kronecker
numbers does not exist.  We define $e_{a,b+1}^{(p)}(z_0)$ using $p$-adic distributions, constructed
originally by Boxall \cite{Box1}\cite{Box2}, Schneider-Teitelbaum \cite{ST}, Fourquaux \cite{Fou} and Yamamoto \cite{Yam}, which
interpolate in one-variable Eisenstein-Kronecker numbers for fixed $b \geq 0$.  The latter construction is valid even when 
$p$ is ordinary.   In this case, the definition is equivalent to the one given above.   The definition of 
Eisenstein-Kronecker numbers for ordinary $p$ also extends to non-zero torsion points $z_0$ of order prime to $\frp$
(in other words, the annihilator of $z_0$ as an element in an $\cO_\bsK$-module is prime to $\frp$).
In both constructions, the crucial fact that  we use is the main result of \cite{BK1}, which states that the 
generating function for Eisenstein-Kronecker numbers is given 
by the Kronecker theta function $\Theta_{z_0,w_0}(z,w)$.  

In \S 3, we give the definition of the $p$-adic elliptic polylogarithm functions, which are overconvergent functions
on the elliptic curve minus the residue disc around the identity characterized as the solutions of a certain 
differential equation.  Then we give the relation of these functions to the $p$-adic Eisenstein-Kronecker numbers.

Finally, in \S 4, we construct and explicitly calculate the $p$-adic elliptic polylogarithm.
Let $K$ be a finite unramified extension of $\bsKfrp$.  
We denote by $\cO_K$ the ring of integers of $K$ and by $k$ its residue field.  
The rigid cohomology $H^1_\rig(E_k/K)$ of $E_k := E \otimes k$
is a Frobenius $K$-module with Hodge filtration coming from the Hodge filtration of de Rham cohomology
of $E_K := E \otimes K$ through the canonical isomorphism
$$
	H^1_\dR(E_K/K) \cong H^1_\rig(E_k/K).
$$
This cohomology group is a $K$-vector space with certain basis $\ul\omega$ and $\uomegas$.
We let $\sH$ be the filtered Frobenius module dual to $H^1_\rig(E_k/K)$, and we denote by $\ul\omega^\vee$
and $\uomegasv$ the dual basis.

Let $S(E)$ be the category of filtered overconvergent $F$-isocrystals on $E$, referred to
as the category of syntomic coefficients in our previous papers,  which plays a rough $p$-adic analogue
of the category of variations of mixed Hodge structures on $E$.   We denote by $S(\sV)$ the
same category on $\sV := \Spec\, \cO_K$, which is simply the category of filtered Frobenius modules.
The elliptic logarithm sheaf $\sLog$ is a pro-object $\sLog = \varprojlim \sLog^\theN$ in $S(E)$.
One of its main features is the splitting principle, given as follows.  

\begin{lemma}[=Lemma \ref{lem: splitting}]
	Let $z_0 \in E(K)$ be a torsion point of order prime to $\frp$.  Then we have a canonical
	 isomorphism
	$$
		i_{z_0}^* \sLog \cong \prod_{j \geq 0} \Sym^j \sH
	$$
	as filtered Frobenius modules  in $S(\sV)$.
\end{lemma}

We let $\sH^\vee$ be the dual of $\sH$, and we denote by $\sH_E$ and $\sH^\vee_E$ the pull-backs of 
$\sH$ and $\sH^\vee$ to $E$ by the structure morphism.
We let $U = E \setminus [0]$, where $[0]$ is the identity element of $E$.   The elliptic polylogarithm class
is an element $\pol$ in the rigid syntomic cohomology group
$$
	 \pol \in H^1_\syn(U, \sH_{E}^\vee \otimes \sLog(1))
$$
characterized by a certain residue condition.  The importance of this element is that it is expected to be
the image by the syntomic regulator of the motivic elliptic polylogarithm in a suitable motivic cohomology.
Our main theorem, Theorem \ref{thm: describe polylog}, is an explicit description of the $p$-adic
elliptic polylogarithm sheaf, which is an extension of $\sLog(1)$ by $\sH_{E}$ in $S(U)$ whose extension
class corresponds to $\pol$.

Let $z_0 \in E(K)$ be a non-zero torsion point of order prime to $\frp$, and let 
$i_{z_0} : \Spec\,\cO_K \rightarrow U$ be the inclusion induced by $z_0$.
By the splitting principle, the pull-back of $\pol$ to $z_0$ gives an element
$$
	i_{z_0}^*\pol \in \prod_{j \geq 0} H^1_\syn(\sV, \sH^\vee \otimes \Sym^j \sH(1)).
$$
The calculation of syntomic cohomology gives an isomorphism
\begin{equation}\label{eq: one syn}
	H^1_\syn(\sV, \sH^\vee \otimes \Sym^j \sH(1)) \cong \sH^\vee \otimes \Sym^j \sH / 
		K \ul\omega \otimes \ul\omega^{*\vee j}.
\end{equation}

Our main result is the following.

\begin{theorem}[=Theorem \ref{thm: main}]
	Suppose $p \geq 5$ is good, i.e. does not ramify in $\cO_\bsK$.  Then
	the image of $i_{z_0}^* \pol$ in $H^1_\syn(\sV, \sH^\vee \otimes \Sym^j \sH(1))$
	through the isomorphism \eqref{eq: one syn} is
	$$
		-\sum_{\substack{m+k = j \\ m \geq 1, k \geq 0}} 
		\frac{e^{(p)}_{-m,k+1}(z_0)}{\Omega_\frp^{k-m}} \ul\omega \otimes\ul\omega^{m,k}
		-
		\sum_{\substack{m+k = j \\ m \geq 0, k \geq 1}} 
		\frac{e^{(p)}_{-m-1,k}(z_0)}{\Omega_\frp^{k-m-2}} \ul\omega^* \otimes \ul\omega^{m,k},
	$$
	where $\ul\omega^{m,k}:= \ul\omega^{\vee m} \ul\omega^{*\vee k}$ and $\Omega_\frp$ is a
	certain $\frp$-adic period in $\bbC_p$.
\end{theorem}

Our main theorem shows that the $p$-adic elliptic polylogarithm specializes to give the 
$p$-adic Eisenstein-Kronecker numbers $e^{(p)}_{-m,k+1}(z_0)$ for $m \geq 0$.   This is a $p$-adic analogue of the 
result of Beilinson-Levin and Wildeshaus in the Hodge theoretic case.  Since the elliptic polylogarithm is motivic in origin, 
and since the $p$-adic Eisenstein-Kronecker numbers are values of certain $p$-adic distributions giving rise to the $p$-adic
$L$-functions in this case, one may regard this as a result of $p$-adic Beilinson conjecture type 
relating motivic elements to special values of $p$-adic $L$-functions.

%
\subsection{Acknowledgment}
%

Part of this research was conducted while the first author was visiting the \'Ecole Normale
Sup\'erieure in Paris, and the second author Institut de Math\'ematiques de Jussieu.
We would like to thank our hosts Yves Andr\'e and Pierre Colmez for hospitality.  The authors would also
like to thank John Coates and Hidekazu Furusho for pointing out the importance of the distribution relation 
which was used in our work, and Daniel Bertrand for comments concerning the definition of the
connection function $L_n(z)$.
The first and second authors were supported in part by the JSPS Postdoctoral Fellowships for Research Abroad.

%
%
\section{de Rham realization of the elliptic polylogarithm}
%
%

%
%
\subsection{Kronecker theta function}
%
%

Here, we first review the definition of the Kronecker theta function $\Theta(z,w)$.  Then we define
the connection function $L_n(z)$, which are rational functions on the elliptic curve.  This function
will later be used to describe the connection of the elliptic polylogarithm.
We fix a lattice $\Gamma \subset \bbC$, and we define $A$ to be the fundamental area of $\Gamma$ divided by
$\pi$.  In other words, if $\Gamma = \bbZ \gamma_1 \oplus \bbZ \gamma_2$ such that $\Im(\gamma_1/\gamma_2) > 0$,
then $A = (\gamma_1 \ol \gamma_2 - \gamma_2 \ol\gamma_1)/2\pi i$.

\begin{definition}
	We define $\theta(z)$ to be the reduced theta function on $E:= \bbC/\Gamma$ corresponding to
	the divisor $[0]$, normalized so that $\theta'(0) = 1$.
\end{definition}

The function $\theta(z)$ may be given explicitly as follows.
Let $\sigma(z)$ be the Weierstrass $\sigma$-function
\begin{equation*}
	\sigma(z) = z \prod_{\gamma \in \Gamma \setminus \{ 0 \}}
	\left( 1 - \frac{z}{\gamma} \right) \Exp{ \frac{z}{\gamma} + \frac{z^2}{2 \gamma^2}}, 
\end{equation*}
and let  $e^*_2:= \lim_{u \rightarrow +0} \sum_{\gamma \in \Gamma \setminus \{0\}} \gamma^{-2} |\gamma|^{- 2u}$.
Then we have
\begin{equation}\label{eq: theta sigma}
	\theta(z) = \Exp{\frac{- e_2^* z^2 }{2}} \sigma(z).
\end{equation}
The theta function $\theta(z)$ is a holomorphic function on $\bbC$ 
whose only zeroes are simple zeros at $z \in \Gamma$,
and $\theta(z)$ satisfies the transformation formula
\begin{equation}\label{eq: transform theta}
	\theta(z + \gamma) = \varepsilon(\gamma) \Exp{\frac{\ol\gamma}{A} \left(z + \frac{\gamma}{2} \right)} \theta(z)
\end{equation}
for any $\gamma \in \Gamma$, where $\varepsilon(\gamma) = -1$ if $\gamma \not\in 2 \Gamma$ and 
$\varepsilon(\gamma)=1$ otherwise.  We use $\theta(z)$ to define the Kronecker theta function.

\begin{definition}[Kronecker theta function]
	We define the Kronecker theta function $\Theta(z,w)$ to be the function
	$$
		\Theta(z,w) := \frac{\theta(z+w)}{\theta(z) \theta(w)}.
	$$
\end{definition}

This function is a reduced theta function associated to the Poincar\'e bundle of $\bbC/\Gamma$.
Note that since $\sigma(z)$ hence $\theta(z)$ is an odd function, we have $\theta''(0) = 0$.
We let $F_1(z)$ be the meromorphic function 
$$
	F_1(z) := \lim_{w \rightarrow 0} \left( \Theta(z,w) - w^{-1} \right) = 
	\theta'(z)/\theta(z),
$$
where the last equality follows from the fact that $\theta'(0) = \sigma'(0)=1$ and $\theta''(0) = 0$.
Then $F_1(z)= \zeta(z) - e_2^* z$, where $\zeta(z) := \sigma'(z)/\sigma(z)$ is the Weierstrass zeta function.
$F_1(z)$ satisfies the transformation formula $F_1(z+\gamma) = F_1(z) + \ol\gamma/A$ for any $\gamma \in \Gamma$.

\begin{definition}
	We define the function $\Xi(z,w)$ by
	$$
		\Xi(z,w) = \exp(- F_1(z) w) \Theta(z,w).
	$$	
\end{definition}
From the transformation formula \eqref{eq: transform theta}, we see that
$
	\Xi(z+ \gamma,w) = \Xi(z,w)
$
for any $\gamma \in \Gamma$.  Hence $\Xi(z,w)$ is periodic with respect to $\Gamma$ for the first variable.  
Since $\theta'(0) = 1$, we have $\lim_{w\rightarrow 0} w \Xi(z,w) = 1$.
Hence for any $z \not\in \Gamma$,
the function $\Xi(z,w) - w^{-1}$ is holomorphic with respect to $w$ in a neighborhood of $w=0$.
We define the connection function $L_n(z)$ as follows.

\begin{definition}
	We define the \textit{connection function} $L_n(z)$ to be the function in $z$ 
	whose value at a fixed $z \in \bbC \setminus \Gamma$ is defined as
	the coefficients of the Laurent expansion 
	\begin{equation}\label{eq: expand xi}
		\Xi(z,w) = \sum_{n \geq 0}  L_n(z) w^{n-1}
	\end{equation}
	of $\Xi(z,w)$ with respect to $w$ at $w=0$.
\end{definition}
The connection function is given explicitly as $L_0(z) \equiv 1$ when $n=0$ and
$$
	L_n(z) := \lim_{w \rightarrow 0} \frac{1}{(n-1)!}  \partial^{n-1}_w \left( \Xi(z,w) - w^{-1} \right)
$$
for any $n > 0$, where $\partial_w$ denotes the differential with respect to the variable $w$.   
The existence of $F_1(z) := \lim_{w \rightarrow 0} \left( \Theta(z,w) - w^{-1} \right)$ shows that  the function
$\Theta(z,w) - w^{-1}$ is holomorphic in a neighborhood of $w=0$ for a fixed $z \not\in \Gamma$. 
By exchanging  $z$ and $w$ and combining the results, we see that
the function $\Theta(z,w) - z^{-1} - w^{-1}$ is holomorphic in a neighborhood of $z=w=0$.
Hence we may consider the two-variable Taylor expansion of this function at $(z,w)=(0,0)$.
As a result, we have an expansion
$$
	\Theta(z,w) = \sum_{b \geq 0} F_b(z) w^{b-1},
$$
where $F_0(z) \equiv 1$, $F_1(z)$ is as before, and $F_b(z)$ is  holomorphic in $z$ for $b >1$.  We define the expansion 
\eqref{eq: expand xi} of $\Xi(z,w)$ to be the product of the above expansion with the standard expansion
 $\exp(- F_1(z) w) = \sum_{m \geq 0} (- F_1(z))^m w^m/m!$.
Then $L_n(z)$ is defined explicitly by the formula
$$
	L_n(z) = \sum_{b=0}^n \frac{(-F_1(z))^{n-b}}{(n-b)!} F_b(z).
$$
This shows that $L_n(z)$ is a meromorphic function for the variable $z$.

By the transformation formula for $\Xi(z,w)$, we see that $L_n(z + \gamma) = L_n(z)$ 
for any $\gamma \in \Gamma$.   Hence the connection function is an elliptic function.  Moreover, by construction,
the only poles of $\Xi(z,w)$ for the first variable are at $z \in \Gamma$.  This implies that $L_n(z)$ is holomorphic  
except for poles at $z \in \Gamma$.

We next prove that when the complex torus $\bbC/\Gamma$ has a model over a subfield $F$ of  $\bbC$, 
then the connection function corresponds to a rational function defined over $F$. Let $E$ be an elliptic curve
defined over $F$, given by the Weierstrass equation
\begin{equation}\label{eq: Weierstrass}
	E:  y^2 = 4 x^3 - g_2 x - g_3, \qquad g_2, g_3 \in F.
\end{equation}
Let $\Gamma$ be the period lattice of $E$ with respect to the invariant differential $\omega = dx/y$, and denote
by $\xi$ the complex uniformization
\begin{align}\label{eq: uniformization}
	\xi: \bbC/\Gamma&\xrightarrow\cong E(\bbC), &  z &\mapsto (\wp(z), \wp'(z)).
\end{align}

The following fact results from the recurrence relation for the 
Taylor coefficients of $\sigma(z)$ at the origin.

\begin{lemma}\label{lem: sigma}
	The Taylor expansion of $\sigma(z)$ at $z=0$ has coefficients in $F$.
\end{lemma}

\begin{proof}
	The Taylor expansion of $\sigma(z)$ at $z=0$ is given by
	$$
		\sigma(z) = \sum_{m,n\geq 0}	
		2^{n-m} a_{m,n} g_2^m g_3^n\frac{ z^{4m+6n+1}}{(4m+6n+1)!},
	$$
	where $a_{0,0}=1$, $a_{m,n} = 0$ for $m <0$ or $n<0$, and the other
	values for $a_{m,n}$ are given by the recurrence relation 
	\begin{multline*}
		a_{m,n} = 3(m+1)a_{m+1,n+1} + \frac{16}{3} (n+1) a_{m-2,n+1}\\
		- \frac{1}{3}(2m+3n-1)(4m+6n-1) a_{m-1,n}
	\end{multline*}
	(See \cite{AS} pp. 635-636 or \cite{Wei}).  Our assertion now follows from the fact that 
	$g_2$, $g_3 \in F$.
\end{proof}

Lemma \ref{lem: sigma} 
implies that the Laurent coefficients of $\zeta(z) := \sigma'(z)/\sigma(z)$,
$\wp(z) = - \zeta'(z)$ and $\wp'(z)$ at $z=0$ are also in $F$.
Note that by definition of $\Xi(z,w)$ and \eqref{eq: theta sigma}, we have 
\begin{equation}\label{eq: Xi sigma}
		\Xi(z,w) = \Exp{- \zeta(z) w} \frac{\sigma(z+w)}{\sigma(z) \sigma(w)}.
\end{equation}	
Hence the coefficients of the two-variable expansion of $\Xi(z,w)$ with respect to variables $z$ and $w$ are also in $F$.
For any rational function $f$ on $E$ defined over $F$, we denote by $f(z)$ the pullback of $f$ by $\xi$.

\begin{proposition}[Algebraicity]  
	For any integer $n \geq 0$, the connection function $L_n(z)$ is 
	obtained as a pullback by $\xi$ of a rational function $L_n$ on $E$ defined over $F$.
\end{proposition}

\begin{proof}
	The function $L_n(z)$ is holomorphic outside $[0] \in E(\bbC)$.   Since the sum of all the residues of a meromorphic 
	function is zero, the residue of $L_n(z)$ at $z=0$ must also be zero.  By the previous lemma,
	the Laurent coefficients of $L_n(z)$ at $z=0$ are in $F$.  
	The Weierstrass functions $\wp(z)$ and $\wp'(z)$ are elliptic functions corresponding to rational functions $x$ and
	$y$ of $E$.  They have poles only at $[0]$, of order $2$ and $3$, 
	and the Laurent coefficients of $\wp(z)$ and $\wp'(z)$ at the origin are in $F$.
	We may remove the negative degree of the Laurent expansion of $L_n(z)$ at $0$ by subtracting a suitable function 
	$h(z) \in F[\wp(z), \wp'(z)]$.    This implies that $L_n(z) - h(z)$ is 
	constant, since  it is a periodic function on $E(\bbC)$ without any poles.  Since the Laurent coefficients of
	both $L_n(z)$ and $h(z)$ are in $F$, the constant also must be in $F$.  This shows that $L_n(z) \in F[\wp(z), \wp'(z)]$,
	proving our assertion. 
\end{proof}

%
\subsection{Review of de Rham cohomology}
%

Next, in order to fix notations, we review some basic facts about de Rham cohomology of smooth 
algebraic varieties defined over a field $F$ with characteristic 0.  

Let $X$ be a smooth algebraic variety defined over $F$.  We denote by $\Omega^\bullet_X$ the
de Rham complex on $X$.

\begin{definition}\label{def: M(X)}
	We denote by $M(X)$ the category consisting of the pair $(\cF, \nabla)$, where $\cF$ is a locally
	free module on $X$ and $\nabla$ is an integrable connection
	$$
		\nabla : \cF \rightarrow \cF \otimes \Omega^1_X
	$$
	on $\cF$.  The category $M(X)$ is an abelian category.
\end{definition}

Suppose $\cF$ is an object in $M(X)$.  Then the connection gives rise to the complex
$$
	\Omega^\bullet(\cF) :=  \cF \xrightarrow\nabla \cF \otimes \Omega^1_X \xrightarrow\nabla
			\cF \otimes \Omega^2_X \rightarrow \cdots.
$$

\begin{definition}
	For any object  $\cF$ in $M(X)$, we define the de Rham cohomology $H^i_\dR(X, \cF)$ 
	of $X$ with coefficients in $\cF$ by
	$$
		H^i_\dR(X, \cF) :=R^i \Gamma(X, \Omega^\bullet(\cF)).
	$$
\end{definition}

We have the following proposition.

\begin{proposition}\label{proposition: extension}
	We have a canonical isomorphism
	$$
		\Ext^1_{M(X)}( \cO_X, \cF) \xrightarrow\cong H^1_\dR(X, \cF).
	$$
\end{proposition}

\begin{proof} The canonical homomorphism is given as follows.  Suppose we have an extension
	$$
		0 \rightarrow \cF \rightarrow \cE \rightarrow \cO_X \rightarrow 0
	$$
	in $M(X)$.  Then the boundary morphism of the long exact sequence associated to 
	the above exact sequence defines a map
	$$
		H^0_\dR(X, \cO_X) \rightarrow H^1_\dR(X, \cF).
	$$
	We define the class $[\cE]$ of $\cE$ in  $H^1_\dR(X, \cF)$ to be the image of $1 \in F \subset H^0_\dR(X, \cO_X)$ 
	by the above map.  The inverse homomorphism is defined as follows.
	Suppose $\frU= \{ U_i \}_{i \in I}$ is an affine open covering of $X$.
	Then the de Rham cohomology of $X$ with
	coefficients in $\cF$ may be calculated using the \v Cech resolution for this covering.  In particular, any cohomology class
	$[\xi]$ in $H^1_\dR(X, \cF)$ may be represented by a cocycle 
	$$
		(\xi_i, u_{ij}) \in	\prod_{i\in I} \Gamma(U_i, \cF \otimes \Omega^1_X)
		 \oplus \prod_{i,j \in I} \Gamma(U_i \cap U_j, \cF),
	$$
	satisfying  $d \xi_i = 0$, $d u_{ij} =\xi_j - \xi_i$ and  $u_{ij}  + u_{jk} = u_{ik}$ for any $i,j,k \in I$.  The extension
	$\cE$  whose class in $\Ext^1_{M(X)}(\cO_X, \cF)$ corresponds to 
	$[\xi] \in H^1_\dR(X, \cF)$ through Proposition \ref{proposition: extension} is constructed
	as follows.  We define $\cE_i$ to be the coherent $\cO_{U_i}$-module 
	$\cE_i := \cO_{U_i} \ul e_i \bigoplus \cF|_{U_i}$, with connection $\nabla(\ul e_i) := \xi_i$.  We define $\cE$ to be
	the coherent module with connection on $X$ obtained by pasting together $\cE_i$ on $U_i \cap U_j$
	through the isomorphism
	\begin{align*}
		\cE_i|_{U_i \cap U_j} &\cong \cE_j|_{U_i \cap U_j},& \ul e_i &= \ul e_j - u_{ij}.
	\end{align*}
	The compatibility of the pasting isomorphism with the connection follows from the
	cocycle condition.
\end{proof}

We next review the localization sequence for de Rham cohomology for smooth curves.  
We first review the properties of logarithmic
de Rham cohomology.  Let $X$ be a smooth algebraic variety over $F$ and $D \hookrightarrow X$
be a normal crossing divisor of $X$ over $F$.   We denote by $\Omega_X^1(\log D)$ the sheaf of differentials on 
$X$ with logarithmic poles along $D$.  Then for any $\cF \in M(X)$, we may naturally 
define a logarithmic connection
$$
	\nabla_{\log} : \cF \rightarrow \cF \otimes \Omega^1_X(\log D)
$$
through the natural inclusion $\cF \otimes \Omega^1_X \hookrightarrow \cF \otimes \Omega^1_X(\log D)$.  We let
$$
	\Omega^\bullet_{\log}(\cF) := \cF \xrightarrow{\nabla_{\log}} \cF \otimes \Omega^1_X(\log D) 
	\xrightarrow{\nabla_{\log}} \cF \otimes \Omega^2_X(\log D) \rightarrow\cdots
$$
and we define the logarithmic de Rham cohomology of $X$ with logarithmic poles along $D$
and coefficients in $\cF$ by
$$
	H^i_{\log \dR}(X, \cF) := R^i\Gamma(X, \Omega^\bullet_{\log}(\cF) ).
$$
Then for $U := (X \setminus D) \overset{j}{\hookrightarrow} X$, the natural inclusion 
$\Omega^\bullet_{\log}(\cF)  \rightarrow j_* \Omega^\bullet(\cF)|_U$ induces a canonical isomorphism
\begin{equation}\label{eq: log and open}
	H^i_{\log \dR}(X, \cF) \xrightarrow\cong H^i_{\dR}(U, \cF).
\end{equation}

\begin{proposition}[Localization sequence]\label{prop: loc seq}
	 Suppose $X$ is a smooth curve defined over $F$, and suppose $i: D \hookrightarrow X$ is a smooth 
	 divisor of $X$ defined over $F$.  Then for $\cF$ in $M(X)$ and $U := X \setminus D$, 
	 the long exact sequence associated to the exact sequence of
	 coherent sheaves
	 $$
	 	0 \rightarrow \Omega_X^\bullet(\cF) \rightarrow \Omega^\bullet_{\log}(\cF) \rightarrow i_* i^* \cF[-1] \rightarrow 0
	 $$
	and the isomorphism \eqref{eq: log and open} gives an isomorphism 
	$$H^0_\dR(X, \cF) \xrightarrow\cong H^0_\dR(U, \cF)$$ and an exact sequence
	\begin{align*}
		0 &\rightarrow H^1_\dR(X, \cF) \rightarrow H^1_\dR(U, \cF) \xrightarrow\res H^0_\dR(D, i^* \cF) \\
		   &\rightarrow  H^2_\dR(X, \cF) \rightarrow H^2_\dR(U, \cF) \rightarrow 0.
	\end{align*}
\end{proposition}

The residue map $\res$ in the above proposition is given explicitly as follows.  The map
\begin{equation}\label{eq: res one}
	\res: H^1_{\log \dR}(X, \cF) \rightarrow H^0_\dR(D, i^*\cF)
\end{equation}
is given as follows.
Let $\frU := \{ U_i \}_{i \in I}$ be an affine open covering of $X$ satisfying the property that for any $i \in I$ such that 
$D_i := U_i \cap D \not= \emptyset$, there exists a local parameter $t_i$ on $U_i$ whose zero locus $t_i=0$ 
defines the divisor $D_i$ in $U_i$.  Again, $ H^1_{\log \dR}(X, \cF)$ may 
be calculated using the \v Cech resolution for this covering. Hence any cohomology class $[\xi]$ in 
$H^1_{\log \dR}(X, \cF)$ may be represented by a cocycle
$$
	(\xi_i, u_{ij}) \in	\prod_{i\in I} \Gamma(U_i, \cF \otimes \Omega^1_X(\log D)) \oplus 
	\prod_{i,j \in I} \Gamma(U_i \cap U_j, \cF),
$$
satisfying  $d u_{ij} =\xi_j - \xi_i$ and $u_{ij}  + u_{jk} = u_{ik}$ for any $i,j,k \in I$.  We let 
$\res(\xi_i) \in \Gamma(D_i, i^* \cF)$ be the residue of $\xi_i$ at $D_i$.  In other words,
we let $\res(\xi_i) = 0$ if $D_i = \emptyset$, and for any $i \in I$ such that $D_i \not= \emptyset$, we let
$\res(\xi_i) := \alpha_i \pmod{t_i}$, where $\alpha_i \in \Gamma(U_i, \cF)$ is the element such that
$
	\xi_i :=  \alpha_i \otimes d \log t_i.
$
We define the map \eqref{eq: res one} by associating to the element
$[\xi]$ the element of $H^0_\dR(D, i^* \cF)$ represented by the cocycle
$$
	(\res(\xi_i))_i \in \prod_{i \in I}\Gamma(D_i, i^* \cF).
$$
Then the residue map of Proposition \ref{prop: loc seq} is the composition
$$
	\res: H^1_\dR(U, \cF) \xleftarrow\cong 
	 H^1_{\log \dR}(X, \cF)
	\xrightarrow\res H^0_\dR(D, i^* \cF).
$$

%
\subsection{The Logarithm sheaf}
%

We return to the case of an elliptic curve $E$ defined over a field $F \subset \bbC$ 
as in \eqref{eq: Weierstrass}.  Let
$
	 H^1_{\dR}(E)
$
be the first de Rham cohomology of $E$.  This cohomology may be calculated as
\begin{equation}\label{eq: de Rham second}
	H^1_{\dR}(E) \xrightarrow\cong R^1 \Gamma \left(E, \cO_E([0]) \xrightarrow d  \Omega^1_E(2[0]) \right) 
	\xleftarrow\cong \Gamma(E, \Omega^1_E(2[0]))
\end{equation}
using differentials of the second kind. 
We denote by $\ul\omega$ and $\ul\eta$ the classes in $H^1_{\dR}(E)$
corresponding to  the differentials $\omega: = dx/y$ and $\eta: = x dx/y$ in $\Gamma(E, \Omega^1_E(2[0]))$, which form
a basis $\{ \ul\omega, \ul\eta\} $ of $H^1_\dR(E)$.  The differentials correspond to
$\omega =dz$ and $\eta = \wp(z) dz = - d \zeta(z)$ through the complex uniformization \eqref{eq: uniformization}.
We denote by
$$
	\cH := H^1_{\dR}(E)^\vee
$$
the dual of $H^1_\dR(E)$, and we let  $\{\uomegav, \uetav \}$ be the
dual basis of $\{ \ul\omega, \ul\eta \}$.   For any smooth scheme $X$ over $F$, 
we denote by $\cH_X$ the coherent module $\cH \otimes \cO_X$ on $X$ with connection such that
$
	\nabla(\uomegav) = \nabla(\uetav) = 0.
$
Since the connection on $\cH_E$ is trivial, we have a canonical isomorphism
\begin{equation}\label{equation: leray}
	H^1_{\dR}(E, \cH_E) = H^1_{\dR}(E) \otimes \cH = \Hom(\cH, \cH).
\end{equation}

\begin{definition}\label{def: log dR}
	We define the first logarithm sheaf $\cLog^{(1)}$ to be any extension of $\cH_E$ by $\cO_E$ in $M(E)$,
	 whose extension class in 
	$$
		\Ext^1_{M(E)}(\cO_E, \cH_E) \cong H^1_\dR(E, \cH_E)
	$$
	is mapped by \eqref{equation: leray} to the identity.  We define the $N$-th logarithm sheaf $\cLog^N$ to 
	be the $N$-th symmetric tensor product of $\cLog^{(1)}$.
\end{definition}
There exists a natural projection $\cLog^{\theN+1} \rightarrow \cLog^\theN$ defined as the composite
$$
	\Sym^{\theN+1} \cLog^{(1)} \rightarrow \Sym^{\theN+1} \left( \cLog^{(1)} \bigoplus \cO_E \right) \rightarrow
	\Sym^\theN \cLog^{(1)},
$$
where the first map is the sum of the identity and the projection $\cLog^{(1)} \rightarrow \cO_E$, and the second
map is the projection in the symmetric algebra of a direct sum.

Let $i \colon  [0] \hookrightarrow E$ be the natural inclusion of the identity $[0]$ in $E$.
Since a connection on a point is zero, we have a splitting
\begin{equation}\label{eq: choice of splitting}
	\varphi \colon i^*_{[0]}\cLog^{(1)} \cong F \bigoplus \cH
\end{equation}
on $M(\Spec\,F)$.  A choice of a splitting $\varphi$ as above induces a splitting
$
	\varphi\colon i^*_{[0]} \cLog^\theN \cong \prod_{j=0}^\theN \Sym^j \cH
$ 
on the $N$-th symmetric tensor product, which is compatible with the projection $\cLog^{\theN+1} \rightarrow \cLog^\theN$.

\begin{remark}\label{rem: not canonical}
	Unlike the case for the Hodge or $p$-adic realizations, the first logarithm sheaf $\cLog^{(1)}$ 
	has non-trivial automorphisms as extensions of $\cH_E$ by $\cO_E$.  Hence there exists
	multiple isomorphisms between various choices of $\cLog^{(1)}$.  However, If we choose a
	pair $(\cLog^{(1)},  \varphi)$ consisting of a first logarithm sheaf and a splitting $\varphi$ as in \eqref{eq: choice of splitting},
	then the pair $(\cLog^{(1)},  \varphi)$ is unique up to unique isomorphism preserving the splitting.
\end{remark}


We now give an explicit construction of $\cLog^{(1)}$.  We first construct a cocycle whose 
cohomology class in $H^1_\dR(E, \cH_E)$ is mapped to the identity 
$ \ul\omega^\vee \otimes \ul\omega  + \uetav \otimes \ul\eta$ by \eqref{equation: leray}.  
Note that $\ul\omega \in H^1_\dR(E)$ corresponds 
to a global holomorphic differential $\omega \in \Gamma(E, \Omega^1_E)$, where as $\ul\eta \in H^1_\dR(E)$
does not.  We take an affine open covering $\frU = \{ U_i \}_{i \in I}$ of $E$.  Then there exists an element
\begin{equation}\label{eq: second kind}
	(  \eta_i, u_{ij})  \in  \prod_{i \in I} \Gamma(U_i, \Omega^1_E) \oplus 
		\prod_{i, j \in I} \Gamma(U_i \cap U_j, \cO_E)
\end{equation}
satisfying the cocycle conditions
\begin{align*}
	d u_{ij} &= \eta_j - \eta_i,    &  u_{ij} + u_{jk} &= u_{ik}  
\end{align*}
for any $i$, $j$, $k \in I$, which represents the class of $\ueta$ in $H^1_\dR(E)$.  
Then, if we set $\boldsymbol\nu_{i}  = \uomegav \otimes \omega + \uetav \otimes \eta_i $
and $\boldsymbol u_{ij} = u_{ij} \uetav$, the pair
$$
	(  \boldsymbol\nu_i, \boldsymbol u_{ij})  \in  \prod_{i \in I} \Gamma(U_i, \cH_E \otimes \Omega^1_E) \oplus 
		\prod_{i, j \in I} \Gamma(U_i \cap U_j, \cH_E)
$$
is a cocycle which represents the cohomology class in $H^1_\dR(E, \cH_E)$ that maps to the identity
by \eqref{equation: leray}.  Hence the relation between cocycles and extensions in 
Proposition \ref{proposition: extension} gives the following proposition.


\begin{proposition}\label{pro: log dR}
	On each open affine neighborhood $U_i$, we let
	$$
		\cLog^{(1)}_i := \cO_{U_i} \ul e_i \bigoplus \cH_{U_i}
	$$
	with connection horizontal on $\cH$ and
	$\nabla(\ul e_i) = \boldsymbol\nu_i \in \Gamma(U_i, \cH_E \otimes \Omega^1_E)$.
	Then $\cLog^{(1)}$ is obtained by pasting together $\cLog^{(1)}_i$ 
	through the isomorphism 
	\begin{align}\label{eq: Log one}
		 \cLog^{(1)}_i|_{U_i \cap U_j}  &\cong \cLog^{(1)}_j|_{U_i \cap U_j}, &
		  \ul e_i &= \ul e_j  - \boldsymbol u_{ij} = \ul e_j - u_{ij} \uetav
	\end{align}
	on $U_i \cap U_j$.
\end{proposition}

The $N$-th logarithm sheaf is the $N$-th symmetric tensor product of $\cLog^{(1)}$.  
If we let
$
	\ul\omega^{m,n}_i: = \ul e_i^{a}{\uomegav}^m {\uetav}^n/a!
$  
for $a = N - m - n$, then $\cLog^{\theN}$ is given as follows.

\begin{corollary}\label{cor: Log N}
	On each open affine neighborhood $U_i$, we let
	$$
		\cLog^{\theN}_i := \bigoplus_{0 \leq m+n \leq \theN} \cO_{U_i} \ul\omega_i^{m,n},
	$$
	with connection $\nabla_{\kern-1mm\cL} = d + \boldsymbol\nu_i$.  In other words, the connection is
	defined to satisfy $\nabla_{\kern-1mm\cL}(\ul\omega_i^{m,n}) 
	= \ul\omega_i^{m+1,n} \otimes \omega + \ul\omega_i^{m,n+1} \otimes \eta_i$.
	Then the $N$-th logarithm sheaf $\cLog^{\theN}$ is given by pasting together 
	$\cLog^{\theN}_i$ on $U_i \cap U_j$ through the isomorphism  
	\begin{align}\label{eq: Log N}
	 	\cLog^{\theN}_i|_{U_i \cap U_j}  & \cong \cLog^{\theN}_j|_{U_i \cap U_j}, &
		\ul\omega_i^{m,n} &= \sum_{k=n}^{\theN-m} \frac{(-u_{ij})^{k-n}}{(k-n)!} \ul\omega_j^{m,k}.
	\end{align}
\end{corollary}

\begin{proof}
	This follows by calculating $\ul\omega_i^{m,n}$ in terms of $\ul\omega_j^{m,k}$.
\end{proof}

The projection $\cLog^{\theN+1} \rightarrow \cLog^\theN$ on each $U_i$ is defined by mapping $\ul\omega_i^{m,n}$
to $\ul\omega_i^{m,n}$ if $m+n \leq \theN$, and to zero if $m+n = \theN+1$.

Using differentials of the second kind, one may describe the restriction of $\cLog^{(1)}$ to 
$U = E \setminus [0]$.  By definition, the differential of the second kind $\eta$ represents the class 
$\ueta$ in $H^1_\dR(E)$
through the isomorphism \eqref{eq: de Rham second}.
Since $(\eta_i, u_{ij})$ also represents the same class in $H^1_\dR(E)$,
this implies that there exists an element
\begin{equation}\label{eq: coboundary}
	(u_i)   \in \prod_{i \in I} \Gamma( U_i , \cO_{E}([0]))
\end{equation}
satisfying the coboundary condition $\eta =   \eta_i-d u_{i}$ and $ u_{ij} =  u_j - u_i$ for any $i$, $j \in I$.  
Denote by $\boldsymbol\nu$ the differential
$$
	\boldsymbol\nu = \uomegav \otimes \omega + \uetav \otimes \eta 
	\in \Gamma(E, \cH_E \otimes \Omega^1(2[0])).
$$

\begin{proposition}\label{pro: log dR on U}
	The restriction of $\cLog^{(1)}$ to $U = E \setminus [0]$ is given by
	the free $\cO_U$-module
	$$
		\cL: = \cO_U  \ul e \bigoplus  \cH_U,
	$$
	with connection horizontal on $\cH$ and satisfying $\nabla(\ul e) =  \boldsymbol\nu$.  
\end{proposition}

\begin{proof}
	For each $i \in I$, we have an isomorphism
	$
		\cL |_{U \cap U_i}    \cong  \cLog^{(1)}_i|_{U \cap U_i}
	$
	on $U \cap U_i$ which is the identity on $\cH$ and $\ul e = \ul e_i - u_i \uetav$.
	The coboundary condition for $u_i$ implies that it is compatible with the connection.
\end{proof}
If we let
$
	\ul\omega^{m,n}: =\ul e^{a}{\uomegav}^m {\ueta}^{\vee n}/a!
$  
for $a = N - m - n$, then the restriction of $\cLog^{\theN}$ to $U$ is given as follows.

\begin{corollary}\label{cor: log isom}
	The restriction of $\cLog^\theN$ to $U$ is given by the free $\cO_U$-module
	$$
		\cL^\theN: = \bigoplus_{0 \leq m+n \leq \theN} \cO_U \ul\omega^{m,n}
	$$
	with connection $\nabla_{\kern-1mm\cL} = d + \boldsymbol\nu$.   In other words, the connection
	satisfies 
	$$\nabla_{\kern-1mm\cL}( \ul\omega^{m,n}) =  \ul\omega^{m+1,n} \otimes \omega +  \ul\omega^{m,n+1} \otimes \eta.
	$$
	For any $i \in I$, the isomorphism
	$
		 \cL^\theN|_{U \cap U_i} \cong \cLog^{\theN}_i|_{U \cap U_i} 
	$  
	 is given by
	\begin{equation}\label{eq: comparison epsilon}
		\ul \omega^{m,n} = \sum_{k=n}^{\theN-m} \frac{(-u_i)^{k-n}}{(k-n)!} \ul\omega_i^{m,k}
	\end{equation}
	on $U \cap U_i$.
\end{corollary}

\begin{proof}
	The statement follows from the fact that $\cL^\theN$ and $\cLog^\theN$ are both $N$-th symmetric
	products of $\cL$ and $\cLog^{(1)}$ and explicit calculation of $\ul\omega^{m,n}$ in terms of $\ul\omega_i^{m,k}$.
\end{proof}

To finish this section, we give a choice of a splitting $\varphi \colon i_{[0]}^* \cLog^{(1)} \cong F \bigoplus \cH$
as in \eqref{eq: choice of splitting}.
We fix an affine open covering $\frU = \{ U_i \}_{i \in I}$ of $E$, and
we choose $(\eta_i, u_{ij})$ and  $(u_i)$ as in \eqref{eq: second kind} and  \eqref{eq: coboundary}.   
Consider any $i \in I$ such that $[0] \in U_i$.  Since $\eta_i = \eta + d u_i$ is a meromorphic differential form without poles on $U_i$, 
the function $\zerofunct_i(z): = \zeta(z) - u_i(z)$ is a meromorphic function on $\bbC$, holomorphic on the inverse 
image of $U_i$ by the uniformization $\bbC \twoheadrightarrow \bbC/\Gamma \cong E(\bbC)$,
whose value at $0$ is an element in $F$.
If we let $\wt u_i := u_i + \zerofunct_i(0)$ if $[0] \in U_i$ and $\wt u_i = u_i$ otherwise, and if we replace 
$u_{ij}$ by $\wt u_{ij} = \wt u_j- \wt u_i$, then
$(\eta_i, \wt u_{ij})$ and $(\wt u_i)$ are also elements defined over $F$ satisfying \eqref{eq: second kind} and \eqref{eq: coboundary}.
Hence by replacing $u_i$ by $\wt u_i$ and $u_{ij}$ by $\wt u_{ij}$, we may assume that 
$\zerofunct_i(0) = ( \zeta - u_i)(0) = 0$ if $[0] \in U_i$.

Consider any open set $U_i$, $U_j \in \frU$ such that $[0] \in U_i \cap U_j$.
By our choice of the coboundary $(u_i)$ and cocycle $(\eta_i, u_{ij})$,   we have
$$
	u_{ij}(0) = (u_j - u_i)(0) = (\zerofunct_i - \zerofunct_j)(0) = 0.
$$
Hence the pullback by $i_{[0]}^*$ of the pasting isomorphism  \eqref{eq: Log one} is given by 
\begin{align}\label{eq: basis trivial}
	 i_{[0]}^* \cLog^{(1)}_i &\xrightarrow\cong i_{[0]}^* \cLog^{(1)}_j, & 
	\ul e_i&= \ul e_j.
\end{align}
For any $i \in I$ such that $[0] \in U_i$, we define the splitting
\begin{equation}\label{eq: split log one}
	\varphi: i_{[0]}^* \cLog^{(1)} \cong F \bigoplus \cH
\end{equation}
to be the isomorphism which is the identity on  $\cH \subset  i_{[0]}^* \cLog^{(1)}$ and
mapping $\ul e_i$ to the identity element $1 \in F$.  By \eqref{eq: basis trivial}, this splitting
is independent of the choice of $i \in I$.  

\begin{definition}\label{def: basis}
	We fix $(\eta_i, u_{ij})$ and $(u_i)$ as above.
	We let  $\{ \ul\omega^{m,n}_{[0]} \}$ be the basis of $i_{[0]}^* \cLog^\theN$ obtained as the
	restriction of the basis $\{ \ul\omega_i^{m,n} \}$ of $\cLog^\theN_i$ for some $i \in I$ such that $[0] \in U_i$.  
\end{definition}

By \eqref{eq: basis trivial}, the basis $\{ \ul\omega_{[0]}^{m,n} \}$ is independent of the choice of $i$.   We have a splitting
\begin{equation}\label{eq: our choice}
	i_{[0]}^* \cLog^\theN 
	\cong \prod_{k=0}^\theN \Sym^k \cH
\end{equation}
induced from the splitting \eqref{eq: split log one}, given by mapping $\ul\omega_{[0]}^{m,n}$ to 
$\ul\omega^{\vee m} \ul\omega^{*\vee n}$. When considering the Hodge theoretic (resp. the $p$-adic) case, 
we will give additional structures in such a way that \eqref{eq: our choice} is the unique splitting which is compatible with the 
respective structures (Remark \ref{rem: C split} and Proposition \ref{pro: p split}).
In other words, the above isomorphism will be isomorphisms of mixed Hodge structures (resp.  
filtered Frobenius modules).

%
\subsection{The Polylogarithm sheaf}
%

Here, we define and explicitly describe the elliptic polylogarithm sheaf, originally studied by
Beilinson and Levin \cite{BL}.   Let the notations be as before.  In particular, let $\cLog^\theN$ be the logarithm
sheaf explicitly described in Corollary \ref{cor: Log N} and let $\varphi\colon i_{[0]}^* \cLog^{\theN} \cong \prod_{j =0}^\theN \Sym^j \cH$
be the splitting given in \eqref{eq: our choice}.  The following calculation is important in defining the polylogarithm sheaf.

\begin{lemma}\label{lem: log coh} 
	The projections
	\begin{align*}
		H^0_\dR(E, \cLog^{\theN+1}) &\rightarrow H^0_\dR(E, \cLog^\theN) \\
		H^1_\dR(E, \cLog^{\theN+1}) &\rightarrow H^1_\dR(E, \cLog^\theN) 
	\end{align*}
	are zero maps, and the projection gives an isomorphism
	$H^2_\dR(E, \cLog^{\theN+1}) \xrightarrow\cong H^2_\dR(E, \cLog^\theN) \cong F$.
\end{lemma}

\begin{proof}
	This statement is proved in \cite{BL}.  See also \cite{HK} Lemma A.1.4 or \cite{Ba3} Lemma 3.4.
\end{proof}

Let $D = [0]$ and $U = E \setminus [0]$.  Then by Proposition \ref{prop: loc seq}, we have a long exact sequence
\begin{align*}
	0 &\rightarrow H^1_\dR(E, \cLog^\theN) \rightarrow H^1_\dR(U, \cLog^\theN) \xrightarrow\res
	H^0_\dR(D,  i_{[0]}^* \cLog^\theN) \\
	& \rightarrow H^2_\dR(E, \cLog^\theN) \rightarrow 0.
\end{align*}
By Lemma \ref{lem: log coh}, the residue map induces an isomorphism
\begin{equation*}
	 \varprojlim_\theN H^1_\dR(U, \cLog^\theN)\\ \xrightarrow\cong  \varprojlim_\theN 
	\Ker \left(H^0_\dR(D,  i_{[0]}^* \cLog^\theN)  \rightarrow H^2_\dR(E, \cLog^\theN) \right).
\end{equation*}
The calculation of $H^0_\dR$ using the splitting \eqref{eq: our choice}
and the fact that we have $H^2_\dR(E, \cLog^\theN) = F$ give the isomorphism
\begin{equation}\label{eq: res isom}
	\res\colon \varprojlim_\theN H^1_\dR(U, \cLog^\theN) \cong \prod_{k \geq 1} \Sym^k \cH.
\end{equation}
Note that this isomorphism depends on the choice of \eqref{eq: our choice}.
We may now define the polylogarithm class and the polylogarithm sheaf for de Rham cohomology.

\begin{definition}
	We define the polylogarithm class (for the splitting \eqref{eq: our choice})
	to be a system of classes $\pol^\theN \in  H^1_\dR(U, \cH_U^\vee \otimes \cLog^\theN)$
	such that
	$$
		\pol := \varprojlim_\theN \pol^\theN \in 	\varprojlim_\theN  H^1_\dR(U, \cH_U^\vee \otimes \cLog^\theN)
	$$
	maps through \eqref{eq: res isom} to the identity in $\cH^\vee \otimes \cH \subset \cH^\vee \otimes 
	\prod_{k \geq 1} \Sym^k \cH$.
\end{definition}

\begin{definition}\label{def: de Rham polylogarithm sheaf}
	We define the elliptic polylogarithm sheaf on $U$ to be a system of sheaves $\cP^\theN$ in $M(U)$ given 
	as an extension
	$$
		0 \rightarrow \cLog^\theN \rightarrow \cP^\theN \rightarrow \cH_U \rightarrow 0
	$$
	whose extension class in 
	$$
		\Ext^1_{M(U)}(\cH_U, \cLog^\theN) \cong H^1_\dR(U, \cH^\vee_U \otimes \cLog^\theN)
	$$ 
	is $\pol^\theN$.
\end{definition}

The rest of this section is devoted to explicitly describing the extension class $\pol^\theN$ and the
polylogarithm sheaf $\cP^\theN$.  We first construct certain classes in $H^1_\dR(U,  \cLog^N)$.
To ease the notation, we write $L^0_n(z) := (-1)^n L_n(z)$.

\begin{definition}
	We define $\bsomegav$ and $\bsetav$  to be the sections
	\begin{align*}
		\bsomegav &=  \ul\omega^{0,0} \otimes \eta - \sum_{n=1}^{\theN}  \funcL^0_{n} \ul\omega^{1,n-1} 
		\otimes  \omega, &
		\bsetav &= -\sum_{n=0}^\theN  \funcL_{ n}^{0} \ul\omega^{0,n} \otimes  \omega	
	\end{align*}
	in  $\Gamma(U, \cLog^\theN \otimes \Omega^1_{U})$.
\end{definition}

Since $U$ is affine, the de Rham cohomology of $U$ may be calculated as the cohomology
of the complex
$$
	\Gamma(U,  \cLog^\theN) \xrightarrow{\nabla} \Gamma(U,  \cLog^{\theN} \otimes \Omega^1_U).
$$
Hence $\bsomegav$ and $\bsetav$ define classes
$[\bsomegav]$ and $[\bsetav]$ in 
$
	 H^1_{\dR}(U,  \cLog^\theN).
$
The following theorem is crucial in describing the elliptic polylogarithm.

\begin{theorem}\label{thm: de Rham}
	We let $\bfp^\theN$ be the cohomology class
	$$
		\bfp^\theN := \ul\omega \otimes [\bsomegav]  + \ueta \otimes [\bsetav] 
			 \in  \cH^\vee \otimes H^1_{\dR}(U,  \cLog^\theN).
	$$
	Then the image of $\bfp^\theN \in H^1_{\dR}(U, \cH_U^\vee \otimes \cLog^\theN)$ by the residue map
	\begin{equation}\label{eq: residue}
		\res: H^1_\dR(U, \cH_U^\vee \otimes \cLog^\theN) \rightarrow \cH^\vee \otimes \prod_{k \geq 1} \Sym^k \cH
	\end{equation}
	is the identity 
	$ \ul\omega \otimes \ul\omega^\vee  + \ueta\otimes \ul\eta^\vee$
	 in $\cH^\vee \otimes \cH \subset \cH^\vee \otimes \prod_{k \geq 1} \Sym^k \cH$.  In particular, we have
	 $$
	 	\pol := \varprojlim_\theN \pol^\theN = \varprojlim_\theN \bfp^\theN.
	 $$
\end{theorem}

The proof of this theorem will be given at the end of this section. 
The correspondence between cohomology classes in $H^1_{\dR}(U, \cH_U^\vee \otimes \cLog^\theN)$ and
extensions $\Ext^1_{M(U)}(\cH_U,  \cLog^\theN)$ gives the following corollary.

\begin{corollary}\label{cor: dR P}
	The polylogarithm sheaf $\cP^\theN$ may be constructed as follows.
	As a coherent $\cO_{U}$-module, it is given as the sum
	$$
		\cP^\theN := \cH_U \bigoplus   \cLog^\theN_{U},
	$$
	and the connection $\nabla_{\kern-1mm\cP}$ is given by $\nabla_{\kern-1mm\cL}$ on $\cLog^\theN$ and
	\begin{align*}
		\nabla_{\kern-1mm\cP}(\uomegav) &=  \bsomegav,  & \nabla_{\kern-1mm\cP}(\uetav) &= \bsetav.
	\end{align*}
\end{corollary}

We now prepare some results necessary for the proof of the theorem.
In order to calculate the residue \eqref{eq: residue}, we must express the classes  $[\bsomegav]$ and $[\bsetav]$
in terms of cocycles in logarithmic de Rham cohomology
$$
	H^1_{\log\dR}(E,  \cLog^\theN) := R^1 \Gamma(E, \Omega^\bullet_{\log}(\cLog^\theN)).
$$
We take an open affine covering $\frU = \{ U_i \}$ of $E$, and we fix a cocycle
$(\eta_i, u_{ij})$ and a coboundary $(u_i)$ as in Definition \ref{def: basis}.
We let
$$
	\Xi_i(z,w) := \exp(u_i(z) w) \Xi(z,w) = \Exp{- \zeta_i(z) w} \frac{\sigma(z+w)}{\sigma(z) \sigma(w)}
$$
for $\zeta_i(z) = \zeta(z) - u_i(z)$, where the last equality follows from \eqref{eq: Xi sigma}.

\begin{definition}
	For any integer $k \geq 0$ and $i \in I$, we define $\funcL_{ k, i}(z)$ 
	to be the function in $z$ given by
	$$
		\Xi_i(z,w) = \sum_{k \geq 0} L_{k,i}(z) w^{k-1},
	$$
	and we again let $L^0_{k,i}(z) := (-1)^k L_{k,i}(z)$.
\end{definition}

By definition, $L_{k,i}(z)$ may be expressed in terms of $L_n(z)$ and $u_i(z)$, hence comes from a rational function
$L_{k,i}$ on $E$ defined over $F$.  The importance of this function is the following property.

\begin{proposition}\label{prop: L good}
	The functions $\funcL_{ k, i}$ 
	are holomorphic on $U_i \setminus [0]$.   It has a simple pole of residue one at $[0]$ 
	if $k=1$ and is holomorphic on $U_i$ if $k\not= 1$.
\end{proposition}

\begin{proof}
	By definition, we have $L_{0,i} = L_{0} \equiv 1$, so the statement is true for $k=0$.  Similarly, we have 
	$L_{1,i} = u_i + L_1 = u_i$.   By definition, $u_i$ is holomorphic outside $[0] \in U_i$.
	When $[0] \in U_i$, our choice of $u_i$ implies that the differential form 
	$d \zeta_i(z) = d \zeta(z) - d u_i$ is holomorphic at $z=0$. Hence the residue of $\zeta_i(z) = \zeta(z) - u_i(z)$ at 
	$z=0$ is zero.   Since the residue of  $\zeta(z)$ at $z=0$ is one, the residue of $u_i(z)$ at $z=0$ is also one.
	This gives the assertion for $k=1$.  For the general case, note that for a fixed $w \not\in\Gamma$, the function
	$\sigma(z+w)/\sigma(z)\sigma(w) - z^{-1}$ is holomorphic in a neighborhood of $z=0$.
	Note also that we have chosen $\zeta_i(z)$ so that $\zeta_i(0)=0$.
	Hence the function appearing in the coefficient of $w^{k-1}$ in the expansion of $\Xi_i(z,w)$ for $k>1$ 
	must be holomorphic at $z=0$.  This implies that $L_{k,i}$ is holomorphic on $U_i$ when $k>1$.
	Our assertion is now proved.
\end{proof}

This proposition gives the following corollary.

\begin{corollary}
	Let 
	\begin{align*}
		\bsomegav_{i} &:= \ul\omega_i^{0,0} \otimes \eta_i -
			 \sum_{k=1}^{\theN} \funcL^0_{ k, i} \ul\omega_i^{1, k-1} \otimes  \omega, &
			\bsetav_i &:= - \sum_{k=0}^\theN  \funcL^0_{ k, i} \ul\omega_i^{0, k} \otimes  \omega.
	\end{align*}
	Then we have
	 $\bsomegav_{i}$, $\bsetav_i
	  \in \prod_{i \in I}\Gamma\left(U_i, \cLog^\theN_i \otimes \Omega^1_{E}(\log [0])\right)$.
\end{corollary}
	
Let the notations be as above.

\begin{proposition}
	For any $i,j  \in I$, we define $\boldsymbol\alpha_{ij} \in \Gamma(U_i \cap U_j, \cLog^\theN)$ 
	to be the element
	$$
		\boldsymbol\alpha_{ij} :=  \sum_{k=0}^\theN \frac{u_{ij}^{k+1}}{(k+1)!} \ul\omega_i^{0,k}.
	$$
	Then $(\bsomegav_{i}, \boldsymbol\alpha_{ij} )$ and $(\bsetav_i, 0)$ 
	satisfy the cocycle conditions
	\begin{align*}
		\nabla_{\kern-1mm\cL}(\boldsymbol\alpha_{ij}) &= \bsomegav_j - \bsomegav_i, &
		\boldsymbol\alpha_{ij} + \boldsymbol\alpha_{jk} &= \boldsymbol\alpha_{ik},&
		\bsetav_j - \bsetav_i &= 0, &  &
	\end{align*}
	hence define cohomology classes  in $H^1_{\log\dR}(E,  \cLog^\theN)$.
	These classes coincide with the classes  $[\bsomegav]$ and  $[\bsetav]$ in $H^1_\dR(U,  \cLog^\theN)$
	through the isomorphism
	$
		H^1_{\log\dR}(E,  \cLog^\theN) \cong H^1_{\dR}(U, \cLog^\theN).
	$
\end{proposition}

\begin{proof}
	By \eqref{eq: comparison epsilon} and the definition of $L^0_n$ and $L^0_{k,i}$, we have
	\begin{equation}\label{eq: base change}
		\sum_{n=0}^\theN \funcL_{n}^0 \ul\omega^{0,n} = 
		\sum_{k=0}^\theN \sum_{n=0}^k   \funcL^0_{ n} \frac{(-u_i)^{k-n}}{(k-n)!} \ul\omega_i^{0, k}
		= \sum_{k=0}^\theN \funcL^0_{k,i}\ul \omega_i^{0,k}.
	\end{equation}
	Hence $\bsetav = \bsetav_i = \bsetav_j$, which proves that $\bsetav_j - \bsetav_i = 0$.
	For any $i \in I$, we define $\boldsymbol\alpha_i \in \Gamma(U_i \setminus [0], \cLog^\theN)$
	to be the element
	$$
		\boldsymbol\alpha_i := - \sum_{k=0}^\theN \frac{(-u_i)^{k+1}}{(k+1)!} \ul\omega_i^{0,k}.
	$$
	Then using \eqref{eq: Log N}, we see that
	$
		\boldsymbol\alpha_{ij} = \boldsymbol\alpha_j - \boldsymbol\alpha_i
	$
	if $U_i \cap U_j \not=\emptyset$.
	The equality $\boldsymbol\alpha_{ij} + \boldsymbol\alpha_{jk} = \boldsymbol\alpha_{ik}$
	follows immediately from this fact.
	By definition of the connection, $\nabla_{\kern-1mm\cL}(\boldsymbol\alpha_i)$ is given by
	\begin{multline*}
		\sum_{k=0}^\theN \frac{(-u_i)^{k}}{k!}\ul \omega_i^{0,k} \otimes d u_i -
		 \sum_{k=0}^{\theN-1} \frac{(-u_i)^{k+1}}{(k+1)!}
		\left( \ul\omega_i^{1,k} \otimes \omega + \ul\omega_i^{0,k+1} \otimes \eta_i   \right)
		\\
		= \ul \omega_i^{0,0} \otimes  d u_i - \sum_{k=1}^\theN \frac{(-u_i)^{k}}{k!}\ul \omega_i^{0,k} \otimes \eta
		- \sum_{k=1}^{\theN} \frac{(-u_i)^{k}}{k!} \ul\omega_i^{1,k-1} \otimes \omega.
	\end{multline*}
	Here, we have used the fact that $d u_i = \eta_i - \eta$.  Again by \eqref{eq: comparison epsilon}, we have
	$$
		\sum_{n=1}^{\theN} \funcL^0_{n} \ul\omega^{1,n-1} 
		=\sum_{k=1}^{\theN} \sum_{n=1}^{k} \funcL^0_{ n} \frac{(-u_i)^{k-n}}{(k-n)!} \ul\omega_i^{1,k-1}
		=\sum_{k=1}^{\theN} \left( \funcL^0_{ k, i} -  \frac{(-u_i)^{k}}{k!} \right)\ul \omega_i^{1,k-1}.
	$$
	Applying this equality and  \eqref{eq: comparison epsilon} of Corollary \ref{cor: log isom},
	we see that
	\begin{equation}\label{eq: alpha boundary}
		\nabla_{\kern-1mm\cL}(\boldsymbol\alpha_i)  = \bsomegav_{i} - \bsomegav.
	\end{equation}
	This shows that
	$
		\nabla_{\kern-1mm\cL}(\boldsymbol\alpha_{ij}) 
		= \nabla_{\kern-1mm\cL}(\boldsymbol\alpha_j) - \nabla_{\kern-1mm\cL}(\boldsymbol\alpha_i) = 
		\bsomegav_j - \bsomegav_i.
	$
	The classes $(\bsomegav_{i}, \boldsymbol\alpha_{ij} )$ and $(\bsetav_i, 0)$ 
	coincide with $[\bsomegav]$ and $[\bsetav]$, since by \eqref{eq: alpha boundary},
	the element
	$
		(\boldsymbol\alpha_{i})_{i \in I} \in \prod_i\Gamma(U_i \setminus [0], \cLog^\theN)
	$
	gives the coboundary of the difference between $\bsomegav$ and $\bsomegav_i$ and
	$\bsetav = \bsetav_i$ on each $U_i$.
\end{proof}

\begin{proof}[Proof of Theorem \ref{thm: de Rham}]
	By the previous proposition, we may calculate the residues of
	$[\bsomegav]$ and $[\bsetav]$ using the cocycles $(\bsomegav_{i}, \boldsymbol\alpha_{ij} )$ and $(\bsetav_i, 0)$. 
	By Proposition \ref{prop: L good}, the functions $\funcL_{ k, i}$ are holomorphic on $U_i$ for $k\not= 1$, 
	and the residue of $\funcL_{1,i} = u_i$ at $0$ is one.  Hence from the definition of
	$\bsomegav_i$, $\bsetav_i$ and the definition of the residue map, we have
	$\res([\bsomegav]) =   \omega_{[0]}^{1,0}$ and $\res([\bsetav]) =   \omega_{[0]}^{0,1}$.
	Our theorem now follows from the definition of $\bfp^\theN$.
\end{proof}

%
\subsection{Change of basis}
%

In defining the logarithm and the polylogarithm sheaves, we took the basis $\{ \ul\omega, \ueta \}$
of $H^1_\dR(E)$ over $F$.   In this section, we introduce a different basis which is 
only defined over $F(e_2^*)$, but useful for explicit calculations.   
In what follows, we assume that $e_2^* \in F$.
We define $\omegas$ to be the differential of the second kind 
\begin{equation}\label{eq: omegas}
	\omegas := - \eta - e_2^* \omega \quad \in  \quad \Gamma(E, \Omega^1_E(2[0])),
\end{equation}
which by definition is defined over $F$ since $e_2^* \in F$. 
We let $\uomegas$ be the class represented by $\omegas$ in $H^1_\dR(E)$.  Let $\frU = \{ U_i \}_{i \in I}$
be an affine covering of $E$, and we take $(\eta_i, u_{ij})$ and $(u_i)$ as in Definition \ref{def: basis}.
Let
$$
	\omega^*_i := - \eta_i - e_2^* \omega \quad \in \quad \Gamma(U_i, \Omega^1_E). 
$$
Then $(\omega^*_i, -u_{ij} )$ is a cocycle which represents the same class in $H^1_\dR(E)$
as $\uomegas$, and the coboundary between $(\omega^*_i, -u_{ij} )$ and $(\omega^*,  0)$  is given 
by $(-u_i)$.

The classes $\{ \ul\omega, \uomegas \}$ form a basis of $H^1_\dR(E)$, and 
we denote by $\{ \wuomegav, \wuomegasv \}$ its dual basis.
We first describe the logarithm sheaf using this basis.
The definition \eqref{eq: omegas} of $\omegas$ gives the relation $\wuomegav = \uomegav - e_2^* \uetav$ and $\wuomegasv = - \uetav$.
Using this relation, we have that $\bnu$ and $\bnu_i$ is equal to
\begin{align*}
	\bnu &= \wuomegav \otimes \omega + \wuomegasv \otimes \omegas, &
	\bnu_i &= \wuomegav \otimes \omega + \wuomegasv \otimes \omegas_i.
\end{align*}
Write the restrictions of $\cLog^{(1)}$ to $U$ and $U_i$ as 
\begin{align*}
	\cLog^{(1)}|_{U} &= \cO_{U} \ul e\bigoplus \cH_U, &
	\cLog^{(1)}|_{U_i} &= \cO_{U_i} \ul e_i \bigoplus \cH_{U_i},
\end{align*}
with connection $\nabla(\ul e) = \bnu$ and $\nabla(\ul e_i) = \bnu_i$.
Then the pasting isomorphisms of Proposition \ref{pro: log dR on U}  on $U \cap U_i$ is given by
\begin{align}
	\ul e  &= \ul e_i + u_i \wuomegasv.
\end{align}
The elements $\wt{\ul\omega}^{m,n} := \ul e^a \wt{\ul\omega}^{\vee m} \wt{\ul\omega}^{*\vee n}/a!$ for $a = N - m-n$
form a basis of $\cLog^\theN$ restricted to $U$.

Next we describe the polylogarithm class. Let 
\begin{align*}
	\wt{\bomega}^\vee &=  -\wt{\ul\omega}^{0,0} \otimes \omegas + \sum_{n=1}^{\theN}  \funcL_{n}   \wt{\ul\omega}^{1,n-1} 
	\otimes  \omega, &
	\wt{\bomega}^{*\vee} &= \sum_{n=0}^\theN  \funcL_{ n} \wt{\ul\omega}^{0,n} \otimes  \omega.
\end{align*}

\begin{lemma}
	We have
	$$
		\bfp^\theN = \ul\omega \otimes [\wt{\bomega}^\vee]  + \ul\omega^* \otimes [\wt{\bomega}^{*\vee}].
	$$
\end{lemma}

\begin{proof}
	We have $\wt{\ul\omega}^{0,n} = (-1)^n \ul\omega^{0,n}$ and
	$\wt{\ul\omega}^{1,n} = (-1)^{n} \ul\omega^{1,n}	+ e_2^* (-1)^{n+1}\ul\omega^{0,n+1}$.
	Our assertion follows by directly calculating the right hand side of the above equality
	using the definition of  $\wt{\bomega}^\vee$ and $\wt{\bomega}^{*\vee}$.
\end{proof}

\begin{remark}
	Similarly to Corollary \ref{cor: dR P}, the above lemma shows that the polylogarithm sheaf $\cP^\theN$ may be given
	as a coherent $\cO_{U}$-module 
	$$
		\cP^\theN := \cH_U \bigoplus   \cLog^\theN_{U},
	$$
	with connection $\nabla_{\kern-1mm\cP}$ given by $\nabla_{\kern-1mm\cL}$ on $\cLog^\theN$ and
	\begin{align*}
		\nabla_{\kern-1mm\cP}(\wt{\ul\omega}^\vee) &=  \wt{\bomega}^{\vee},  &
		\nabla_{\kern-1mm\cP}(\wt{\ul\omega}^{*\vee}) &= \wt{\bomega}^{*\vee}.
	\end{align*}
\end{remark}

\begin{remark}\label{rem: choice of basis}
	In the rest of this paper, we only need descriptions corresponding to the basis $\{ \ul\omega, \uomegas \}$.
	Henceforth  for simplicity, we will write $\uomegav$, $\uomegasv$, $\bsomegav$, ${\bomega}^{* \vee}$ and 
	$\ul\omega^{m,n}$
	for $\wuomegav$, $\wuomegasv$, $\wt{\bomega}^\vee$, $\wt{\bomega}^{* \vee}$ and $\wt{\ul\omega}^{m,n}$.
\end{remark}

%
%
\section{Classical and $p$-adic Eisenstein-Kronecker numbers}\label{section: pEK}
%
%

In this section, we will first review the definition of Eisenstein-Kronecker numbers.
We will then review the construction of the $p$-adic distribution interpolating Eisenstein-Kronecker
numbers, when the corresponding elliptic curve has complex multiplication.  We will then use this
distribution to define the $p$-adic Eisenstein-Kronecker numbers.  In what follows, we fix a lattice 
$\Gamma \subset \bbC$.

%
\subsection{Eisenstein-Kronecker numbers}
%

In this subsection, we first review the definition of the Eisenstein-Kronecker-Lerch series and the Eisenstein-Kronecker numbers.
We will then review the fact that generating function of Eisenstein-Kronecker numbers is given by the Kronecker theta function $\Theta(z,w)$.

We let  $\pair{z,w} := \Exp{(z \ol w -\ol z w)/A}$ where $A$ is again the fundamental area of $\Gamma$ divided by $\pi$.

\begin{definition}[\cite{We} VIII \S12]
	Let $z_0$, $w_0 \in \bbC$.
	For any integer $a \geq 0$, we define the Eisenstein-Kronecker-Lerch series $K^*_a(z_0, w_0,s)$ by
	$$
		K^*_a(z_0, w_0,s) := {\sum_{\gamma \in \Gamma}}^* \frac{(\ol z_0 +\ol \gamma)^a}{|z_0 + \gamma|^{2s}} 
		\pair{\gamma, w_0},
	$$
	where $*$ denotes the sum over all $\gamma \in \Gamma$ such that $\gamma \not= -z_0$.
	The series converges for $\Re(s)> \frac{a}{2} +1$, and extends to a meromorphic function on the 
	complex plane by analytic continuation.
\end{definition}

Suppose $z_0$, $w_0 \in \bbC$.  In our pervious paper, we defined Eisenstein-Kronecker numbers $e^*_{a,b+1}(z_0, w_0)$ for 
$a$, $b \geq 0$  (See \cite{BK1}, Definition 1.5).   In this section, we first review in detail the properties of the Eisenstein-Kronecker-Lerch series given in \cite{We} \rm{VIII} \S 13 in order to define  Eisenstein-Kronecker numbers for integers $a<0$.
We will later prove that the specializations of the elliptic polylogarithm are related to such Eisenstein-Kronecker numbers.

Let $a$ be an integer $\geq 0$.  We define $\theta^*_a(t,z_0,w_0)$ to be the function
$$
	\theta^*_a(t,z_0,w_0) = {\sum_{\gamma \in \Gamma}}^* \exp(- t|z_0 + \gamma|^2/A) \pair{\gamma, w_0} (\ol z_0 +\ol \gamma)^a.
$$
Then $\theta^*_a(t,z_0,w_0) = \Theta^*_a(t/A, z_0, w_0)$ using the notations of \cite{We} \rm{VIII} \S 13.
Then \cite{We} \rm{VIII} \S 13 (29) gives the equality
$$
	A^{s} \Gamma(s) K^*_{a}(z_0, w_0, s) =  \int_{0}^\infty \theta_a^*(t,z_0,w_0) t^{s-1} dt
$$
for $\Re(s)$ sufficiently large.   We define $\theta_a(t, z_0, w_0)$ similarly to $\theta_a^*(t, z_0, w_0)$, but by taking the sum over all 
$\gamma \in \Gamma$.  Then $\theta_0(t,z_0,w_0) = \theta_0^*(t,z_0,w_0) + \pair{w_0,z_0}$ if $a=0$ and $z_0 \in \Gamma$,
and  $\theta^*_a = \theta_a$ otherwise.
We have by \cite{We} \rm{VIII} \S 13 (30) the equality
$$
	\theta_a(t, z_0, w_0) = t^{-a-1} \theta_a(t^{-1}, w_0, z_0) \pair{w_0, z_0}.
$$
We let
$$
	I_a(z_0, w_0, s) = \int_{1}^\infty\theta^*_a(t,z_0,w_0) t^{s-1} dt,
$$
which is an analytic function in $s$ defined for any $s \in \bbC$.   Then  \cite{We} \rm{VIII} \S 13 (31) gives the equality
\begin{multline}\label{eq: integral expression}
	A^s \Gamma(s) K_a^*(z_0, w_0, s) = I_a(z_0, w_0, s) - \frac{\delta_{a, z_0}}{s} \pair{w_0,z_0}\\ 
	+ I_{a} (w_0, z_0, a+1-s) \pair{w_0,z_0}  + \frac{\delta_{a, w_0}}{s-1},
\end{multline}
where $\delta_{a,x} = 1$ if $a =0$ and $x \in \Gamma$, and $\delta_{a,x} = 0$ otherwise.   The above integral
expression gives the meromorphic extension of $K_a^*(z_0, w_0, s)$ to the whole complex plane.
We extend the definition of Eisenstein-Kronecker-Lerch series to integers $a < 0$ by 
$$
	K^*_a(z_0, w_0,s) = (-1)^a \overline{K^*_{-a}(-z_0,w_0,\ol s-a)}.
$$
Note that we have
$$
	K^*_a(z_0, w_0,s) = {\sum_{\gamma \in \Gamma}}^* \frac{(\ol z_0 +\ol\gamma)^a}{|z_0 + \gamma|^{2s}} \pair{\gamma, w_0}
$$
if $\Re(s) > a/2+1$.  We have the following proposition.

\begin{proposition}\label{proposition: properties of K}
  Let $a$ be an integer.
  \begin{enumerate}
     \renewcommand{\theenumi}{\roman{enumi}}
     \renewcommand{\labelenumi}{(\theenumi)}
     \item The function $K^*_a(z_0,w_0,s)$  for $s$ continues meromorphically to a function on 
     the whole $s$-plane, with a simple pole only at $s=1$ if $a=0$ and $w_0 \in \Gamma$.
     \item The functions $K^*_a(z_0,w_0,s)$ satisfy the functional equation
      \begin{equation} \label{equation: functional equation for K}
         		\Gamma(s) K^*_a(z_0,w_0,s) = A^{a+1-2s} \Gamma(a+1-s) K^*_a(w_0,z_0,a+1-s) \pair{w_0,z_0}.
      \end{equation}
\end{enumerate}
\end{proposition}

\begin{proof}
	The result for $a \geq 0$ is given in \cite{We} \rm{VIII} \S 13 and is proved as follows.
	The fact that $K^*_a(z_0,w_0,s)$ is holomorphic in $s$ for $a > 0$ follows from the integral expression 
	\eqref{eq: integral expression}.  The same integral expression shows that the only poles of
	$\Gamma(s) K^*_0(z_0,w_0,s)$ are a simple pole at $s=0$ if $z_0 \in \Gamma$ and at 
	$s=1$ if $w_0 \in \Gamma$.  Since $\Gamma(s)$ has a simple pole at $s=0$ and $\Gamma(1) =1$,
	we see that the unique pole of $K^*_0(z_0,w_0,s)$ is a simple pole at $s=1$ if $w_0 \in \Gamma$.
	The functional equation follows again from \eqref{eq: integral expression}.
	For integers $a<0$, the fact that $K^*_a(z_0,w_0,s)$ is holomorphic in $s$ follows from the definition and
	the statement for $a>0$.
	Noting that $\ol{\Gamma(s)} = \Gamma(\ol s)$,
	the functional equation applied to $ \Gamma(\ol s-a) K^*_{-a}(-z_0, w_0, \ol s-a)$
	gives the functional equation
	\begin{multline} \label{eq: wow wow}
         		\Gamma(s-a) K^*_a(z_0,w_0,s) \\= (-1)^a A^{a+1-2s} \Gamma(1-s) K^*_a(w_0,z_0,a+1-s) \pair{w_0,z_0}.
	\end{multline}
	Since $\Gamma(s)\Gamma(1-s) = \pi/\sin \pi s$, we have 
	$$
		\Gamma(s)/\Gamma(s-a) = (-1)^a \Gamma(a+1-s)/\Gamma(1-s).
	$$
	The functional equation for $a<0$ follows by multiplying the above quotient to both sides of \eqref{eq: wow wow}.
\end{proof}

We define the Eisenstein-Kronecker number $e^*_{a,b}(z_0, w_0)$ as follows.

\begin{definition}
	Let $z_0$, $w_0 \in \bbC$, and we let $a$ and $b$ be integers such that $(a,b) \not=(-1,1)$ if $w_0 \in \Gamma$.  
	We define the Eisenstein-Kronecker number $e^*_{a,b}(z_0,w_0)$ by
	$$
		e^*_{a,b}(z_0,w_0) = K^*_{a+b}(z_0,w_0, b).
	$$
\end{definition}

The Eisenstein-Kronecker number for $z_0 \in \Gamma$ and $(a,b)=(0,0)$ may be calculated explicitly
as follows.  The integral expression \eqref{eq: integral expression} shows that the residue of 
$\Gamma(s) K^*_0(z_0, w_0, s)$ at $s=0$ is 
$-\pair{w_0,z_0}$.  Since $\Gamma(s)$ has a simple pole of residue one at $s=0$, this implies that 
\begin{equation}\label{eq: value 00}
	e^*_{0,0}(z_0,w_0) := K^*_0(z_0,w_0,0) = - \pair{w_0,z_0}.
\end{equation}
In this paper, we use the following version of Eisenstein-Kronecker numbers.

\begin{definition}\label{def: EK}
	Suppose $a$ and $b$ are integers, and let $z_0 \in \bbC$ such that $z_0 \not\in\Gamma$ if $(a,b) = (-1,1)$.
	We let
	$$
		e^*_{a,b}(z_0) := e^*_{a,b}(0, z_0) = K^*_{a+b}(0, z_0, b).
	$$
\end{definition}

W will next consider the generating function of Eisenstein-Kronecker numbers and relate this function
to the connection function $L_n(z)$.   For any $z_0$, $w_0 \in \bbC$, we let 
$$
	\Theta_{z_0, w_0}(z,w) := \Exp{- \frac{z_0 \ol w_0}{A}} \Exp{- \frac{z \ol w_0 + w \ol z_0}{A}}
	\Theta(z+z_0,w+w_0).
$$
The following theorem is fundamental in relating  $\Theta(z,w)$ to the Eisenstein-Kronecker numbers.

\begin{theorem}\label{thm: generating function}
	The Laurent expansion of $\Theta_{z_0, w_0}(z,w)$ at the origin is given by
	$$
		\Theta_{z_0, w_0}(z,w) = \pair{w_0, z_0}\frac{\delta_{z_0}}{z} 
		+ \frac{\delta_{w_0}}{w}  + \sum_{a, b\geq 0}  
		(-1)^{a+b} \frac{e^*_{a, b+1}(z_0, w_0)}{a! A^a}  z^{b} w^{a},
	$$
	where $\delta_x = 1$ if $x \in \Gamma$ and $\delta_x =0$ otherwise.
\end{theorem}

\begin{proof}
	This is obtained from \S 1.4 Theorem 1.17 of \cite{BK1}, by replacing the index $a$, $b$ by $a$, $b+1$.
\end{proof}

By definition of $\Theta_{z_0,w_0}(z,w)$, we have $\Theta_{z_0,w_0}(z,w) = \pair{w_0,z_0} \Theta_{w_0,z_0}(w,z)$.
Hence by Theorem \ref{thm: generating function}, we have
\begin{equation*}
	\Theta_{z_0, 0}(z,w) = \Theta_{0,z_0}(w,z) = \frac{\delta_{z_0}}{z} + \frac{1}{w} + \sum_{a, b \geq 0} (-1)^{a+b}
	\frac{e^*_{a, b+1}(z_0)}{a! A^a} z^a w^b.
\end{equation*}
Using this function, we next define the function $F_{z_0, b}(z)$,
which will be used in the construction of the elliptic polylogarithm.

\begin{definition}
	For any $z_0 \in \bbC$, we let $F_{z_0,b}(z)$ be the function such that
	$$
		\Theta_{z_0,0}(z,w) = \sum_{b \geq 0} F_{z_0,b}(z) w^{b-1}.
	$$
\end{definition}

We have in particular 
$$
	F_{z_0,1}(z) := \lim_{w \rightarrow 0}\left( \Theta_{z_0,0}(z,w) - \frac{1}{w} \right) = F_1(z + z_0) - \ol z_0/A.
$$
The function $\Theta_{z_0, 0}(z,w)$ satisfies 
\begin{multline*}
	\Theta_{z_0 + \gamma, 0}(z,w) = \Exp{- \frac{w(\ol z_0 + \ol \gamma)}{A}} \Theta(z + z_0 + \gamma, w) \\
	= \Exp{- \frac{w \ol z_0}{A} } \Theta(z + z_0, w) = \Theta_{z_0, 0}(z,w)
\end{multline*} 
for any $\gamma \in \Gamma$.  Hence $F_{z_0, b}(z)$ depends only on the choice of $z_0$ modulo $\Gamma$.
Henceforth, the $z_0$ of $F_{z_0,b}$ will either denote an element in $\bbC$ or a class in 
$\bbC/\Gamma$.  We have from Theorem \ref{thm: generating function} the following corollary.

\begin{corollary}\label{corollary: generating}
	For any $b \geq 0$, the Laurent expansion of $F_{z_0, b}(z)$ at $0$ is given by
	$$
		F_{z_0, b}(z) = \frac{\delta_{b-1, z_0}}{z}  + \sum_{a \geq 0} (-1)^{a+b-1} \frac{e^*_{a, b}(z_0)}{a! A^a} z^a,
	$$
	where $\delta_{b,x}=1$ if $b=0$ and $x \in \Gamma$, and is zero otherwise.
\end{corollary}

\begin{proof}
	The statement for $b > 0$ is Theorem \ref{thm: generating function}.  The statement for $b=0$ follows from the fact
	that $e^*_{0,0}(z_0) = -1$ and $e^*_{a,0}(z_0) = 0$ for any $a > 0$ (See Remark \ref{rem: zero} for a proof 
	of this fact).
\end{proof}

The above result shows that $F_{z_0, b}(z)$ for $b \geq 0$ is a one-variable generating function for the 
Eisenstein-Kronecker numbers.  We will use this function to construct the $p$-adic distribution,
originally due to Manin-Vishik and Katz if $p$ is ordinary and Fourquaux and Yamamoto when $p$ is
supersingular,  interpolating such values.  We will then use this function to define and investigate the 
$p$-adic elliptic polylogarithm function.

The relation of this function to the connection function is given by the equality
\begin{equation}\label{eq: relation to connection}
	\Theta_{z_0,0}(z,w) := \exp\left[ \frac{- w \ol z_0}{A}\right] \Theta(z + z_0,w) = \exp( F_{z_0,1}(z) w ) \Xi(z + z_0,w).
\end{equation}

%
\subsection{$p$-adic Eisenstein-Kronecker numbers}\label{ss: EK number}
%

We next give the definition of $p$-adic Eisenstein-Kronecker numbers, when the corresponding elliptic
curve has complex multiplication.   We assume that $\Gamma$ has complex multiplication by the ring of 
integers of an imaginary quadratic field $\bsK$ of class number one, and suppose $p \geq 5$ is a prime which does not 
ramify in $\bsK$.   If we replace $\Gamma$ by a suitable constant multiple, the theory of complex multiplication
asserts that there exists a Weierstrass equation
\begin{equation}\label{eq: integral Weierstrass}
	E : y^2 = 4 x^3 - g_2 x - g_3, \qquad \omega = dx/y
\end{equation}
with coefficients in the ring of integers $\cO_{\bsK}$ of $\bsK$ and good reduction at the primes above $p$,
such that $\Gamma$ is the period lattice of $E$ with respect to the invariant differential $\omega=dx/y$.
We assume that $\Gamma$ satisfies this additional condition.
In this case, the following theorem was proved by Damerell.
\begin{theorem}[Damerell]
	Let $a$ and $b$ be integers $\geq 0$, and let $z_0 \in \Gamma \otimes \bbQ$,
	which corresponds to a torsion point of $E$.  Then we have
	$$
		e^*_{a,b}(z_0)/A^a  \in \ol\bsK.
	$$
	Moreover, we have $e_2^*: = e_{0,2}^*(0) \in \bsK$.  
\end{theorem}
See \cite{BK1} Corollary 2.10 for a proof using $\Theta(z,w)$.  
The above theorem shows that the Laurent coefficients of $\Theta_{z_0,0}(z,w)$ at the origin are in $\ol \bsK$.
 In what follows, we assume that $z_0 \in \Gamma \otimes \bbQ$. 

We denote by $\psi := \psi_{E/\bsK}$ the Gr\"ossencharacter of $\bsK$ associated to $E$.
We fix a prime $\frp$ of $\cO_\bsK$ over $p$, and we let $\pi := \psi_{E/\bsK}(\frp) \in \cO_{\bsK}$.
Then $\pi$ is a generator of $\frp$. 

Henceforth, we fix an embedding $i_p: \ol\bsK \hookrightarrow  \bbC_p$ such that the completion of $\bsK$ in 
$\bbC_p$ is $\bsKfrp$.    We let $K$ be a finite unramified extension of $\bsKfrp$ in $\bbC_p$. 
In what follows, we let $\wh E$ be the formal group associated to 
$E \otimes_{\cO_\bsK} \cO_{\bsKfrp}$ at the origin $s= -2x/y$.  Then $\wh E$ is a 
Lubin-Tate group over $\cO_{\bsKfrp}$, of height one or two depending on whether $E$ has ordinary or 
supersingular reduction at $\frp$.  We denote by $\lambda(t)$ the formal logarithm of $\wh E$ normalized so 
that $\lambda'(0) = 1$, and we denote by $\wh\Theta_{z_0, 0}(s,t)$ the formal composition of the two-variable 
Laurent expansion of $\Theta_{z_0, 0}(z,w)$ at $z=0$ and $w=0$ with the formal power series $z = \lambda(s)$, 
$w = \lambda(t)$. We let
$
	\partial_{s, \log} := \lambda'(s)^{-1} \partial_s= \partial_z.
$
We first investigate the $p$-adic properties of the functions $F_{z_0,b}(z)$,
which will be necessary to define the $p$-adic distribution used in defining the $p$-adic
Eisenstein-Kronecker numbers.

We denote by $\wh F_{z_0, b}(s)$ the formal composition of the 
Taylor expansion of $F_{z_0, b}(z)$ at $z=0$ with the power series $z = \lambda(s)$. 
Note that by definition and \eqref{eq: relation to connection}, we have
\begin{equation}\label{eq: F and L}
	\wh F_{z_0,b}(s) = \sum_{n=0}^b \frac{\wh F_{z_0,1}(s)^{b-n}}{(b-n)!} \wh L_{z_0, n}(s),
\end{equation}
where $\wh L_{z_0,n}(s)  := L_n(z + z_0)|_{z = \lambda(s)}$. 

\begin{remark}\label{rem: expansion}
	Suppose $f(z)$ is a meromorphic function on $\bbC/\Gamma$ corresponding
	to a rational function on $E$.  Then the power series $\wh f(s):=f(z)|_{z=\lambda(s)}$ is 
	in fact the expansion of the rational function of $f$ with respect to the formal parameter $s = -2x/y$
	at the identity of the elliptic curve.
\end{remark}

In investigating $\frp$-integrality of power series, the following lemma will play an important role.

\begin{lemma}\label{lem: prev}
	Let $K$ be a finite extension of $\bsKfrp$, and suppose that the meromorphic function $f(z)$ on 
	$\bbC/\Gamma$ corresponds to a rational function $f$ on $E$ defined over $K$, without any pole on 
	$\wh E(\frm_{\bbC_p})\setminus[0]$.  Then $\wh f(s) := f(z)|_{z=\lambda(s)}$
	has bounded coefficients.
\end{lemma}

\begin{proof}
	Consider the embedding $K(E) \hookrightarrow \operatorname{Frac}(\cO_K[[s]])$
	of the functional field of $E$ to the fractional field of the ring of formal power series with respect to $s$.
	Note that the image $\wh f(s)$ of $f$ is of the form $\wh f(s) = \wh P(s)/\wh Q(s)$, where
	$\wh P(s)$, $\wh Q(s) \in \cO_K[[s]]$ are relatively prime.  By the $p$-adic Weierstrass preparation
	theorem, $\wh Q(s)$ is of the form $\wh Q(s) = p^m \wh U(s) R(s)$, where $\wh U(s) \in \cO_K[[s]]^\times$
	and $R(s)$ is a \textit{distinguished polynomial}, in other words $R(s) \equiv s^N$ modulo $\frp$ for $N= \deg \, R(s)$.
	If $R(s) \not\equiv s^N$, then the non-zero roots of this polynomial would correspond to poles of $f$ on 
	$\wh E(\frm_{\bbC_p})$  other than zero. Hence by our assumption, we must have $R(s) = s^N$.  This shows that 
	$p^m \wh f(s) \in \cO_K[[s]][s^{-1}]$, proving the lemma.
\end{proof}

For any torsion point $z_0 \in E(K)$, we define the order of $z_0$ to be the annihilator of $z_0$ as an
element in an $\cO_\bsK$-module.  
The following result is fundamental.

\begin{proposition}
	Suppose $z_0 \in E(\ol\bbQ)$ is a non-zero torsion point of order $\frn$ prime to $\frp$.
	Then the power series
	$$
		\wh F_{z_0, b}(s) \in \ol\bsK[[s]]
	$$
	converges on the open unit disc $B^-(0,1) := \{ s \in \bbC_p \mid |s|_p < 1 \}$.
	In particular, this series defines a rigid analytic function on $B^-(0,1)$.
\end{proposition}

\begin{proof}
	When $b=0$ then there is nothing to prove. Since $F_1(z) = \zeta(z) - e_2^* z$, we have 
	$\partial_z F_1(z) = - \wp(z) - e_2^*$.  Hence 
	$$
		\partial_z F_{z_0,1}(z) = \partial_z F_1(z+z_0) = - \wp(z + z_0) - e_2^*.
	$$
	By Lemma \ref{lem: prev}, $\wh\wp_{z_0}(s) := \wp(z + z_0)|_{z=\lambda(s)}$ is known to have bounded coefficients
	(in fact, one may prove that the coefficients are $\frp$-integral).  Hence this power series converges on $B^-(0,1)$,
	which implies that $\wh F_{z_0,1}(s)$ also converges on $B^-(0,1)$. The assertion for general $b$ follows 
	from \eqref{eq: F and L}, noting that by Lemma \ref{lem: prev}, $\wh L_{z_0,n}(s)$ also has bounded coefficients.
\end{proof}

\begin{remark} When $p$ is a prime of ordinary reduction, then we may prove that the coefficients of  $\wh F_{z_0,1}(s)$,
	hence that of $\wh F_{z_0,b}(s)$, is bounded, as follows. By \cite{BK1} Corollary 2.17, the Laurent expansion of 
	$\wh \Theta_{z_0,0}(s,t)$ is $\frp$-integral if the order of $z_0$ is prime to $p$.  The result is valid even for the case 
	$z_0 = 0$.    Hence the coefficients of $\wh F_{z_0,1}(s)$ are $\frp$-integral in this case.  For a general 
	$z_0 \in 	\Gamma \otimes \bbQ$ of order prime to $\frp$, let $\alpha \in \cO_\bsK$ be an element prime to $\frp$ such that 
	$\alpha z_0 \in \Gamma$.  If we let 
	$F^\alpha_1(z) = F_1(z) - \ol\alpha^{-1} F_1(\alpha z) = \partial_z \log 
	(\theta(z)^{N(\alpha)}/\theta(\alpha z))/N(\alpha)$, then 
	this function is known to be an elliptic function defined over $\bsK$ (See for example \cite{CW} \S5).
	Then by definition and our choice of $\alpha$, we have $F_{z_0,1}(z) = \ol\alpha^{-1} F_1(\alpha z) + F^\alpha_1(z+z_0)$.
	 Our assertion now follows from Lemma \ref{lem: prev} and the assertion for $z_0=0$, noting that $\alpha$ is a $\frp$-adic
	 unit and $F_1(\alpha z)|_{z=\lambda(s)} = \wh F_1([\alpha] s)$.
\end{remark}

We will use $\wh F_{z_0, b}(s)$ to construct our $p$-adic distribution.  Since the formal group 
$\wh E$ is a Lubin-Tate group,  it has an action of $\cO_{\bsKfrp}$.  We have 
$\cO_{\bsKfrp}$-linear isomorphisms
$$
	\Hom_{\cO_{\bbC_p}}( \wh E, \wh \bbG_m) \cong \Hom_{\bbZ_p}(T_p \wh E, T_p  \wh\bbG_m) 
	\xleftarrow\cong \cO_{\bsKfrp}.
$$
The last isomorphism is not canonical, and depends on the choice of a $p$-adic period as follows.
There exists $\Omega_\frp \in \bbC_p^\times$ such that the formal power series $\exp(\lambda(s)/\Omega_\frp)$
is an element in $\cO_{\bbC_p}[[s]]$.  The second isomorphism is given by associating to any 
$x \in \cO_{\bsKfrp}$ the homomorphism of formal groups defined by $\exp(x \lambda(s)/\Omega_\frp)$.  
This isomorphism depends on the choice of $\Omega_\frp$.
The notation $\exp(x \lambda(s)/\Omega_\frp)$ needs some care.  If $s_m$ is a 
primitive $p^m$-torsion  point in $\wh E(\frm_{\cO_{\bbC_p}})$, then $\lambda(s_m) = 0$ but 
$\exp(x \lambda(s)/\Omega_\frp)|_{s=s_m}$ is a primitive $p^m$-th root of unity.

In what follows, we fix once and for all a choice of a $p$-adic period $\Omega_\frp$.
Let $C^\an(\cO_{\bsKfrp}, \bbC_p)$ be the set consisting of locally $\bsKfrp$-analytic functions on $\cO_{\bsKfrp}$.
We define our $p$-adic distribution $\mu_{z_0,b}$ as follows.

\begin{definition}\label{def: distribution}
	Let $z_0$ be a non-zero torsion point in $E(\ol\bbQ)$ of order prime to $\frp$.  
	For any integer $b \geq 0$, we define $\mu_{z_0, b}$ to be the $p$-adic distribution 
	on $C^\an(\cO_{\bsKfrp}, \bbC_p)$ associated to $\wh F_{z_0, b}(s)$.  
	Such distribution satisfies the relation
	$$
		\int_{\cO_{\bsKfrp}} \exp\left(x \lambda(s)/\Omega_\frp\right) d \mu_{z_0, b}(x) = \wh F_{z_0, b}(s).
	$$
	When $p$ is ordinary, then this is the $p$-adic measure associated to bounded power series.
	When $p$ is supersingular, then this is the $p$-adic distribution associated to rigid analytic functions
	on the open unit disc constructed in \cite{ST} Theorem 2.3 and Theorem 3.6.
\end{definition}

When $p$ is ordinary, the above distribution is related to the two-variable measure used by Manin-Vishik 
and Katz in defining the two-variable $p$-adic $L$-function interpolating special values of Hecke $L$-function
of imaginary quadratic fields (See Proposition \ref{pro: relation}).  
When $p$ is supersingular, the above distribution was considered by Boxall \cite{Box1} \cite{Box2}
and Schneider-Teitelbaum \cite{ST}
for the case $b=0$, and by  Fourquaux \cite{Fou} and Yamamoto \cite{Yam} for any $b \geq 0$.

We may now define the $p$-adic Eisenstein-Kronecker numbers using this distribution.

\begin{definition}\label{def: p EK}
	Let $z_0$ be a non-zero torsion point in $E(\ol\bbQ)$ of order prime to $\frp$.  
	For any integers $a$ and $b$ such that $b \geq 0$, we define the \textit{$p$-adic Eisenstein-Kronecker number}
	$
		e^{(p)}_{a,b}(z_0)
	$
	by
	$$
		e^{(p)}_{a,b}(z_0) :=  \Omega_\frp^{b-1} \int_{\cO_{\bsKfrp}^\times} x^a d \mu_{z_0, b}(x),
	$$
	where we denote again by $\mu_{z_0, b}$ the restriction of $ \mu_{z_0, b}$ to $\cO_{\bsKfrp}^\times$.
\end{definition}

The above definition is justified  from the fact that the distribution $\mu_{z_0, b}$ interpolates Eisenstein-Kronecker numbers
$e^*_{a,b}(z_0)$ for $a, b \geq 0$.
This fact will be proved in \S \ref{subsection: interpolation} Corollary \ref{cor: interpolation}.

%
\subsection{Interpolation property of the $p$-adic distribution}\label{subsection: interpolation}
%

Here we will prove an interpolation property of the distribution  $\mu_{z_0, b}$,
justifying the definition of the $p$-adic Eisenstein-Kronecker numbers.  We will then prove
in Proposition \ref{pro: relation} that when $p \geq 5$ is a prime of ordinary reduction, 
then the definition given in Definition \ref{def: p EK} coincides with the definition given in the introduction.

We first begin with the distribution property of $\Theta_{z_0,  0}(z,w)$.

\begin{proposition}[Distribution relation]\label{pro: distribution relation}
	For any $z_0 \in \bbC$, we have
	$$
		 \Theta_{\pi^m z_0,  0}(\pi^m z, \ol\pi^{-m} w) = \frac{1}{\pi^m}
		 	\sum_{z_m \in \frac{1}{\pi^{m}} \Gamma/\Gamma} \Theta_{z_0 + z_m, 0}(z,w).
	$$
\end{proposition}

\begin{proof}
	Note that we have
	$$
		 \Theta_{\pi^m z_0,  0}(\pi^m z, \ol\pi^{-m} w; \Gamma) = \Theta_{N \frp^m z_0,  0}(N \frp^m z, w; \ol\frp^{m} \Gamma).
	$$
	Our assertion is a special case of \cite{BK1} Proposition 1.16.
\end{proof}

\begin{corollary}\label{cor: distribution}
	The function $F_{z_0, b}(z)$ satisfies the relation
	$$
		F_{\pi^m z_0, b}(\pi^m z) = \frac{\ol\pi^{mb}}{\bN^m} \sum_{z_m \in  \frac{1}{\pi^{m}}
		 \Gamma/\Gamma}  F_{z_0 + z_m, b}(z).
	$$
\end{corollary}

\begin{proof}
	The statement is trivial when $b=0$.  The case for $b \geq 1$
	follows from the distribution relation Proposition \ref{pro: distribution relation}
	for the Kronecker theta function and the definition of $F_{z_0, b}(z)$.
\end{proof}

In what follows, we again let $z_0$ be a non-zero torsion point of $E$ of order prime to $\frp$.
The power series $\wh F_{z_0, b}(s)$ satisfies the following translation formula with
respect to $\frp^m$-torsion points.

\begin{lemma}[Translation]\label{lem: translation}
	Recall that $\wh F_{z_0,b}(s)$ is the formal power series composition of the Taylor expansion of $F_{z_0,b}(z)$ 
	at the origin with $z = \lambda(s)$.  Then we have
	$$
		\wh F_{z_0, b}(s \oplus s_m) = \wh F_{z_0 + z_m, b}(s),
	$$
	where $s_m$ is a torsion point in $\wh E[\frp^m]$, and $z_m$ is the image of $s_m$ through the inclusion 
	$\wh E(\frm_p)_{\tor} \subset E(\ol\bbQ)_\tor \subset \bbC/\Gamma$.
\end{lemma}

\begin{proof}
	The statement is trivial when $b=0$.  When $b=1$,  the function 
	$$	
		F(z) := F_{z_0,1}(z) - \ol\pi^{-m} F_{z_0,1}(\pi^m z)
	$$
	is elliptic, hence satisfies 	$\wh F(s \oplus s_m) = F(z + z_m)|_{z = \lambda(s)}$. 
	By the equalities 
	$
		F_{z_0,1}(\pi^m z)|_{z=\lambda(s)} = \wh F_{z_0,1}([\pi^m]s)
	$ and
	$
		\wh F_{z_0,1}([\pi^m](s \oplus s_m)) = \wh F_{z_0,1}([\pi ^m] s),
	$
	we have 
	$$
		\wh F(s \oplus s_m) = \wh F_{z_0,1}(s \oplus s_m) -  \ol\pi^{-m} \wh F_{z_0,1}([\pi^m] s).
	$$
	In addition, we have by definition
	$F_{z_0,1}(z + z_m) = F_{z_0+z_m,1}(z) + \ol z_m/A$.  Hence applying the equality $F_{z_0 + \pi^m z_m,1}(\pi^m z) =
	 F_{z_0,1}(\pi^m z)$ for $\pi^m z_m \in \Gamma$, we have
	$$	
		 F(z + z_m) := F_{z_0 + z_m,1}(z) - \ol\pi^{-m} F_{z_0,1}(\pi^m z).
	$$
	Our assertion follows by combining the above results.
	The case for $b>1$ follows from \eqref{eq: F and L},
	applying our lemma for $b=1$ and noting that $\wh L_{z_0,n}(s \oplus s_m) = \wh L_{z_0+z_m,n}(s)$ since
	$L_n(z)$ corresponds to a rational function.
\end{proof}

Using the above lemma, we have the following.

\begin{proposition}\label{prop: distribution}
	The distribution $\mu_{z_0, b}$ restricted to $\cO_{\bsKfrp}^\times$ satisfies
	$$
		\int_{\cO_{\bsKfrp}^\times} \exp\left(x \lambda(s)/\Omega_\frp \right) d \mu_{z_0, b}(x) 
			= \wh F_{z_0, b}(s) -  \frac{1}{\ol \pi^b} \wh F_{\pi z_0, b}([\pi] s). 
	$$
\end{proposition}

\begin{proof}
	Note that for any primitive $\frp$-torsion point $s_1$ in $E[\frp]$, the value $\exp(\lambda(s)/\Omega_\frp)|_{s=s_1}$
	is a primitive $p$-th root of unity.  Hence standard argument shows that
	the restriction of distributions from $\cO_{\bsKfrp}$ to $\cO_{\bsKfrp}^\times$ is given by
	$$
		\int_{\cO_{\bsKfrp}^\times} \exp\left(x \lambda(s)/\Omega_\frp \right) d \mu_{z_0, b}(x) 
		= \wh F_{z_0, b}(s) - \frac{1}{\bN} \sum_{s_1\in E[\frp]} \wh F_{z_0, b}(s \oplus s_1).
	$$
	The formula for the formal translation by $\frp$-power torsion point in Lemma \ref{lem: translation} shows that 
	$\wh F_{z_0, b}(s \oplus s_1) = \wh F_{z_0 + z_1,b}(s)$, hence we have
	$$
		\int_{\cO_{\bsKfrp}^\times}  \exp\left(x \lambda(s)/\Omega_\frp \right)  d \mu_{z_0, b}(x) 
		=\wh F_{z_0, b}(s) - \frac{1}{\bN} \sum_{z_1 \in \frac{1}{\pi}\Gamma/\Gamma} \wh F_{z_0 + z_1, b}(s).
	$$
	Our result now follows from the distribution relation (Corollary \ref{cor: distribution}).
\end{proof}

The expansion of $F_{z_0,b}(z)$ given in Corollary \ref{corollary: generating} shows that
we have the following corollary, which implies that our distribution $\mu_{z_0,b}$ interpolates
the Eisenstein-Kronecker numbers.

\begin{corollary}\label{cor: interpolation}
	The distribution $\mu_{z_0, b}$ satisfies
	$$
		\Omega_\frp^{-a} \int_{\cO_{\bsKfrp}^\times} x^a d \mu_{z_0, b}(x) 
			= (-1)^{a+b-1}\left( \frac{e^*_{a,b}(z_0)}{A^a} -  \frac{\pi^a e^*_{a,b}(\pi z_0)}{\ol\pi^bA^a} \right)
	$$
	for any integer $a \geq 0$.  In particular, we have
	$$
		\frac{e^{(p)}_{a,b}(z_0)}{\Omega_\frp^{a+b-1}}
			= (-1)^{a+b-1} \left( \frac{e^*_{a,b}(z_0)}{A^a} -  \frac{\pi^a e^*_{a,b}(\pi z_0)}{\ol \pi^{b} A^a}  \right).
	$$
\end{corollary}
The above result shows that the $p$-adic Eisenstein-Kronecker numbers are related to usual Eisenstein-Kronecker
numbers when $a, b \geq 0$.

We now relate the definition of $e^{(p)}_{a,b}(z_0)$ given in Definition \ref{def: p EK} with the definition
given in the introduction, when the prime $p$ is ordinary.
Assume just for now that $p$ is a prime for which \eqref{eq: integral Weierstrass} has good ordinary reduction at $p$,
and suppose that the order of $z_0$ is prime to $p$.
In \cite{BK1}, we defined a two-variable $p$-adic measure $\mu_{z_0,0}$ on 
$\bbZ_p \times \bbZ_p$.    By substituting $\eta_{\frp}(t) = \exp(\lambda(t)/\Omega_\frp)-1$ into Definition 3.2 of \cite{BK1},
we see that the $p$-adic measure $\mu_{z_0,0}$ is defined to satisfy
$$
	\int_{\bbZ_p \times \bbZ_p} \exp(x \lambda(s)/\Omega_\frp) \exp(y \lambda(t)/\Omega_\frp) 
	d \mu_{z_0,0}(x,y) = \wh\Theta_{z_0,0}^*(s,t),
$$
where $\wh\Theta^*_{z_0,0}(s,t) := \wh\Theta_{z_0,0}(s,t) - t^{-1}$.  By taking $\partial_{t,\log}^{b}$ of both sides and substituting
$t=0$, we obtain the equality
$$
	\frac{1}{\Omega_{\frp}^{b}} \int_{\bbZ_p \times \bbZ_p} \exp(x \lambda(s)/\Omega_\frp) y^b 
	d \mu_{z_0,0}(x,y) =b! \wh F_{z_0,b+1}^*(s),
$$
where 
$$
	\wh F^*_{z_0, b+1}(s) = \frac{1}{b!}\lim_{t \rightarrow 0} \partial^{b}_{t, \log} 
	 \wh \Theta^*_{z_0, 0}(s,t).
$$
By definition, the relation between $\wh F_{z_0,b+1}(s)$ and $\wh F_{z_0,b+1}^*(s)$ for $b \geq 0$ is given by 
$
	\wh F_{z_0,b+1}^*(s) = \wh F_{z_0,b+1}(s) + c_{b+1},
$	
where $c_{b+1}$ is the constant 
$$
	c_{b+1} : = 	 \frac{1}{b!} 
	\lim_{t \rightarrow 0}\partial_{t,\log}^{b} \left(\lambda(t)^{-1} -t^{-1} \right)
$$
in $K$.  The usual formula for restriction to $\bbZ_p^\times \times \bbZ_p$ gives the equality
\begin{multline*}
	\frac{1}{b!\Omega_{\frp}^{b}} \int_{\bbZ_p^\times \times \bbZ_p} \exp(x \lambda(s)/\Omega_\frp) y^b 
	d \mu_{z_0,0}(x,y) \\= \wh F_{z_0,b}(s) -\frac{1}{\bN} \sum_{s_1\in E[\frp]} \wh F_{z_0, b}(s \oplus s_1)
		= \wh F_{z_0, b}(s) -  \frac{1}{\ol \pi^b} \wh F_{\pi z_0, b}([\pi] s).
\end{multline*}
Note that the $*$ is not required in the middle term since the constant $c_{b+1}$ cancels in the sum.
This shows that we have
$$
	\int_{\bbZ_p^\times \times \bbZ_p} x^a y^b 
	d \mu_{z_0,0}(x,y) = b! \Omega_{\frp}^{b} \int_{\bbZ_p^\times} x^a d \mu_{z_0,b+1}(x).
$$
Then the definition of  $e^{(p)}_{a,b}(z_0)$ gives the following.

\begin{proposition}\label{pro: relation}
	Suppose $p$ is an ordinary prime.
	Let $e^{(p)}_{a,b}(z_0)$ be the $p$-adic Eisenstein-Kronecker number defined in  Definition \ref{def: p EK}.
	Then we have
	$$
		e^{(p)}_{a,b+1}(z_0) =\frac{1}{b!} \int_{\bbZ_p^\times \times \bbZ_p} x^a y^b 
	d \mu_{z_0,0}(x,y) 
	$$
	for any integer $a,b$ such that $b \geq 0$.
\end{proposition}

%
%
%
\section{$p$-adic elliptic polylogarithm functions}
%
%
%

In this section, we define the $p$-adic elliptic polylogarithm function.  We then prove the relation
of this function to the $p$-adic distribution defined in the previous section.  We keep the notations
of \S \ref{ss: EK number}  We assume again that $p \geq 5$ is a prime of good reduction.

%
\subsection{$p$-adic elliptic polylogarithm functions}
%

In this section, we define the $p$-adic elliptic polylogarithm function.  
We first begin by defining a $p$-modified variant of the connection function.
For any $z_0 \in \bbC$, let 
$$
	\Theta^{(p)}_{z_0,0}(z,w) := \Theta_{z_0,0}(z,w) - \frac{1}{\ol\pi} 
	\Theta_{\pi z_0,0}(\pi z, \ol\pi^{-1} w),
$$
and $\Theta^{(p)}(z,w) := \Theta^{(p)}_{0,0}(z,w)$.

\begin{definition}
	Let $\Xi^{(p)}(z,w) := \exp(-F_{1}(z) w )  \Theta^{(p)}(z, w)$.
	For any integer $n \geq 0$, we define $L^{(p)}_{n}(z)$ to be the function given by
	$$
		\Xi^{(p)}(z,w)  = \sum_{n \geq 0} L^{(p)}_{n}(z)  w^{n-1}.
	$$
\end{definition}

We let $F^{(p)}_1(z)$ be the meromorphic function
$$
	F^{(p)}_1(z) := F_1(z) - \frac{1}{\ol\pi} F_1(\pi z).
$$
Then the transformation formula $F_1(z + \gamma) = F_1(z) + \ol \gamma/A$ for
$\gamma \in \Gamma$ shows that $F_1^{(p)}(z)$ is a rational function defined over $\bsK$.
By definition, $F_{z_0,1}(z) = F_1(z+z_0) - \ol z_0/A$.  Hence for any $z_0 \in \bbC$, we have
\begin{equation}\label{eq: funny}
		 F^{(p)}_1(z + z_0) 
		 = F_{z_0,1}(z) -\ol\pi^{-1} F_{\pi z_0,1}(\pi z).
\end{equation}

\begin{lemma}\label{lem: L p}
	Let $n$ be an integer $\geq 0$.   The function $L^{(p)}_{n}(z)$ is  periodic with respect to 
	$\Gamma$ and is holomorphic outside $\pi^{-1} \Gamma$. The meromorphic function $L^{(p)}_n(z)$
	corresponds to a rational function on $E$ defined over $\bsK$.
	In addition, we have
	\begin{equation}\label{eq: L translate}
		\Xi^{(p)}_{z_0}(z,w)  = \sum_{n \geq 0} L^{(p)}_{n}(z + z_0)  w^{n-1},
	\end{equation}
	where  $	\Xi^{(p)}_{z_0}(z,w) := \exp(-F_{z_0,1}(z)w) \Theta^{(p)}_{z_0,0}(z,w)$. 
\end{lemma}

\begin{proof}
	From the equality $\Theta(z,w) = \exp(F_1(z) w) \Xi(z,w)$, we have 
	$$
		\Theta^{(p)}(z,w) = \exp(F_1(z) w) \Xi(z,w) - \ol\pi^{-1} \exp(F_1(\pi z) w/\ol \pi)
		\Xi(\pi z, \ol\pi^{-1} w).
	$$
	Since $F_1^{(p)}(z) = F_1(z) - \ol\pi^{-1} F_1(\pi z)$, the above equality and
	the definition of $\Xi^{(p)}(z,w)$ shows that we have
	$$
		\Xi^{(p)}(z,w) = \Xi(z,w) - \ol\pi^{-1} \exp( - F_1^{(p)}(z)  w   )   \Xi(\pi z, \ol \pi^{-1} w).
	$$
 	This implies that $L_n^{(p)}(z)$ may be expressed as the difference between $L_n(z)$
	and a sum over products of $\ol\pi$, $F_1^{(p)}(z)$ and $[\pi]^* L_k(z)$.
	This proves that $L_n^{(p)}(z)$ is defined over $\bsK$ and holomorphic outside $\pi^{-1} \Gamma$,
	since the same holds true for $L_n(z)$, $F_1^{(p)}(z)$ and $[\pi]^* L_k(z)$.
	The last equality follows from the fact that we have $\Xi^{(p)}(z,w) = \Xi^{(p)}_{z_0}(z + z_0,w)$
	for $z_0 \in \bbC$,  which proves our assertion.
\end{proof}

Let $\cU$ be the formal completion of $U$ with respect to the special fiber, and we denote by $\cU_K$ the 
rigid analytic space associated to $\cU$.   Let $R = \Gamma(U, \cO_U)$, and we denote by $R^\dagger$ the weak completion 
of $R$ (\cite{MW} Definition 1.1). Then $R^\dagger_K : = R^\dagger \otimes K$ is the ring of overconvergent functions on $\cU_K$.
The following are the $p$-adic elliptic polylogarithm functions.

\begin{theorem}\label{thm: existence two}
	Let $D^{(p)}_{0,n} = L^{(p)}_n$ for any integer $n>0$, and $D^{(p)}_{m,n} = 0$ if 
	$n \leq 0$.
	Then for integers $m, n >0$,
	there exists a unique system of overconvergent functions $D^{(p)}_{m,n}$ on $\cU_K$
	iteratedly satisfying the differential equation
	\begin{equation}\label{eq: differential elliptic polylog}
		d D^{(p)}_{m,n} = - D^{(p)}_{m-1,n} \omega - D^{(p)}_{m,n-1} \omegas.
	\end{equation}
	We call the overconvergent functions $D^{(p)}_{m,n}$ the $p$-adic elliptic polylogarithm functions.
\end{theorem}

See Lemma \ref{lemma: first} for the differential equations satisfied by the elliptic polylogarithm functions
in the Hodge case.
The above theorem follows from the following stronger statement. 

\begin{theorem}\label{thm: existence one}
	We let $D^{(p)}_{m,n}$ for $m \leq 0$ or $n \leq 0$ as in Theorem \ref{thm: existence two}.
	Then for integers $m$, $n >0$,
	there exists a unique system of overconvergent functions $D^{(p)}_{m,n}$ on $\cU_K$
	iteratedly satisfying
	\begin{enumerate}
	\item The differential equation
	$
		d D^{(p)}_{m,n} = - D^{(p)}_{m-1,n} \omega - D^{(p)}_{m,n-1} \omegas.
	$ 
	\item The distribution relation 
	\begin{equation}\label{eq: distribution}
		\sum_{z_1 \in E[\frp](K) } G^{(p)}_{m,n}(z + z_1) = 0, 
	\end{equation}
	where 
	\begin{equation}\label{eq: E}
		G^{(p)}_{m,n} = \sum_{k=0}^n \frac{(F_1^{(p)})^{n-k}}{(n-k)!} D^{(p)}_{m,k}.
	\end{equation}
	\end{enumerate}
\end{theorem}

We will give the proofs of Theorem \ref{thm: existence one} and Theorem \ref{thm: existence two} at the end of this subsection.
We first give a lemma.

\begin{lemma}
	The distribution relation \eqref{eq: distribution} is true if $m=0$ or $n=0$.
\end{lemma}

\begin{proof}
	The statement for $n=0$ follows from the fact that $G_{m,0}^{(p)} = D^{(p)}_{m,0} = 0$.
	We prove the statement when $m=0$ and $n>0$. 
	By the definition of $\Theta^{(p)}_{z_0,w_0}(z,w)$, 
	the distribution relation in Proposition \ref{pro: distribution relation} for
	$\Theta_{z_0,w_0}(z,w)$ gives the relation
	\begin{equation}\label{eq: to zero}
		\sum_{z_1 \in E[\frp]} \Theta^{(p)}_{z_1,0}(z,w) = 0.
	\end{equation}
	If $z_1 \in E[\frp]$, then \eqref{eq: funny} gives the equality
	\begin{equation}\label{eq: fun}
		 F^{(p)}_1(z + z_1) = F_{z_1,1}(z) - \ol\pi^{-1} F_1(\pi z).
	\end{equation}
	We see from \eqref{eq: L translate} and the above equality that
	\begin{multline*}
		\Theta^{(p)}_{z_1,0}(z,w) = 
		\exp\left(F_{z_1,1}(z) w\right) \sum_{n \geq 0} L^{(p)}_{n}(z+z_1)  w^{n-1} \\
		=\exp\left( \frac{1}{\ol\pi} F_1(\pi z) w\right)
		\exp\left(F^{(p)}_1(z+z_1) w\right) 
		 \sum_{n \geq 0} L^{(p)}_{n}(z +z_1)  w^{n-1}
	\end{multline*}
	for any $z_1 \in E[\frp]$. Then \eqref{eq: to zero} translates to the equality
	$$
		\sum_{z_1 \in E[\frp]}
		\exp\left(F^{(p)}_1(z+z_1) w\right) 
		 \sum_{n \geq 0} L^{(p)}_{n}(z+z_1 ) w^{n-1} =0.
	$$
	Writing out the coefficient of $w^{n-1}$, we obtain 
	$$
		\sum_{z_1 \in E[\frp]} G^{(p)}_{0,n}(z+ z_1) =
		\sum_{z_1 \in E[\frp]} \sum_{k=0}^n \frac{F_1^{(p)}(z+ z_1)^{n-k}}{(n-k)!} L_k^{(p)}(z+ z_1) = 0,
	$$
	which is the desired equality.
\end{proof}

Berthelot defined rigid cohomology for any scheme of finite type over a field of characteristics $p > 0$ (\cite{Ber0}, \cite{Ber2}).
When the scheme is affine and smooth over the base, then rigid cohomology is canonically isomorphic to Monsky-Washnitzer 
cohomology  (\cite{Ber2} Proposition 1.10), which by definition may be calculated using overconvergent functions and differentials.  
Rigid cohomology $H^i_\rig(U_k/K)$ of $U_k$ may be calculated as
$$
	H^i_\rig(U_k/K) \cong H^i \left[  R^\dagger_K \xrightarrow d R^\dagger_K \otimes \Omega^1_{R}  \right].
$$
In order to determine if a differential form in $R^\dagger \otimes \Omega^1_{R}$ is integrable
by overconvergent functions, it is sufficient to determine the vanishing of 
the corresponding cohomology class in $H^1_\rig(U_k/K)$.  We have isomorphisms
$$
	H^1_\rig(U_k/K) = H^1_\dR(U_K) \cong H^1_\dR(E_K) = K \ul\omega \oplus K \uomegas
$$
for this cohomology group.

\begin{proof}[Proof of Theorem \ref{thm: existence one}]
	The statement is true if $m=0$ or $n=0$ by the previous lemma.  We define $D^{(p)}_{m,n}$ by 
	induction on $\theN = m+n$.  Suppose $\theN \geq 2$ and $D^{(p)}_{a,b}$ exists for integers $a$, $b$ such that $a + b < \theN$.
	Let $m$ and $n$ be integers $\geq 1$ such that $m + n = \theN$.  Then since $D^{(p)}_{m-1,n}$ and $D^{(p)}_{m,n-1}$
	are overconvergent functions,  the overconvergent differential form
	$$
		- D^{(p)}_{m-1,n} \omega - D^{(p)}_{m,n-1} \omegas
	$$
	defines a cohomology class in $H^1_\rig(U_k/K)$.  Since 
	$
		H^1_\rig(U_k/K)  = K \ul\omega \oplus K \uomegas,
	$
	there exists uniquely determined constants $c_{m,n}$, $c^*_{m,n} \in K$ such that the 
	cohomology class of
	$$
		- D^{(p)}_{m-1,n} \omega - D^{(p)}_{m,n-1} \omegas + c_{m,n} \omega + c^*_{m,n} \omegas
	$$
	vanishes in  $H^1_\rig(U_k/K)$.  This implies that for any $m$ and $n$ such that $m + n = \theN$,
	there exist overconvergent functions $\wt D_{m,n}$ on $\cU_K$ such that
	$$
		d \wt D_{m,n} = - D^{(p)}_{m-1,n} \omega - D^{(p)}_{m,n-1} \omegas + c_{m,n} \omega + c^*_{m,n} \omegas.
	$$
	We define $\wt G_{m,n}$ as in \eqref{eq: E}, with the highest term 
	$D^{(p)}_{m, n}$ replaced by $\wt D_{m,n}$.
	This function is again overconvergent, satisfying
	\begin{multline*}
			d \wt G_{m,n} = \sum_{k=0}^{n-1} \frac{(F_1^{(p)})^{n-k-1}}{(n-k-1)!} D^{(p)}_{m,k} 
			\left(1 - \frac{[\pi]^*}{\ol\pi} \right) \omega^*  
			- \sum_{k=0}^n \frac{(F_1^{(p)})^{n-k}}{(n-k)!} D^{(p)}_{m-1,k} \omega\\
			- \sum_{k=1}^n \frac{(F_1^{(p)})^{n-k}}{(n-k)!} D^{(p)}_{m,k-1} \omega^*  
			+ c_{m,n} \omega + c^*_{m,n} \omega^*\\
		= - \ol\pi^{-1} G^{(p)}_{m,n-1}  [\pi]^*\omega^*
		- G^{(p)}_{m-1,n} \omega
		+ c_{m,n} \omega + c^*_{m,n} \omega^*.
	\end{multline*}
	Let $C_{m,n} :=  \sum_{z_1 \in E[\frp]} \tau^*_{z_1} \left(\wt G_{m,n}\right)$.
	The differential forms $\omega$ and $[\pi]^* \omega^*$ are invariant under translations $\tau_{z_1} \colon z \mapsto z + z_1$
	for $z_1 \in E[\frp]$. Hence if we sum the above equation with respect to translations $\tau_{z_1}$,
	then the distribution relation for $G^{(p)}_{m,n-1}$ and $G^{(p)}_{m-1,n}$ give the relation
	$$
		dC_{m,n} = 
			\sum_{z_1 \in E[\frp]} \tau^*_{z_1} (d \wt G_{m,n})
			= \bN c_{m,n} \omega + 	c^*_{m,n} \sum_{z_1 \in E[\frp]} \tau_{z_1}^* (\omega^*).
	$$
	Since the cohomology class of $\tau_{z_1}^* (\omega^*)$ is $\ul\omega^*$ in 
	$H^1_\rig(U_k/K)$, and since $C_{m,n}$ is an overconvergent function, the above formula 
	implies that the cohomology class 
	$\bN ( c_{m,n} \ul\omega +  c^*_{m,n} \ul\omega^* )$ is zero in 
	$H^1_\rig(U_k/K)$.  Since $\ul\omega$ and $\ul\omega^*$ form
	a $K$-basis of $H^1_\rig(U_k/K)$, this implies that $c_{m,n} = c_{m,n}^* = 0$.
	As a consequence, we see that the function $C_{m,n}$ is  constant.   We let
	$$
		D^{(p)}_{m,n} := \wt D_{m,n} - (C_{m,n}/N(\frp)),
	$$ 
	which satisfies the differential equation \eqref{eq: differential elliptic polylog}.
	If we define $G^{(p)}_{m,n}$ as in \eqref{eq: E}, then we have 
	$G^{(p)}_{m,n} = \wt G_{m,n} - (C_{m,n}/\bN)$.
	Then by the definition of $C_{m,n}$, the function $G^{(p)}_{m,n}$ satisfies the distribution 
	relation \eqref{eq: distribution}.   The uniqueness of $D^{(p)}_{m,n}$ is assured by the distribution relation.
\end{proof}

\begin{proof}[Proof of Theorem \ref{thm: existence two}]
	The theorem follows from Theorem \ref{thm: existence one}, except the uniqueness.  	
	Suppose $D^{(p)}_{m,n}$ and $\wt D^{(p)}_{m,n}$ are two system of solutions satisfying 
	\eqref{eq: differential elliptic polylog}.  We prove by induction on $\theN = m+n$
	that $D^{(p)}_{m,n} = \wt D^{(p)}_{m,n}$. 
	Assume that $D^{(p)}_{a,b} = \wt D^{(p)}_{a,b}$ for $a+b < \theN$.	 
	We have by assumption
	$$
		d D^{(p)}_{m,n} =  - D^{(p)}_{m-1,n} \omega - D^{(p)}_{m,n-1} \omegas = d \wt D^{(p)}_{m,n}.
	$$
	for integers $m$, $n \geq 1$ satisfying $m+n = \theN$.
	Hence there exist constants $c_{m,n}$ such that $D^{(p)}_{m,n} = \wt D^{(p)}_{m,n} + c_{m,n}$.
	We extend the definition of this constant to $n=0$ by taking $\wt D^{(p)}_{\theN,0} = D^{(p)}_{\theN,0}$ and $c_{\theN,0}=0$.
	Then again by \eqref{eq: differential elliptic polylog} and the induction hypothesis, we have
	\begin{multline*}
		d \left( D^{(p)}_{m+1,n} - \wt D^{(p)}_{m+1,n} \right) = (D^{(p)}_{m,n} - \wt D^{(p)}_{m,n} ) \omega
		+ (D^{(p)}_{m+1,n-1} - \wt D^{(p)}_{m+1,n-1}) \omega^*
		\\= c_{m,n} \omega + c_{m+1,n-1} \omega^*.
	\end{multline*}
	This implies that the class of $c_{m,n} \omega + c_{m+1,n-1} \omega^*$ in 
	$H^1_\rig(U_k/K) = K \ul\omega \oplus K\uomegas$ is zero,
	hence we have $c_{m,n} = 0$ as desired.  Our assertion now follows by induction.
\end{proof}

%
\subsection{$p$-adic Eisenstein-Kronecker numbers and the elliptic polylogarithm functions}
%

We next compare the $p$-adic elliptic polylogarithm function $D^{(p)}_{m,n}$ constructed in the previous 
section with the $p$-adic distribution used in defining the $p$-adic Eisenstein-Kronecker numbers. 
We first begin by describing the residue discs of $E$.
Let $E^\an_{\bbC_p}$ be the extension to $\bbC_p$ of the rigid analytic space $E_K^\an$.  There is a morphism
$$
	\red : E_{\bbC_p}^\an \rightarrow E_{\ol\bbF_p}
$$
called the reduction map.  The inverse image of a point in $E_{\ol\bbF_p}$ is called a residue disc, 
which is an admissible open in $E^\an_{\bbC_p}$ (See [BGR, Sec. 9.1.4, Prop. 5]).
The formal parameter $s = -2x/y$ at the identity  of the elliptic curve parameterizes the residue disc 
around $[0]$, and we have a natural inclusion
$$
	\iota : B^{-}(0,1) \hookrightarrow E_{\bbC_p}^\an,
$$
where $B^-(0,1)$ is the rigid analytic open disc $B^-(0,1) := \{s \in \bbC_p \mid |s| < 1 \}$.

\begin{definition}
	Suppose we are given a rigid analytic function $f(z)$ on $E_{\bbC_p}^\an$. 
	 Then following the convention 
	of the complex case (see Remark \ref{rem: expansion}), we denote by $f(z)|_{z=\lambda(s)}$
	the rigid analytic function on $B^{-}(0,1)$ obtained as the pull-back of $f(z)$ by $\iota$.
\end{definition}

We denote by $\bsD$ the translation invariant derivation on $E$ defined by $df = \bsD(f) \omega$ for any 
overconvergent function $f$.  Then $\bsD$ restricts to the derivation
\begin{equation}\label{eq: operator}
	 \partial_{s, \log} := \lambda'(s)^{-1} \partial_s
\end{equation}
on each residue disc.

We define a $p$-modified variant of the generating function $F_{z_0, b}(z)$ using the 
right hand side of Proposition \ref{prop: distribution}.

\begin{definition} 
	For any integer $b \geq 0$, we define $\GenFunc_{z_0, b}^{(p)}(z)$ to be the function
	$$
		\GenFunc_{z_0, b}^{(p)}(z) = F_{z_0,b}(z)  - \ol\pi^{-b} F_{\pi z_0, b}(\pi z).
	$$
\end{definition}

Note that we have
$$
	\Theta^{(p)}_{z_0,0}(z,w) = \sum_{b \geq 0} \GenFunc_{z_0, b}^{(p)}(z) w^{b-1}.
$$
The relation $\Theta^{(p)}_{z_0,0}(z,w) =\exp \left(F_{z_0,1}(z)w\right) \Xi^{(p)}_{z_0}(z,w) $
of Lemma \ref{lem: L p} gives the relation
\begin{equation}\label{eq: F and L p}
	F^{(p)}_{z_0,b}(z) = 	\sum_{n=0}^b \frac{F_{z_0,1}(z)^{b-n}}{(b-n)!} L^{(p)}_n(z + z_0)
\end{equation}
between $F^{(p)}_{z_0,b}(z)$ and $L^{(p)}_{n}(z)$.
When $b=1$, then we have $F^{(p)}_{z_0, 1}(z) = F^{(p)}_1(z+z_0)$, which corresponds to
a rational function defined over $\bsK(z_0)$.  Assume now that $z_0$ is a torsion point in 
$E(\bbC_p)$ of order prime to $\frp$.
Then Proposition \ref{prop: distribution} implies that we have
$$
	\wh F^{(p)}_{z_0,b}(s) = \int_{\cO_{\bsKfrp}^\times}
	 \exp\left(x \lambda(s)/\Omega_\frp\right) d \mu_{z_0,b}(x),
$$
where $\wh F^{(p)}_{z_0,b}(s):= F^{(p)}_{z_0,b}(z)|_{z=\lambda(s)}$. 

\begin{definition}
	Suppose $z_0$ is a torsion point of $E(\bbC_p)$ of order prime to $\frp$.
	For any integer $m$ and $b$ such that $b \geq 0$, we define $\wh E^{(p)}_{z_0, m, b}(s)$ to be the function on $B^{-}(0,1)$
	given by the power series
	$$
		\wh E^{(p)}_{z_0, m, b}(s)  :=  (-\Omega_\frp)^{m} \int_{\cO_{\bsKfrp}^\times} x^{-m}
		 \exp\left(x \lambda(s)/\Omega_\frp\right)  d \mu_{z_0, b}(x).
	$$
\end{definition}

By the definition of $p$-adic Eisenstein-Kronecker numbers, we have
$$
 	\wh E^{(p)}_{z_0, m, b}(0) = \Omega_\frp^{m-b+1} e^{(p)}_{-m,b}(z_0)
$$
for any integer $m$ and $b$ such that $b \geq 0$.  
The purpose of this section 
is to establish a relation between $\wh E^{(p)}_{z_0, m, b}$ and $D^{(p)}_{m,n}$.
By construction, the function $\wh E^{(p)}_{z_0, m, b}(s)$ satisfies the differential equation
$$
	\partial_{s, \log} \wh E^{(p)}_{z_0, m, b}(s) = - \wh E^{(p)}_{z_0, m-1, b}(s).
$$
 Since  $\wh E^{(p)}_{z_0,m,b}(s)$ is defined as integration on $\cO_{\bsKfrp}^\times$, the part corresponding 
 to integration on $\frp\cO_{\bsKfrp}$ must be zero.  By calculating the restriction of distributions
on $\cO_{\bsKfrp}$ to $\frp\cO_{\bsKfrp}$, we obtain the  distribution relation
$$
	\sum_{s_1\in E[\frp]} \wh E^{(p)}_{z_0,m,b}(s\oplus s_1) = 0.
$$ 
We will use this property to characterize the function $\wh E^{(p)}_{z_0,m,b}$.  

\begin{lemma}\label{lem: characterization}
	Suppose $\wh F(s)$ is a function on $B^-(0,1)$ given by a power series in $K[[s]]$,
	satisfying the differential equation
	$
		\partial_{s, \log} \wh F(s) = \wh E_{z_0,m,n}(s)
	$
	and the distribution relation
	$$
		\sum_{s_1\in E[\frp]} \wh F(s\oplus s_1) = 0.
	$$ 
	Then we have $\wh F(s) = \wh E^{(p)}_{z_0,m+1,b}(s)$.
\end{lemma}

\begin{proof}
	Since 
	$
		\partial_{s, \log} \wh F(s) = \wh E_{z_0,m,n}(s) = \partial_{s,\log}  \wh E_{z_0,m+1,n}(s),
	$
	we see that $\wh F(s) = \wh E^{(p)}_{z_0,m+1,b}(s) + c$ for some constant $c \in K$.
	The distribution relation for $\wh F(s)$ and $\wh E^{(p)}_{z_0,m+1,b}(s)$ shows that
	 $\bN c = 0$, proving our assertion.
\end{proof}

The relation between $\wh E^{(p)}_{z_0,m,b}$ and $D^{(p)}_{m,n}$ may now be
given as follows.

\begin{proposition}
	For any integer $m \geq 0$, we have the equality
	\begin{equation}\label{eq: D explicit}
		\wh E^{(p)}_{z_0,m,b}(s)
		 =  \sum_{n=0}^b \frac{ \wh F_{z_0, 1}(s)^{b-n}}{(b-n)!} \left. D^{(p)}_{m,n}(z + z_0)
		 \right.|_{z=\lambda(s)}.
	\end{equation}
\end{proposition}

\begin{proof}
	We let $\wt E_{z_0,m,b}(s)$ be the right hand side of \eqref{eq: D explicit}.
	We prove the statement by induction on $m \geq 0$.
	When $m=0$, then $D^{(p)}_{0,n} = L^{(p)}_n$, and the relation between
	$L^{(p)}_n(z)$ and $F^{(p)}_{z_0,b}(z)$ of \eqref{eq: F and L p} gives the relation
	$$
		\wt E_{z_0,0,b}(s) =   \sum_{n=0}^b
		 \frac{ \wh F_{z_0, 1}(z)^{b-n}}{(b-n)!}\left.L^{(p)}_n(z + z_0)\right|_{z = \lambda(s)}
		 =  \wh F^{(p)}_{z_0,b}(s).
	$$
	This proves our assertion for $m=0$, since
	$
		\wh E^{(p)}_{z_0,0,b}(s) = \wh F^{(p)}_{z_0,b}(s)
	$
	by the definition of $\wh E^{(p)}_{z_0,0,b}(s)$.  	
	Suppose now that the statement is true for $m \geq 0$.  
	Let $\wh D^{(p)}_{z_0,m,n}(s) = D^{(p)}_{m,n}(z+z_0)|_{z=\lambda(s)}$.
	Note that since $F_{z_0,1}(z) = F_1(z+z_0) - \ol z_0/A$, we have
	$d F_{z_0,1} = d (\tau^*_{z_0}  F_1) = \tau^*_{z_0}(\omega^*)$.
	By the differential equation satisfied by $D^{(p)}_{m+1,n}$,
	we see that
	\begin{multline*}
		d \wt E_{z_0,m+1,b}(s) = 
		\sum_{n=0}^{b-1} \frac{ \wh F_{z_0, 1}(s)^{b-1-n}}{(b-1-n)!} 
		\wh D^{(p)}_{z_0,m+1,n}(s)  \tau^*_{z_0}(\omega^*)
		\\
		- 	\sum_{n=0}^b \frac{ \wh F_{z_0, 1}(s)^{b-n}}{(b-n)!} \wh D^{(p)}_{z_0,m,n}(s) \omega
		-    	\sum_{n=1}^b \frac{ \wh F_{z_0, 1}(s)^{b-n}}{(b-n)!} \wh D^{(p)}_{z_0,m+1,n-1} (s)
		\tau^*_{z_0}(\omega^* )\\
		= - \wt E_{z_0,m,b}(s) \omega = - \wh E^{(p)}_{z_0,m,b}(s) \omega.
	\end{multline*}
	Hence $\partial_{s,\log} \wt E_{z_0,m+1,b}(s)  = - \wh E^{(p)}_{z_0,m,b}(s)$.
	Furthermore, we have from the definition of $G^{(p)}_{m+1,k}$ and \eqref{eq: funny} that
	$$
		\wt E_{z_0,m+1,b}(s) =    \sum_{k=0}^b
		 \frac{ \wh F_{\pi z_0, 1}([\pi]s)^{b-k}}{\ol\pi^{b-k}(b-k)!}\left.
		 G^{(p)}_{m+1,k}(z+z_0)\right|_{z = \lambda(s)}.
	$$
	Note that we have $\wh F_{\pi z_0,1}([\pi](s \oplus s_1)) = \wh F_{\pi z_0,1}([\pi]s)$.
	Hence the distribution relation for $G^{(p)}_{m+1,k}$ gives the distribution relation
	$$
		\sum_{s_1\in E[\frp]} 	\wt E_{z_0,m+1,b}(s \oplus s_1) = 0
	$$
	for $\wt E_{z_0,m,b}(s)$.   Our assertion now follows from Lemma \ref{lem: characterization}.
\end{proof}

%
%
%
\section{$p$-adic realization of the elliptic polylogarithm}
%
%
%

We keep the notation of \S \ref{ss: EK number}.  In this section, we explicitly determine the $p$-adic
elliptic polylogarithm sheaf.  In particular, we will show that the functions $\pPolFunc_{m,n}$
give the Frobenius action of the $p$-adic elliptic polylogarithm sheaf.

%
\subsection{Rigid syntomic cohomology}
%

We first briefly recall the theory of filtered overconvergent  $F$-isocrystals (or syntomic coefficients) 
and rigid syntomic cohomology developed in \cite{Ba1}.  Let $K$ be a finite unramified extension of $\bbQ_p$ 
with ring of integers $\cO_K$ and residue field $k$.  We fix an integer $q = p^m$ for some $m\geq1$, 
and we denote by $\Frob_q$ the Frobenius $x \mapsto x^q$ on $k$.  We let $\sigma$ be the 
extension to $\cO_K$ and $K$ of the Frobenius $\Frob_q$ on $k$. 
 
Let $X$ be a smooth scheme of finite type over $\cO_K$, with smooth compactification $j: X \hookrightarrow \ol X$ 
over $\cO_K$ such that the complement $D:= \ol X \setminus X$ is a relative strict normal crossing divisor over 
$\cO_K$.  Denote by $\cX$ and $\ol\cX$ the formal completion of $X$ and $\ol X$ with respect to the special fiber.
We assume in addition that there exists a Frobenius $\phi: \ol\cX \rightarrow \ol\cX$ lifting the 
Frobenius $\Frob_q$ on $\ol X_k := \ol X \otimes k$, such that $\phi(\cX) \subset \cX$.  Then the triple
$$
	\sX := (X, \ol X, \phi)
$$
is a syntomic datum in the sense of \cite{Ba1} Definition 1.1.  

Suppose we are given a coherent $\cO_{\ol X_K}$-module $\cM$ with integrable connection 
$\nabla : \cM \rightarrow \cM \otimes \Omega^1_{\ol X_K}(\log D)$ with logarithmic poles along $D$.
Then as in \cite{Ba1}, paragraph before Definition 1.8, we let
$$
	\sM_\rig := j^\dagger \cM^\an
$$
with connection 
$$
	\nabla_\rig: \sM_\rig \rightarrow \sM_\rig \otimes \Omega^1_{\ol\cX_K}.
$$ 
A Frobenius structure on $\sM_\rig$ is an isomorphism
$\Phi: \phi^*(\sM_\rig) \xrightarrow\cong \sM_{\rig}$ of $j^\dagger\cO_{\ol\cX_K}$-modules
compatible with the connection.  If such a structure exists on $\sM_\rig$, then by \cite{Ber1} Theorem 2.5.7,
the connection $\nabla_\rig$ on $\sM_\rig$ becomes overconvergent (in the sense of
\cite{Ber1} Definition 2.2.5), and the pair $(\sM_\rig, \Phi)$ gives  a realization
(in the sense of \cite{Ber1} p.68) of an overconvergent $F$- isocrystal in
$F$-$\Isoc^\dagger(X_k/K)$.

\begin{definition}
	We define the category of filtered overconvergent $F$- isocrystals on $\sX$ to be the category $S(\sX)$
	consisting of objects the 4-uple $\sM:= (\cM, \nabla, F^\bullet, \Phi)$ consisting of:
	\begin{enumerate}
		\item $\cM$ is a coherent $\cO_{\ol X_K}$-module of finite rank.
		\item $\nabla : \cM \rightarrow \cM \otimes \Omega^1_{\ol X_K}(\log D)$ is an integrable connection
		on $\cM$ with logarithmic singularities along $D$.
		\item $F^\bullet$ is a descending exhaustive separated filtration by coherent
		 $\cO_{\ol X_K}$-submodules on $\cM$, satisfying the Griffiths transversality
		$$
			\nabla(F^m \cM) \subset  F^{m-1} \cM \otimes \Omega^1_{\ol X_K}(\log D).
		$$
		\item $\Phi: \phi^*(\sM_\rig) \xrightarrow\cong \sM_{\rig}$ is a Frobenius structure on $\sM_\rig$.
	\end{enumerate}
\end{definition}

An important example of filtered overconvergent $F$-isocrystals is given by the Tate objects.

\begin{definition}
	We define the \textit{Tate object} in $S(\sX)$
	to be the filtered overconvergent $F$-isocrystal $K(j) = (\cM, \nabla, F^\bullet, \Phi)$,
	such that $\cM = \cO_{\ol X_K} e_j$, the connection is given by $\nabla(e_j) = 0$, the Hodge filtration is
	such that $F^{-j} \cM = \cM$ and $F^{-j+1} \cM = 0$, and the Frobenius is $\Phi(e_j) = p^{-j} e_j$.
\end{definition}

For any filtered overconvergent $F$-isocrystal $\sM$ in $S(\sX)$, we may define the de Rham and
rigid cohomology $H^i_\dR(\sX, \cM)$ and $H^i_\rig(\sX,\sM_\rig)$ of $\sX$ with coefficients in $\cM$ and $\sM_\rig$.
The de Rham cohomology has a Hodge filtration induced from the Hodge filtration on $M$, and
rigid cohomology has a Frobenius $\phi: H^i_\rig(\sX,\sM_\rig) \otimes_\sigma K \xrightarrow\cong H^i_\rig(\sX,\sM_\rig)$
induced from the Frobenius  on $\sM_\rig$.   There exists a natural homomorphism
\begin{equation}\label{eq: theta}
	\theta: H^i_\dR(\sX, \cM) \rightarrow H^i_\rig(\sX,\sM_\rig)
\end{equation}
(See for example \cite{BC} \S 2 or \cite{Sh}).

\begin{remark}\label{rem: admissible}
	In order to define rigid syntomic cohomology, we consider filtered overconvergent $F$- isocrystals $\sM$ on $\sX$ 
	satisfying the following additional conditions.
	\begin{enumerate}
		\item The Hodge to de Rham spectral sequence
		$$
			E_1^{i,j} = R \Gamma^{i+j} ( \ol X_K, \Gr_F^i( \cM \otimes \Omega^\bullet_{\ol X_K}(\log D) ))
			\Rightarrow H^{i+j}_\dR(\sX, \cM)
		$$
		degenerates at $E_1$.
		\item The $\theta$ of \eqref{eq: theta} is an isomorphism.
	\end{enumerate}
\end{remark}

In what follows, suppose $\sM$ satisfies the condition of Remark \ref{rem: admissible}.

\begin{definition}\label{def: p-adic cohomology}
	We define the filtered Frobenius module $H^i(\sX,\sM)$ of $\sX$ with coefficients in $\sM$ to be the cohomology group
	$H^i_\rig(\sX,\sM_\rig)$ with its natural Frobenius and Hodge filtration induced 
	from the Hodge filtration on de  Rham cohomology $H^i_\dR(\sX, \cM)$ through the isomorphism $\theta$.
\end{definition}

\begin{remark}
	In our previous paper, we imposed an additional condition, that the filtered Frobenius module $H^i(\sX,\sM)$
	defined above is a weakly admissible filtered Frobenius module in the sense of Fontaine (\cite{Fon} 4.1.4).  
	This would insure that morphisms between $H^i(\sX,\sM)$ are strictly compatible with the Hodge filtration.
	However, in the current paper, especially in the supersingular case, we are considering Frobenius $\sigma$ which is 
	not necessarily absolute, hence the notion of weakly admissible is not useful.  Hence morphisms between $H^i(\sX,\sM)$
	may not a priori be strictly compatible with the Hodge filtration.  For our application, any such morphism we use
	will be strictly  compatible with the Hodge filtration, since it underlies a morphism of mixed Hodge structures.
\end{remark}

We denote by $H^i_\syn(\sX, \sM)$ the rigid syntomic cohomology (or simply syntomic cohomology)
of $\sX$ with coefficients in $\sM$, defined in \cite{Ba1} Definition 2.4.  The most important fact concerning 
syntomic cohomology that we will use in this paper is the following.

\begin{proposition}[\cite{Ba1} Theorem 1]\label{prop: ext isom}
	We have a canonical isomorphism
	$$
		\Ext^i_{S(\sX)}(K(0), \sM) \xrightarrow\cong H^i_{\syn}(\sX, \sM)
	$$
	when $i=0$, $1$, and $K(0)$ is the Tate object in $S(\sX)$.
\end{proposition}

Denote by $\sV = (\Spec\, \cO_K, \Spec\, \cO_K, \sigma)$ the trivial syntomic datum.
In this case, any object $\sM$ in $S(\sV)$ is simply a filtered Frobenius module, and we have
\begin{equation}\label{eq: syn base}
	H^i_{\syn}(\sV, \sM) = H^i \left[ F^0 \cM \xrightarrow{1-\phi^*} \cM \right],
\end{equation}
where the $F^0\cM$ in the complex on the right is in degree zero.
We have the following relation between rigid cohomology and syntomic cohomology.

\begin{lemma}\label{lem: SES}
	We have the short exact sequence
	\begin{multline*}
		0 \rightarrow H^1_\syn(\sV, H^i(\sX, \sM)) \rightarrow H^{i+1}_\syn(\sX, \sM) \\
		\rightarrow H^0_\syn(\sV, H^{i+1}(\sX, \sM)) \rightarrow 0.
	\end{multline*}
\end{lemma}
We have an inclusion 
$H^0_\syn(\sV, H^1(\sX, \sM)) \hookrightarrow H^1_\rig(\sX, \sM_\rig)$.  Let
$$
	 H^{1}_\syn(\sX, \sM) \rightarrow
	 H^{1}_\rig(\sX, \sM_\rig) \overset{\theta^{-1}}{\cong} H^1_\dR(\sX, \cM) \rightarrow H^1_\dR(X_K, \cM)
$$ 
be the morphism induced from the surjection of the short exact sequence in Lemma \ref{lem: SES}.
The relation between this map and the isomorphism of Proposition \ref{prop: ext isom} is given by the following.

\begin{lemma}[\cite{Ba1} Proposition 4.4]\label{lem: compatible}
	The following diagram is commutative
	$$
		\begin{CD}
				\Ext^1_{S(\sX)}(K(0), \sM)  @>\For>> 	\Ext^1_{M(X_K)}(K(0), \cM)  \\
				@V{\cong}VV   @V{\cong}VV\\
				 H^{1}_\syn(\sX, \sM)  @>>> 	 H^{1}_\dR(X_K, \cM),
		\end{CD}
	$$
	where $M(X_K)$ is the category of coherent modules with integrable connection on $X_K$, and 
	$\For: S(\sX) \rightarrow M(X_K)$ is the functor obtained by forgetting the Hodge filtration and the Frobenius structure.
\end{lemma}

%
\subsection{The Logarithm sheaf}
%

We now describe the logarithm sheaf.  We fix a model $E$ over $\cO_{\bsK}$
of our elliptic curve as in \eqref{eq: integral Weierstrass}.  We again let $K$ be a finite unramified extension of
$\bsKfrp$ in $\bbC_p$.  Let $\pi := \psi_{E/\bsK}(\frp)$ as before, and let $[\pi]: E \rightarrow E$ be the multiplication
induced on $E$.  We define $\phi: = [\pi] \otimes \sigma$ to be the Frobenius on  
$E_{\cO_K} := E \otimes_{\cO_{\bsK}} \cO_K$ induced from $[\pi]$ and the Frobenius on $\cO_K$.  
Then $\phi$  is a Frobenius morphism of degree $[K_\frp: \bbQ_p]$, and the triple $\sE = (E_{\cO_K}, E_{\cO_K}, \phi)$ is a syntomic 
datum.   

By Damerell's theorem, we have $e_2^* \in \bsK \subset K$.   We let the notations be as in 
Remark \ref{rem: choice of basis}.  We let $H^1(\sE)$ the filtered Frobenius module defined in
Definition \ref{def: p-adic cohomology}  associated
to the trivial object in $S(\sE)$.  Then $H^1(\sE)$ is a $K$-vector space of rank two
 $H^1(\sE) = K \ul\omega \oplus K\uomegas$, with Hodge filtration such that
 $F^{0} H^1(\sE)  = H^1(\sE)$, $F^1 H^1(\sE)  =  K \ul\omega$, and $F^2  H^1(\sE)  =0$.
The action of the Frobenius $\phi^*:  H^1(\sE) \rightarrow  H^1(\sE)$ is given by
\begin{align*}
	\phi^*(\ul\omega) &= \pi \ul\omega, & 	\phi^*(\uomegas) &= \ol\pi\uomegas.
\end{align*}
Hence the action on  $\sH := H^1(\sE)^\vee$ is given by
\begin{align*}
		\phi^*(\uomegav) &= \pi^{-1}\uomegav, & \phi^*(\uomegasv) &= \ol\pi^{-1} \uomegasv.
\end{align*}
For any syntomic datum $\sX$, we denote by $\sH_\sX$ the constant filtered overconvergent $F$-isocrystal on $\sX$
obtained as the pull-back of $\sH$ to $\sX$.
We have the natural isomorphism
$$
	H^0_\syn(\sV, H^1(\sE, \sH_\sE)) = H^0_\syn(\sV, \sH^\vee \otimes \sH) = \Hom_{S(\sV)}(\sH, \sH).
$$
The short exact sequence of Lemma \ref{lem: SES} in this case gives
\begin{equation}\label{eq: p-adic deg}
	0 \rightarrow H^1_\syn(\sV, \sH) \rightarrow H^1_\syn(\sE, \sH_\sE) \rightarrow
	\Hom_{S(\sV)}( \sH, \sH) \rightarrow 0.
\end{equation}
The pull back  $i_{[0]}^*$  by the identity $[0]: \sV \rightarrow \sE$ of the elliptic curve gives a splitting
$i_{[0]}^*:  H^1_\syn(\sE, \sH_\sE) \rightarrow H^1_\syn(\sV, \sH)$ of the above exact sequence.

\begin{definition}\label{def: p-adic log}
	We define the first logarithm sheaf  $\sLog^{(1)}$ to be any extension of $\sH_\sE$ by $K(0)$ in $S(\sE)$,
	whose extension class in
	$$
		\Ext^1_{S(\sE)}(K(0), \sH_\sE) \cong H^1_\syn(\sE, \sH_\sE)
	$$
	is mapped by the surjection of \eqref{eq: p-adic deg} to the identity and to zero by the splitting $i_{[0]}^*$.
	We define the $N$-th logarithm sheaf $\sLog^\theN$ to be the $N$-th symmetric product of $\sLog^{(1)}$.
\end{definition}

Next, we explicitly construct $\sLog^{(1)}$ in the $p$-adic case.  We take an affine open covering 
$\frU = \{ U_i \}_{i \in I}$ of $E$, and we fix $(\eta_i, u_{ij})$ and $(u_i)$ as in Definition \ref{def: basis}.
By abuse of notation, we again denote by $U_i$ the base extension $U_i \otimes_{\cO_\bsK} \cO_K$ of
$U_i$ to $\cO_K$. Then $\{ U_i \}_{i \in I}$ is an affine covering of $E_{\cO_K}$. 
We denote by $\cLog^{(1)}$ the module defined in Proposition \ref{pro: log dR},
whose restriction to $U_i$ is given by the module $\cLog_i^{(1)}$ defined as
$$
		\cLog^{(1)}_i := \cO_{U_i} \ul e_i \bigoplus \cH_{U_i},
$$
with Hodge filtration as the direct sum and connection $\nabla(\ul e_i) =\uomegav \otimes \omega + \uomegasv \otimes \omegas_i$.

We define the Frobenius $\Phi: \phi^*(\sLog_\rig^{(1)}) \rightarrow \sLog_\rig^{(1)}$ as follows.  
The connection of $\phi^*\cLog^{(1)}_i$ on $\phi^{-1}(U_i)$ is given by
$$
	\nabla_{\kern-0.5mm\phi}(\ul e^\phi_{i}) = 
	\uomegav \otimes [\pi]^*\omega+  \uomegasv \otimes [\pi]^*\omegas_i.
$$
We let  
$$
	F_1^{(p)} := \left(  1 - \frac{[\pi^*]}{\ol\pi} \right) F_1 = \frac{1}{\bN} \left( 
	\log \left( \theta(z)^{\bN} / \theta(\pi z) \right) \right)'
$$
as in the previous section.
We define  $\xi^{(p)}_{ij}$ to be the rational function
$\xi^{(p)}_{ij} := F_1^{(p)} - u_j + ([\pi]^*u_i/\ol\pi)$  on  $\phi^{-1}(U_i) \cap U_j$.  
Then the differential $d \xi^{(p)}_{ij} = \omega_j^* -  [\pi]^* \omega^*_i/\ol \pi$ is holomorphic on  $\phi^{-1}(U_i) \cap U_j$,
hence $\xi^{(p)}_{ij}$ is also holomorphic on $\phi^{-1}(U_i) \cap U_j$.
We define the Frobenius $\Phi: \phi^*(\sLog_\rig^{(1)}) \rightarrow \sLog_\rig^{(1)}$ 
locally on $\phi^{-1}(U_i) \cap U_j$ as the morphism
$$
	\Phi_{ij} \colon\phi^*(\cLog^{(1)}_i) \xrightarrow\cong \cLog^{(1)}_j
$$
given by
\begin{align*}
	\Phi_{ij}(\ul e_i^\phi) &=\ul e_j  -\xi_{ij}^{(p)} \uomegasv, &
	\Phi_{ij}(\uomegav) &= \frac{1}{\pi} \uomegav, &
	\Phi_{ij}(\uomegasv) &= \frac{1}{\ol\pi} \uomegasv.
\end{align*}
Then $\Phi$ gives an isomorphism compatible with the connection.  

\begin{proposition}\label{pro: p split}
	We define $\sLog^{(1)}$ to be the filtered overconvergent $F$-isocrystal whose
	underlying coherent module with connection is given by $\cLog^{(1)}$, whose Hodge filtration
	is given as the direct sum of the Hodge filtration on each $U_i$, and the Frobenius
	$\Phi: \phi^*(\sLog_\rig^{(1)}) \rightarrow \sLog_\rig^{(1)}$ is given as above.
	Then this satisfies the property of the first logarithm sheaf.
\end{proposition}

\begin{proof}
	By construction, the underlying coherent module with connection of $\sLog^{(1)}$ maps
	to that of de Rham cohomology, which by definition corresponds to
	the identity.  Hence by Lemma \ref{lem: compatible}, the surjection of \eqref{eq: p-adic deg} 
	maps $\sLog^{(1)}$ to the identity.   It remains to show that $\sLog^{(1)}$ splits as an extension of 
	$\sH$ by $K(0)$ when pulled back by $i_{[0]}^*$.   From our choice of $u_i$, for $\zeta_i (z)= F_1(z) - u_i(z)$,
	we have $\zeta_i(0) = 0$.  Hence $\xi^{(p)}_{ij}(0) = 0$, since we have from the
	definition $\xi^{(p)}_{ij}(z) = \zeta_j(z) - [\pi]^* \zeta_i(z)/\ol\pi$.  
	This implies that the splitting $\varphi: i_{[0]}^* \cLog^{(1)} \cong K \bigoplus \cH$
	 given in \eqref{eq: split log one} of underlying $K$-vector spaces gives a splitting compatible with the Hodge
	 filtration and the Frobenius.
\end{proof}

Let $U_{\cO_K} = E_{\cO_K} \setminus [0]$ and $j : U_{\cO_K} \hookrightarrow E_{\cO_K}$ be the natural inclusion.
The Frobenius $\Phi_{ij}$ on each $\phi^{-1}(U_i) \cap U_j$ paste together to give the Frobenius
$$
	\Phi_\sL(\ul e^\phi) = \ul e - F_1^{(p)} \uomegasv
$$ 
for
$
 	\sL := j^\dagger \cO_{\cE_K} \ul e \bigoplus  j^\dagger \cO_{\cE_K} \uomegav
  	\bigoplus j^\dagger \cO_{\cE_K} \uomegasv
$
on $\cE_K$.  This gives the following.

\begin{corollary}
	The base extension of the overconvergent $F$-isocrystal underlying the logarithm sheaf $\sLog^\theN$ to 
	$j^\dagger \cO_{\cE_K}$ is given as follows.
	The underlying module is
	$$
		\cL^\theN = \bigoplus_{0 \leq m+n \leq \theN} j^\dagger \cO_{\cE_K} \ul\omega^{m,n},
	$$
	with connection given by $\nabla_{\kern-1mm\sL} = d + \boldsymbol{\nu}$ and the Frobenius isomorphism 
	given by
	$$
		\Phi_{\sL}(\ul\omega^{m,n}_{\phi}) = 
		\frac{1}{\pi^m \ol \pi^n} 
		\sum_{k=n}^{\theN-m}  \frac{(-\pGenFunc_1)^{k-n}}{(k-n)!} \ul\omega^{m, k}.
	$$
\end{corollary}
The above overconvergent $F$-isocrystal underlies the filtered overconvergent $F$-isocrystal 
obtained as the restriction of $\sLog^\theN$ to $\sU = (U_{\cO_K}, E_{\cO_K}, \phi)$.

%
\subsection{The Polylogarithm sheaf}
%

In this section, we define and calculate the $p$-adic polylogarithm sheaf.  We again let
$U_{\cO_K} = E_{\cO_K} \setminus [0]$ and we define $\sU$ to be the syntomic datum $\sU = (U_{\cO_K}, E_{\cO_K}, \phi)$.
We also let $\sD = ([0], [0], \sigma)$.
Similarly to the Hodge case, we have the following.

\begin{lemma}\label{lem: rig coh} 
	We have isomorphisms of filtered Frobenius modules
	\begin{align*}
		\varprojlim_\theN H^0(\sU, \sLog^\theN(1))
			&		\xleftarrow\cong
		 \varprojlim_\theN H^0(\sE, \sLog^\theN(1))  = 0 \\
		\res: \varprojlim_\theN H^1(\sU, \sLog^\theN(1)) & \xrightarrow\cong
		 \varprojlim_\theN H^0(\sD,  i_{[0]}^* \sLog^\theN)= \varprojlim_\theN i_{[0]}^* \sLog^\theN.
	\end{align*}
\end{lemma}

\begin{proof}
	By Tsuzuki \cite{Ts}, the localization maps in rigid cohomology are compatible with the Frobenius.
	The isomorphisms follow from the calculation for de Rham cohomology given in Lemma \ref{lem: log coh}.
	The Hodge filtration is strictly compatible, since the isomorphism for de Rham cohomology underlies an 
	isomorphism of mixed Hodge structures.
\end{proof}

The exact sequence of Lemma \ref{lem: SES} gives the short exact sequence
\begin{multline*}
	0 \rightarrow H^1_\syn(\sV, H^0(\sU, \sH_\sU^\vee \otimes \sLog^\theN)(1))
	\rightarrow H^1_\syn(\sU, \sH_\sU^\vee \otimes \sLog^\theN(1)) \\\rightarrow
	H^0_\syn(\sV, H^1(\sU, \sH_\sU^\vee \otimes \sLog^\theN)(1)) \rightarrow 0.
\end{multline*}
Since the projective limit of $H^0(\sU, \sLog^\theN)$ is zero, this gives a natural isomorphism
\begin{equation}\label{eq: p-pol isom one}
	\varprojlim_\theN
	H^1_\syn(\sU, \sH_\sU^\vee \otimes \sLog^\theN(1))  \xrightarrow\cong
	\varprojlim_\theN
	H^0_\syn(\sV, H^1(\sU, \sH_\sU^\vee \otimes \sLog^\theN(1))).
\end{equation}
We have a natural isomorphism
\begin{equation}\label{eq: p-pol isom two}
	H^0_\syn(\sV, H^0(\sD, \sH_\sE^\vee \otimes i_{[0]}^* \sLog^\theN))\\
	\xrightarrow\cong H^0_\syn(\sV, \sH^\vee \otimes i_{[0]}^* \sLog^\theN)
	\cong \bbQ_p,
\end{equation}
where the last isomorphism is obtained by mapping the identity element in 
$\sH^\vee \otimes \sH = \Hom(\sH, \sH)$ to $1 \in \bbQ_p$.  Combining \eqref{eq: p-pol isom one} 
and \eqref{eq: p-pol isom two} with the residue isomorphism, we obtain an isomorphism
\begin{equation}\label{eq: p-res isom two}
	\varprojlim_\theN H^1_\syn(\sU, \sH_\sU^\vee \otimes \sLog^\theN(1)) \xrightarrow\cong \bbQ_p.
\end{equation}

\begin{definition}\label{def: p-pol class}
	We define the polylogarithm class to be a system of classes $\pol^\theN \in  
	H^1_\syn(\sU, \sH_\sU^\vee \otimes \sLog^\theN(1))$
	such that
	$$
		\pol := \varprojlim_\theN \pol^\theN \in 	\varprojlim_\theN  H^1_\syn(\sU, \sH_\sU^\vee \otimes \sLog^\theN(1))
	$$
	maps through \eqref{eq: p-res isom two} to $1$.
\end{definition}

\begin{definition}
	We define the elliptic polylogarithm sheaf on $\sU$ to be a system of filtered overconvergent $F$-isocrystals
	$\sP^\theN$ on $\sU$ given as an extension
	$$
		0 \rightarrow \sLog^\theN(1) \rightarrow \sP^\theN \rightarrow \sH_\sU \rightarrow 0
	$$
	whose extension class corresponds to $\pol^\theN$ in 
	$$
		\Ext^1_{S(\sU)}(\sH_\sU, \sLog^\theN(1)) \cong H^1_\syn(\sU, \sH_\sU^\vee \otimes \sLog^\theN(1)).
	$$ 
\end{definition}

We now prove the main result of this paper, explicitly describing the $p$-adic elliptic polylogarithm sheaf
$\sP^\theN$ on $\sU$. 

\begin{theorem}\label{thm: describe polylog}
	The elliptic polylogarithm sheaf $\sP^\theN$ on $\sU$ is given by the filtered overconvergent $F$-isocrystal
	$\sP^\theN := (\cP^\theN, \nabla, F^\bullet, \Phi)$ defined as follows.
	\begin{enumerate}
		\item $\cP^\theN$ is the coherent module $\cP^\theN = \cH_E \bigoplus \cLog^\theN$, with integrable
			connection $\nabla(\ul\omega^\vee) = \bomega^\vee$ and $\nabla(\ul\omega^{*\vee}) = \bomegas^{\vee}$.
		\item $F^\bullet$ is the filtration on $\cP^\theN$ given by
			$
				F^m \cP^\theN = F^m \cH_E \bigoplus F^{m+1} \cLog^\theN.
			$
		\item We denote by $\Phi_{\!\sP}$ the morphism
			$
				\Phi_{\!\sP} : \phi^*\sP^\theN_\rig \rightarrow \sP^\theN_\rig
			$
			which extends the Frobenius $\Phi_{\sL(1)} := \bN^{-1} \Phi_\sL$ on $\sLog_\rig^\theN(1)$ and is 
			given by
			\begin{align*}
				\Phi_{\!\sP}(\pi \uomegav) &=  \uomegav 
				-\sum_{n=0}^\theN \frac{(- F_1^{(p)})^{n+1}}{(n+1)!}  \ul\omega^{0,n}
				+  \sum_{m=1}^\theN \sum_{n=0}^{\theN-m}  \pPolFunc_{m,n+1} \ul\omega^{m,n}\\
			\Phi_{\!\sP}(\ol\pi \uomegasv) &= \uomegasv +
			 \sum_{m=0}^\theN \sum_{n=0}^{\theN-m}\pPolFunc_{m+1,n} \ul\omega^{m,n}. 
			\end{align*}
	\end{enumerate}
	By definition, $\Phi_{\!\sP}$ is compatible with the projection $\sP^{\theN+1} \rightarrow \sP^\theN$.
\end{theorem}

\begin{proof}
	Since the polylogarithm sheaf as an extension class is uniquely characterized by the property of
	Definition \ref{def: p-pol class}, it is sufficient to prove that $\sP^\theN$ is a filtered overconvergent $F$-isocrystal
	satisfying the required property.
	By Lemma \ref{lem: compatible}, we may take the underlying module with connection 
	of the polylogarithm sheaf to be that of the de Rham realization of \S 1, and we may take the
	Hodge filtration as in \S 3, as the direct sum.   Hence it is sufficient to prove that $\Phi_{\!\sP}$
	is compatible with the connection and indeed gives a Frobenius structure on the
	overconvergent isocrystal underlying $\sP^\theN$.  This fact is 
	Proposition \ref{pro: F compatible} below.
\end{proof}

The rest of this subsection is devoted to proving the following proposition.
\begin{proposition}\label{pro: F compatible}
	Let the notations be as in Theorem \ref{thm: describe polylog}.
	Then the  morphism $ \Phi_{\!\sP}\colon\phi^*\sP_\rig^\theN \rightarrow \sP_\rig^\theN$ is compatible with the connection on 
	$\sP_\rig^\theN$.
\end{proposition}

The proof of Proposition \ref{pro: F compatible} will be given at the end of this section.   
We first start with a lemma concerning the action of $\Phi_{\!\sP}$ on 
$\bomega^\vee$ and $\bomega^{*\vee}$.  Recall that the Frobenius $\phi$ on $\cE$ is defined as 
$\phi:= [\pi] \otimes \sigma$.

\begin{lemma}
	We have
	\begin{align*}
		\left( 1 - \pi \Phi_{\sL(1)}\right) {\bomega}^{\vee} &= 
		- \ul\omega^{0,0} \otimes
		\left( 1 - \frac{[\pi]^*}{\ol\pi}\right) \omegas
		+ \sum_{k=1}^\theN \frac{(-\pGenFunc_1)^{k}}{k!}\ul\omega^{0,k} \otimes \frac{[\pi]^*}{\ol\pi}\omegas
		\\
		& \qquad +\sum_{k=0}^{\theN-1} \left( L^{(p)}_{k+1} +  \frac{(-\pGenFunc_1)^{k+1}}{(k+1)!}\right) \ul\omega^{1,k} 
			\otimes \omega,\\
		\left( 1 - \ol\pi\Phi_{\sL(1)}\right) {\bomega}^{*\vee} &= 
		\sum_{k=1}^\theN L^{(p)}_{k} \ul\omega^{0,k} \otimes \omega.
	\end{align*}
\end{lemma}

\begin{proof}
	The action of the Frobenius on ${\bomega}^\vee$ is given by
	\begin{multline*}
		\Phi_{\sL(1)}({\bomega}^\vee) 
		%
		%
		%
		%
		%
		%
	 	=-\frac{1}{\bN} \sum_{k=0}^\theN \frac{(-\pGenFunc_1)^{k}}{k!}\ul\omega^{0,k} \otimes [\pi]^*\omegas 
		\\+ \frac{1}{\pi}
		\sum_{k=0}^{\theN-1} \sum_{n=0}^k \frac{([\pi]^* L_{n+1})(-\pGenFunc_1)^{k-n}}{\ol\pi^{n+1}(k-n)!}
		 \ul\omega^{1,k} \otimes \omega.
	\end{multline*}
	The last sum may be expressed as 
	$$
		 \sum_{n=0}^k\frac{([\pi]^* L_{n+1})(-\pGenFunc_1)^{k-n}}{\ol\pi^{n+1}(k-n)!} 
		 = \left( \sum_{n=0}^{k+1}\frac{([\pi]^* L_{n})(-\pGenFunc_1)^{k+1-n}}{\ol\pi^{n}(k+1-n)!}  \right) -
		 \frac{(-\pGenFunc_1)^{k+1}}{(k+1)!}.
	$$
	Expanding the $L_n$ in the above sum, we have
	\begin{multline*}
		\sum_{n=0}^{k+1}\frac{([\pi]^* L_{n})(-\pGenFunc_1)^{k+1-n}}{\ol\pi^{n}(k+1-n)!} 
		%
		%
		 %
		 %
		 =\sum_{b=0}^{k+1} \sum_{n=b}^{k+1}  \frac{(-[\pi]^* F_1)^{n-b} (-F_1^{(p)})^{k+1-n}  }{\ol\pi^{n-b} (n-b)!(k+1-n)!} 
		 \frac{[\pi]^*F_b}{\ol\pi^b} \\
		 =  \sum_{b=0}^{k+1} \frac{(-F_1)^{k+1-b}}{(k+1-b)!} 	 \frac{[\pi]^*F_{b}}{\ol\pi^{b}}.
	\end{multline*}
	Hence
	$$
		L_{k+1} - \sum_{n=0}^{k} 
		\frac{([\pi]^* L_{n+1})(-\pGenFunc_1)^{k-n}}{\ol\pi^{n+1}(k-n)!}
		= L_{k+1}^{(p)} + \frac{(-\pGenFunc_1)^{k+1}}{(k+1)!}
	$$
	as desired.   The second equality may be proved in a similar fashion, again by direct 
	calculation.
%
%
%
\end{proof}

\begin{proof}[Proof of Proposition \ref{pro: F compatible}]
	For the basis $\pi\uomegav$,  if we calculate the composition of the connection with the Frobenius, 
	then we have
	$
		 \Phi_{\!\sP} \circ \nabla_\sP(\pi\uomegav) = \pi\Phi_{\sL(1)}({\bomega}^\vee)
	$
	and  
	\begin{multline}\label{eq: short calc}
		\nabla_\sP \circ \Phi_{\!\sP}(\pi\uomegav)\\
		=    \nabla_\sP
	 	\left( \uomegav - \sum_{n=0}^\theN \frac{(- F_1^{(p)})^{n+1}}{(n+1)!} \ul\omega^{0,n}
			+ \sum_{m=1}^\theN \sum_{n=0}^{\theN-m} \pPolFunc_{m,n+1} \ul\omega^{m,n} \right).
	\end{multline}
	Since $\omegas = d F_1$, we have
	\begin{multline*}
		 \sum_{n=0}^\theN \frac{(- F_1^{(p)})^{n}}{n!} \ul\omega^{0,n} \otimes d F_1^{(p)}
		- \sum_{n=0}^{\theN-1} \frac{(- F_1^{(p)})^{n+1}}{(n+1)!} \ul\omega^{0,n+1} \otimes \omegas \\
		= \ul\omega^{0,0} \otimes
		\left( 1 - \frac{[\pi]^*}{\ol\pi}\right) \omegas
		- \sum_{k=1}^\theN \frac{(-\pGenFunc_1)^{k}}{k!}\ul\omega^{0,k} \otimes \frac{[\pi]^*}{\ol\pi}\omegas,
	\end{multline*}
	Using the equality and the differential equation satisfied by $\pPolFunc_{m,n}$, we see that 
	\eqref{eq: short calc} is equal to
	\begin{multline*}
	 	{\bomega}^{\vee} 
		 +\sum_{n=0}^\theN \frac{(- F_1^{(p)})^{n}}{n!} \ul\omega^{0,n} \otimes d F_1^{(p)}
		- \sum_{n=0}^{ \theN-1} \frac{(- F_1^{(p)})^{n+1}}{(n+1)!} \ul\omega^{0,n+1} \otimes \omegas\\
		-\sum_{k=0}^{\theN-1} \left( L^{(p)}_{k+1} +  \frac{(-\pGenFunc_1)^{k+1}}{(k+1)!}\right) \ul\omega^{1,k}.
 	\end{multline*}
 	Hence the previous lemma gives the compatibility 
	$$ 
 		\Phi_{\!\sP} \circ \nabla_\sP(\pi\uomegav) = \nabla_\sP \circ \Phi_{\!\sP}(\pi\uomegav)
	$$
	of the Frobenius with the connection on $\pi \uomegas$.  
	Similarly, if we consider the basis $\ol\pi\uomegas$, then we have
	$
		 \Phi_{\!\sP} \circ \nabla_\sP(\ol\pi\uomegasv) = \ol\pi\Phi_{\sL(1)}({\bomega}^{*\vee}), 
	$
	and, if we calculate the composition of the Frobenius with the connection, we have
	\begin{multline*}
	 \nabla_\sP \circ \Phi_{\!\sP}(\ol\pi\uomegasv) = \nabla_\sP
	 \left( \uomegasv + \sum_{m=0}^\theN \sum_{n=0}^{\theN-m} \pPolFunc_{m+1,n} \ul\omega^{m,n} \right)\\
	 = { \bomega}^{*\vee} + \sum_{m=0}^\theN \sum_{n=0}^{\theN-m}
		\ul\omega^{m,n} \otimes d \pPolFunc_{m+1,n} 
	 	+ \sum_{m=1}^\theN \sum_{n=0}^{\theN-m} \pPolFunc_{m,n} \ul\omega^{m,n} \otimes \omega \\
	  	+ \sum_{m=0}^\theN \sum_{n=0}^{\theN-m}\pPolFunc_{m+1,n} \ul\omega^{m,n+1} \otimes \omegas	
	 =   {\bomega}^{*\vee}- \sum_{k=1}^\theN L_k^{(p)} \ul\omega^{0,k} \otimes \omega,
	\end{multline*}
	where the last equality follows from the differential equation satisfied by $\pPolFunc_{m+1,n}$ 
	defined in Theorem \ref{thm: existence one}. Hence the  previous lemma gives the compatibility
	$
		\Phi_{\!\sP} \circ \nabla_\sP(\ol\pi\uomegasv) = \nabla_\sP \circ \Phi_{\!\sP}(\ol\pi\uomegasv)
	$
	of the Frobenius with the connection on $\ol\pi \uomegasv$.
	This finishes the proof of our assertion.
\end{proof}

%
\subsection{Specialization of the elliptic polylogarithm}\label{subsection: second to last}%
%

We next calculate the specialization of the $p$-adic elliptic polylogarithm to non-zero torsion point $z_0$ of 
$E(K)$ of order prime to $\frp$.  Note that by the theory of complex multiplication, we have a commutative
diagram
$$
	\begin{CD}
		\Spec\,\cO_K @>{i_{z_0}}>> E \otimes_{\cO_\bsK} \cO_K   \\
		@A{\sigma}AA  @A{\phi}AA  \\
		\Spec\,\cO_K @>{i_{z_0}}>> E \otimes_{\cO_\bsK} \cO_K.
	\end{CD}
$$
Hence the induced map $i_{z_0} : \sV \rightarrow \sE$ is a morphism of syntomic data.
Denote by $i_{z_0}^* \sLog^\theN$ the filtered Frobenius module defined as 
the pullback of the logarithm sheaf $\sLog^\theN$ to $z_0$.  Then
$$
	i_{z_0}^* \sLog^\theN = \bigoplus_{0 \leq m+n \leq \theN} K \ul\omega^{m,n},
$$
with the Hodge filtration given by the direct sum and the Frobenius 
$
	 \Phi : K  \otimes_{\sigma}i_{z_0}^* \sLog^\theN \xrightarrow\cong i_{z_0}^* \sLog^\theN
$
given by 
$$
	\Phi(\ul\omega^{m,n}) := \frac{1}{\pi^{m} \ol\pi^{n}} \sum_{k=n}^{\theN-m} 
	\frac{(-\pGenFunc_1(z_0))^{k-n}}{(k-n)!} \ul\omega^{m,k}.
$$
In order to calculate the specialization of the polylogarithm sheaf, we will first
describe the splitting of the sheaf $i_{z_0}^* \sLog^\theN$.

\begin{lemma}\label{lem: coleman}
	We have	$F_{z_0,1}(0) \in K$ and
	$$
		\left( 1  - \frac{\sigma}{\ol\pi} \right) F_{z_0,1}(0)= \pGenFunc_1(z_0).
	$$
\end{lemma}

\begin{proof}
	Let $n$ be the smallest integer $>1 $ such that $\pi^n z_0 = z_0$ in $E(K)$.  Then repeatedly 
	using the equality $F_{z_0,1}(z) = F^{(p)}_1(z + z_0) + \ol\pi^{-1} F_{\pi z_0,1}(\pi z)$ of \eqref{eq: funny},
	we have
	$$
		F_{z_0,1}(0) = (1-\ol\pi^{-n})^{-1} \sum_{k=0}^{n-1} \ol\pi^{-k} F^{(p)}_1(\pi^k z_0).
	$$
	Since $F^{(p)}_1(z)$ is a rational function defined over $\bsK$, the above equality
	shows that $F_{z_0,1}(0) \in K$.  Furthermore, we have by the theory of
	complex multiplication $(F_1^{(p)}(z_0))^\sigma = F_1^{(p)}( \pi z_0)$.  Then the above formula
	shows that
	$
		\sigma(F_{z_0,1}(0)) = F_{\pi z_0,1}(0).
	$
	Hence our assertion now follows from the fact that
	$F^{(p)}_1(z_0) = F_{z_0,1}(0) - \ol\pi^{-1}F_{\pi z_0,1}(0)$.
\end{proof}

We first describe the unique splitting of $i_{z_0}^* \sLog^{(1)}$ as a filtered Frobenius module.  We define an element $\ul e'$ 
which corresponds to the basis denoted by the same character in the Hodge case.

\begin{lemma}
	Let $\ul e' = \ul e - F_{z_0,1}(0) \uomegasv  \in F^0(i_{z_0}^* \sLog^{(1)})$.  Then the Frobenius acts
	on this element as
	$
		\Phi(\ul e') = \ul e'.
	$
	In particular, by mapping the basis of $K(0)$ to $\ul e'$, we have a splitting 
	as filtered Frobenius modules of the sequence
	$$
		0 \rightarrow \sH \rightarrow i^*_{z_0}\sLog^{(1)} \rightarrow K(0) \rightarrow 0.
	$$
\end{lemma}

\begin{proof}
	Since we have $\Phi(\uomegasv) = \uomegasv/\ol\pi$, the previous lemma gives the equality
	$\sigma ( F_{z_0,1}(0)) \Phi(\uomegasv) = (F_{z_0,1}(0) - F^{(p)}_{1}(z_0)) \uomegasv$.
	Then the action $\Phi(\ul e) = \ul e - F^{(p)}_1(z_0) \uomegasv$ gives the equality
	$$
		\Phi(\ul e') = \Phi(\ul e) - \sigma ( F_{z_0,1}(0)) \Phi(\uomegasv) = \ul e - F_{z_0,1}(0) \uomegasv = \ul e'
	$$
	as desired.
\end{proof}

We let $\ul e^{m,n} := \ul e'^a \ul\omega^{\vee m} \ul\omega^{*\vee n}/a!$, where $a = \theN - m - n$.   
This basis gives a splitting of $i_{z_0}^* \sLog^\theN$ as filtered Frobenius modules.
See Lemma \ref{lem: split hodge} for the splitting principle in the Hodge case.

\begin{lemma}[Splitting Principle]\label{lem: splitting}
	We have a splitting of filtered Frobenius modules
	$$
	 	i_{z_0}^* \sLog^\theN  \cong \prod_{j = 0}^\theN \Sym^j \sH
	 $$
	 by mapping  $\ul e^{m,n}$ to $\ul\omega^{\vee m} \ul\omega^{* \vee n}$.
\end{lemma}

\begin{proof}
	The splitting is the $N$-th symmetric product of the splitting given in the previous lemma.
\end{proof}

Let $z_0$ be a torsion point in $E(K)$ of order prime to $\frp$ as above.
We denote by $i_{z_0}^* \sP^\theN$ the pullback of the polylogarithm sheaf $\sP^\theN$ to $z_0$.
We now describe $i_{z_0}^* \sP^\theN$ using the basis $\ul e^{m, k}$.   
By construction of the polylogarithm sheaf, $i_{z_0}^* \sP^\theN$ is a $K$-vector space
$$
	i_{z_0}^* \sP^\theN = \sH \bigoplus	i_{z_0}^* \sLog^\theN,
$$
endowed with a Frobenius $\Phi : K \otimes_\sigma i_{z_0}^* \sP^\theN \xrightarrow\cong i_{z_0}^* \sP^\theN$
induced from that of $\sP^\theN$.  Then we have the following.

\begin{proposition}\label{prop: Frob splitting} 
	The Frobenius on $i_{z_0}^* \sP^\theN$ is expressed as 
	\begin{align*}
		\Phi(\pi \uomegav) &=  \uomegav 
		+ \sum_{k=0}^\theN \left( 1-\frac{\sigma}{\ol\pi^{k+1}} \right)
			\frac{F_{z_0,1}(0)^{k+1}}{(k+1)!} \ul e^{0,k}\\
			& \hspace{3cm} 
			+\sum_{m=1}^\theN \sum_{k=0}^{\theN-m}
			  \wh E^{(p)}_{z_0, m, k+1}(0)
			 \ul e^{m,k},\\
		\Phi(\ol\pi \uomegasv) &=  \uomegasv +  \sum_{m=0}^\theN \sum_{k=0}^{\theN-m} 
		\wh E^{(p)}_{z_0, m+1, k}(0) \ul e^{m,k},
	\end{align*}
	using the basis $\ul e^{m, k}$ of $i_{z_0}^* \sLog^\theN$.
\end{proposition}

\begin{proof}
	By definition of $\ul \omega^{m,n}$ and $\ul e'$, we have
	$
		\ul \omega^{m,n} := \ul e^a \ul\omega^{\vee m} \ul\omega^{*\vee n}/a!
			= (\ul e' + F_{z_0,1}(0) \uomegasv)^a \ul\omega^{\vee m} \ul\omega^{*\vee n}/a!
	$
	for $a = N - m - n$.  This implies that
	$$
			\ul \omega^{m,n} = \sum_{k=n}^{\theN - m} \frac{F_{z_0,1}(0)^{k-n}}{(k-n)!}  \ul e^{m,k}.
	$$
	Using this formula, we have
	$$
		\sum_{n=0}^{ \theN} \frac{(-\pGenFunc_1(z_0))^{n+1}}{(n+1)!} \ul \omega^{0,n}
	 = 	\sum_{k=0}^{ \theN}  \sum_{n=0}^{k} 	
	  \frac{(-\pGenFunc_1(z_0))^{n+1} F_{z_0,1}(0)^{k-n}}{(n+1)!(k-n)!}
		 \ul e^{0,k}.
	$$
	By Lemma \ref{lem: coleman}, the sum in the coefficient of $\ul e^{0,k}$ is thus
	$$
		 \frac{F_{z_0,1}(0)^{k+1}}{(k+1)!} -  \frac{\sigma(F_{z_0,1}(0))^{k+1}}{\ol\pi^{k+1}(k+1)!}  
		=\left( 1 - \frac{\sigma}{\ol\pi^{k+1}} \right)
			\frac{F_{z_0,1}(0)^{k+1}}{(k+1)!}.
	$$
	For the coefficients for the other basis, we have
	$$
		  \sum_{n=0}^{\theN-m} \pPolFunc_{m,n+1}(z_0) \ul\omega^{m,n}
		  =  \sum_{k=0}^{\theN - m} \kern-1mm \sum_{n=0}^{k}
		 \pPolFunc_{m,n+1}(z_0)
		  \frac{F_{z_0,1}(0)^{k-n}}{(k-n)!} \ul e^{m,k}.
	$$
	By \eqref{eq: D explicit}, the coefficient of $\ul e^{m,k}$ is given by
	$$
			 \sum_{n=1}^{k+1}  \pPolFunc_{m,n}(z_0)  \frac{F_{z_0,1}(0)^{k+1-n}}{(k+1-n)!} \\
			   = \wh E^{(p)}_{z_0, m, k+1}(0),
	$$
	where the equality holds since $D^{(p)}_{m,0} \equiv 0$.
	This gives the first equality.    
	Next, we consider the second equality.  We have
	$$
		 \sum_{n=0}^{\theN-m} \pPolFunc_{m+1,n}(z_0) \ul\omega^{m,n} 
		  =  \sum_{k=0}^{\theN-m} \sum_{n=0}^{k}  \pPolFunc_{m+1,n}(z_0) 
		   \frac{F_{z_0,1}(0)^{k-n}}{(k-n)!} \ul e^{m,k}.
	$$
	Again by \eqref{eq: D explicit}, the coefficient of $\ul e^{m,k}$ above becomes
	$$
		\sum_{n=0}^{k}  \frac{F_{z_0,1}(0)^{k-n}}{(k-n)!}   \pPolFunc_{m+1,n}(z_0) 
		=\wh E^{(p)}_{z_0, m+1,k}(0),
	$$
	which proves the second assertion.
\end{proof}

%
\subsection{Specialization in syntomic cohomology}
%

We now calculate the extension class of the specialization of the elliptic polylogarithm in syntomic cohomology.
Let $\sV = \Spec\, \cO_K$, and let $z_0 \in E(K)$ be as in \S \ref{subsection: second to last}.
The specialization of the elliptic polylogarithm $i_{z_0}^* \sP^\theN$ gives a cohomology class 
$$
	 \pol^\theN_{z_0} \in \Ext^1_{S(\sV)}( \sH, \sLog(1)) = H^1_{\syn}(\sV, \sH^\vee \otimes i_{z_0}^* \sLog^\theN(1)).
$$
We use the results of the previous section to calculate this element explicitly.

\begin{lemma}
	We have an isomorphism of $K$-vector spaces
	\begin{equation}\label{eq: one isom}
		H^1_{\syn}(\sV, \sH^\vee \otimes i_{z_0}^* \sLog^\theN(1))  \xrightarrow\cong
		 \sH^\vee \otimes i_{z_0}^* \sLog^\theN /
		 \bigoplus_{n=0}^\theN K \ul \omega \otimes \ul e^{0,n}
	\end{equation}
\end{lemma}

\begin{proof}
	By \eqref{eq: syn base}, we have $H^1_\syn(\sV, M) \xrightarrow\cong M/ (1 - \Phi) F^0 M$ for any filtered 
	Frobenius module $M$ in $S(\sV)$, where the map is given by mapping the extension class $[M']
	\in \Ext^1_{S(\sV)}(K(0), M) = H^1_{\syn}(\sV, M)$ of an extension
	$$
		0 \rightarrow M \rightarrow M' \rightarrow K(0) \rightarrow 0
	$$
	to
	$
		(1 - \Phi) e \in M,
	$
	where $e$ is any lifting of the fixed basis of $K(0)$ to $F^0 M'$.
	The lemma follows from the fact that $F^0( \sH^\vee \otimes i_{z_0}^* \sLog^\theN(1)) 
	=  \bigoplus_{n=0}^\theN K \ul \omega \otimes \ul e^{0,n}$, and that  $(1 - \Phi)$ 
	acts as an isomorphism on this subspace.
\end{proof}

\begin{theorem}\label{thm: main}
	Let $z_0$ be a non-zero torsion point in $E(K)$ of order prime to $\frp$.  In other words, $z_0$ is a non-zero point
	whose annihilator as a $\cO_\bsK$-module is prime to $\frp$.  Then the image of $ [\pol^\theN_{z_0}]$ with respect 
	to the isomorphism
	$$
		 \Ext^1_{S(\sV)}( \sH, i_{z_0}^* \sLog^\theN(1))  \xrightarrow\cong  \sH^\vee \otimes i_{z_0}^* \sLog^\theN
		 / \bigoplus_{n=0}^\theN K \ul \omega \otimes \ul e^{0,n}
	$$
	of \eqref{eq: one isom} is 
	$$
		- \sum_{m=1}^\theN \sum_{k=0}^{\theN-m} \frac{e^{(p)}_{-m,k+1}(z_0)}{ \Omega_\frp^{k-m} }
		 \ul\omega \otimes \ul e^{m,k}
		- \sum_{m=0}^\theN \sum_{k=1}^{\theN-m}
		\frac{ e^{(p)}_{-m-1,k}(z_0)}{\Omega_\frp^{k-m-2}} \uomegas \otimes \ul e^{m,k},	
	$$
	where the $e^{(p)}_{a,b}$ are the $p$-adic Eisenstein-Kronecker numbers defined in Definition \ref{def: p EK}.
\end{theorem}

\begin{proof}
	The theorem follows from Proposition \ref{prop: Frob splitting}, \eqref{eq: one isom} and the definition of the 
	$p$-adic Eisenstein-Kronecker numbers.  The terms containing $\ul\omega \otimes \ul e^{0,n}$
	maps to zero, since $\ul\omega \otimes \ul e^{0,n} \in F^0(\sH^\vee \otimes i_{z_0}^* \sLog^\theN(1))$.
\end{proof}

The above result is analogous to the Hodge theoretic calculation given in Theorem \ref{thm: Beilinson and Levin}.

\appendix

%
%
%
\section{the real Hodge realization}
%
%
%

The Hodge realization of the elliptic polylogarithm was calculated in the original paper by Beilinson and Levin \cite{BL},
as well as \cite{W}.  In the above papers, the Hodge realization was expressed in terms of the $q$-averaged polylogarithm
function, obtained from the classical polylogarithm function on $\bbP^1 \setminus \{0,1,\infty\}$.  In the Appendix, closely
following our method of the $p$-adic case, we will describe how to explicitly describe the $\bbR$-Hodge realization of the 
elliptic polylogarithm using functions given by certain iterated integrals starting from the connection functions $L_n(z)$.

%
\subsection{Real analytic elliptic polylogarithm functions}
%

In this subsection, we will define the Eisenstein-Kronecker functions $E_{m,b}(z)$ and the real analytic elliptic polylogarithm
function $G_{m,b}(z)$.
We first investigate the properties of Eisenstein-Kronecker-Lerch series viewed as a function in $z$ and $w$.
For $z, w \in (\bbC \setminus \Gamma)$, we let
$$
	K_{a}(z,w,s) := K^*_a(z,w,s).
$$
Then the integral expression \eqref{eq: integral expression} gives
\begin{equation}\label{eq: int express two}
	A^s \Gamma(s) K_a(z,w,s) = I_a(z,w,s) + I_{a}(w,z,a+1-s) \pair{w,z}.
\end{equation}
if $a \geq 0$.
Using this integral expression, we may prove that  $K_a(z,w,s)$ is a $\sC^\infty$-function for $(z,w)$, at first for
$a \geq 0$, then for any integer $a$ by the definition of $K_a(z,w,s)$ and reduction to the case for $a>0$.
Moreover, the derivatives of $K_a(z,w,s)$ with respect to $\partial_z$, $\partial_{\ol z}$, $\partial_{w}$, $\partial_{\ol w}$
are all analytic in $s$.

The Eisenstein-Kronecker-Lerch series for various integers $a$ and $s$  are related 
by the following differential equations. 
\begin{lemma} \label{lemma: calculation of differentials}
   Let $a$ be an integer, and consider $K_a(z,w,s)$ to be a function for $z$, $w \in (\bbC \setminus \Gamma)$.
   Then $K_a(z,w,s)$ satisfies the differential equations
   \begin{equation*}
       \begin{aligned}  
           \partial_{z} K_a(z,w,s) &= -s K_{a+1}(z,w,s+1) \\
           \partial_{\ol{z}}K_{a}(z,w,s) &= (a-s)K_{a-1}(z,w,s) \\
           \partial_{w} K_a(z,w,s) &= - ( K_{a+1}(z,w,s) - \ol{z} K_{a}(z,w,s) ) /A \\	
           \partial_{\ol{w}}K_{a}(z,w,s) &=  (K_{a-1}(z,w,s-1) - z K_{a}(z,w,s) )/A.
        \end{aligned}
   \end{equation*}
\end{lemma}


We next define the Eisenstein functions, which will be used to express the complex periods of the 
elliptic polylogarithm sheaf.

\begin{definition}
	For any integer $m$ and $b$, we define the Eisenstein function $E_{m, b}(z)$
	to be the $\sC^\infty$-function on $\bbC \setminus \Gamma$ given by
	$$
		E_{m, b}(z) := K^*_{b-m}(0,z,b).
	$$
\end{definition}
Compare \cite{We}VIII \S 14, where a slightly different definition is used.  By definition, we have
\begin{equation}\label{eq: Eis and EK}
	E_{m,b}(z_0) = e^*_{-m,b}(z_0)
\end{equation}
for $z_0 \in \bbC \setminus \Gamma$.  
Note that $E_{m,b}(z)$ satisfies
$
	\overline{E_{m,b}(z)} = (-1)^{b-m} E_{b,m}(z)
$
under complex conjugation. By definition, we have
$
	E_{m,b}(z + \gamma) = E_{m,b}(z)
$
for any $\gamma \in \Gamma$.  Calculations similar to Lemma \ref{lemma: calculation of differentials} give the equality
\begin{align}\label{eq: diff E}
	\partial_z E_{m+1, b}(z)  &=  - E_{m,b}(z)/A,  &
	\partial_{\ol z} E_{m, b+1}(z) &=  E_{m, b}(z)/A.
\end{align}

\begin{remark}\label{rem: zero}
	Note that by \eqref{eq: value 00}, we have $E_{0,0}(z) \equiv -1$, and \eqref{eq: diff E} shows that we have $
	E_{-a,0}(z) = 0$ for any $a > 0$.   This implies that for any $z_0 \in \bbC$, we have
	$e_{0,0}^*(z_0)=-1$ and $e^*_{a,0}(z_0) = 0$ if $a>0$.
\end{remark}

By Lemma \ref{lemma: calculation of differentials}, we may prove that 
$$
	 K_1(z,w,1) = \sum_{b \geq 0} (-1)^{b-1} E_{0,b}(w) z^{b-1}.
$$
Using the Kronecker theorem
$$
	\Theta(z,w) = \Exp{\frac{z \ol w}{A}} K_1(z,w,1)
$$
 (see \cite{We} VIII \S4 p. 71 (7) or \cite{BK1} Theorem 1.13 for a proof),
 Definition 1.3, Definition 1.4 and the fact that $E_{0,1}(z) = F_1(z) - \ol z/A$, 
we may prove that the relation between the connection function $L_n(z)$ and $E_{0,b}(z)$ is given by
\begin{equation}\label{eq: Xi and E}
	L_n(z) = (-1)^{n-1} \sum_{b=0}^n \frac{E_{0,1}^{n-b}(z)}{(n-b)!} E_{0,b}(z).
\end{equation}

%
%

We next define real analytic functions which corresponds to the periods of the $\sC^\infty$-sheaf associated
to the elliptic polylogarithm sheaf.  In the next section, we will use these functions to construct
holomorphic functions which are periods of the elliptic polylogarithm sheaf.

\begin{definition}[Real analytic elliptic polylogarithm function]
	We define the multivalued functions $G_{m,b}(z)$ on $\bbC \setminus \Gamma$ by first letting
	\begin{align*}
		G_{m,-1}(z) &=0, &
		G_{0,b}(z) &=  E_{0,b}(z),
	\end{align*}
	for integers $m$, $b \geq 0$, and we then iteratedly define the functions $G_{m,b}(z)$ for $m$, $b \geq 0$ to be any function satisfying 
	\begin{align*}
		\partial_z G_{m+1, b}(z) &= - G_{m, b}(z), &
		\partial_{\ol z} G_{m, b+1}(z) &= G_{m, b}(z)/A.  
	\end{align*}
\end{definition}

If we assume the existence of $G_{m+1,b}(z)$ and $G_{m,b+1}(z)$, 
then  we have
$$\partial_z (G_{m+1, b}(z)/A) = -G_{m,b}(z)/A = - \partial_{\ol z}G_{m, b+1}(z),
$$
which implies that $G_{m, b+1}(z) dz  + G_{m+1, b}(z) d \ol z/A $ is a closed form.
Hence the function $G_{m+1, b+1}(z)$ exists in this case.
We fix a choice of $G_{m,b}(z)$ satisfying the condition in 
Proposition \ref{prop: relation to Eisenstein} below.

\begin{proposition}\label{prop: relation to Eisenstein}
	We may iteratedly choose $G_{m,b}(z)$ so that	
	$$
		A^{b} G_{m,b}(z) + (-1)^{m+b} A^{m}\overline{G_{b,m}(z)} =  
		A^{m+b} E_{m,b}(z) - \frac{(-z)^m \ol z^b}{m! b!}
	$$
	for any integers $m$, $b \geq 0$.
\end{proposition}

%
\subsection{Elliptic polylogarithm functions}
%

We next construct the elliptic polylogarithm function $D_{m,n}(z)$ and $D_{m,n}^*(z)$, 
which are holomorphic multi-valued functions on $\bbC \setminus \Gamma$.

\begin{definition}
	For integers $m$, $n \geq 0$, we define the elliptic polylogarithm functions $D_{m,n}(z)$ and
	$D_{m,n}^*(z)$ by
	\begin{align*}
		D_{m,n}(z) &=(-1)^{n-1} \sum_{k=0}^n \frac{E_{0,1}(z)^{n-k}}{(n-k)!}  G_{m,k}(z) \\
		D_{m,n}^*(z) &= D_{m,n}(z) - \frac{(-1)^{m+n}}{m!n!} z^m F_1(z)^n.
	\end{align*}
\end{definition}

\begin{lemma}
	Both $D_{m,n}(z)$ and $D_{m,n}^*(z)$ are holomorphic functions on the universal covering of $\bbC \setminus \Gamma$.
\end{lemma}

\begin{proof}
	The statement for   $D_{m,n}(z)$ follows from the fact that 
	the functions in the sum are defined on $\bbC \setminus \Gamma$ and that
	$\partial_{\ol z} D_{m,n}(z) = 0$, which follows from the fact that
	$\partial_{\ol z} E_{0,1}(z) = -1$, $\partial_{\ol z} G_{m,0}(z) = 0$, and
	$\partial_{\ol z} G_{m,k}(z) =  G_{m,k-1}(z)/A$ for $k \geq 1$.  
	The statement for $D_{m,n}^*(z)$
	follows from the fact that $F_1(z)$ is holomorphic on $\bbC \setminus \Gamma$.
\end{proof}

\begin{lemma}\label{lemma: first}
	The functions $D_{m,n}(z)$ and $D_{m,n}^*(z)$ for $m, n \geq 0$ satisfy
	\begin{align*}
		d D_{m+1,n}(z) &=  - D_{m,n}(z) dz   -  D_{m+1,n-1}(z) d F_1  \\
		d D_{m+1,n}^*(z) &= - D_{m,n}^*(z) dz - D_{m+1,n-1}^*(z) d F_1,
	\end{align*}
	where we let $D_{m,-1}(z) = D^*_{m,-1}(z)  =0$.
\end{lemma}

The $p$-adic analogues of $D_{m,n}(z)$ and $D^*_{m,n}(z)$ are the overconvergent functions
$D^{(p)}_{m,n}$ given in Theorem \ref{thm: existence two}.  By definition, we have
\begin{align*}
		D_{0,n}(z) &= (-1)^{n-1} \sum_{k=0}^n \frac{E_{0,1}(z)^{n-k}}{(n-k)!}  E_{0,k}(z)  = L_n(z),\\
		D_{0,n}^*(z) &= D_{0,n}(z) - \frac{(-1)^n}{n!} F_1(z)^n =  L_n(z) - \frac{(-1)^n}{n!} F_1(z)^n.
\end{align*}
Hence the elliptic polylogarithm function $D_{m,n}(z)$ is obtained from $L_n(z)$ by iterated integration.

%

%
\subsection{Review of absolute Hodge cohomology}
%

In this section, we will freely use terminology concerning variation of mixed Hodge structures and absolute Hodge cohomology.
See for example \cite{HW}  Appendix A for details.   For any variety $X$ smooth and separated of finite type over $\bbC$,
we denote by $\VMHS_\bbR(X)$ the category of polarizable admissible variation of mixed $\bbR$-Hodge 
structures  on $X$.   For any $\sF$ in $\VMHS_\bbR(X)$,  we denote by $\cF$ the underlying locally free 
$\cO_X$-module with integrable connection,  and by $\sF_\bbR$ and $\sF_\bbC$ the underlying $\bbR$- and $\bbC$-local systems.   
If we let  $X^\an := X(\bbC)$, then the theorem of de Rham gives a canonical isomorphism
$$
	H^i_\dR(X, \cF) \cong H^i_B(X^\an, \sF_\bbR ) \otimes_\bbR \bbC
$$
between de Rham and Betti cohomologies.  We denote by $H^i(X^\an, \sF)$ the mixed $\bbR$-Hodge structure
induced from the above isomorphism,
and by $H^i_\AH(X, \sF)$ the  $m$-th absolute Hodge cohomology of $X$ with coefficients in
$\sF$.  Note that for $i=0,1$, there exists a canonical isomorphism
\begin{equation}
		\Ext^i_{\VMHS_\bbR(X)}(\bbR(0), \sF) \xrightarrow\cong H^i_\sA(X, \sF).
\end{equation}

Let $S = \Spec\, \bbC$.  Then $\VMHS_\bbR(S)$ is the category of polarizable mixed 
$\bbR$-Hodge structures $\MHS_\bbR$.  
The Leray spectral sequence
$$
	E_{2}^{p,q} = H^p_\sA(S, H^q(X^\an, \sF)) \Rightarrow H^{p+q}_\sA(X, \sF)
$$
degenerates to give the short exact sequence
\begin{equation}\label{eq: C SES}
	0 \rightarrow H^1_\sA( S, H^0(X^\an, \sF)) \rightarrow 
	H^1_\sA(X, \sF) \\
	\rightarrow H^0_\sA(S, H^1(X^\an, \sF)) \rightarrow 0.
\end{equation}
Consider the map
$$
	H^0_\sA(S, H^1(X^\an, \sF)) \hookrightarrow 
	H^1_B(X^\an, \sF_\bbC) \cong H^1_\dR(X, \cF),
$$
where the last isomorphism is de Rham's theorem.  Then we have a commutative diagram
\begin{equation}\label{eq: C compatible}
		\begin{CD}
				\Ext^1_{\VMHS_\bbR(X)}(\bbR(0), \sF)  @>\For>> 	\Ext^1_{M(X)}(\cO_{X}, \cF)  \\
				@V{\cong}VV   @V{\cong}VV\\
				 H^{1}_\sA(X, \sF)  @>>> 	 H^{1}_\dR(X, \cF),
		\end{CD}
\end{equation}
where $M(X)$ is the category of locally free $\cO_X$-modules with integrable connection
defined in Definition \ref{def: M(X)}, and $\For: \VMHS_\bbR(X) \rightarrow M(X)$ is the functor 
associating to any $\sF$ the underlying coherent module with connection $\cF$ on $X$,
forgetting the Hodge filtration and the $\bbR$-structure.

%
\subsection{The Logarithm sheaf}
%

We now define the logarithm sheaf.
Let $E$ be an elliptic curve defined over $S = \Spec\, \bbC$.  We let $H^1(E, \bbR)$ be the pure $\bbR$-Hodge 
structure of weight $1$ given by $H^1_\dR(E) \cong H^1_B(E^\an, \bbR) \otimes \bbC$, and we let
$\sH := H^1(E, \bbR)^\vee$ be the dual Hodge structure.
We let the notations be as in Remark \ref{rem: choice of basis}. 
Then the $\bbR$-Hodge structure $\sH$ is of pure weight $-1$, with $\bbR$-structure given by 
$\sH_\bbR = H_1(E^\an, \bbR)$ and Hodge filtration defined by $F^{-1}(\sH) = \sH$, $F^0(\sH) = \bbC \uomegasv$ and $F^1(\sH ) = 0$.

We denote by $\phi_\infty$ the complex conjugation on $\sH_\bbC = \sH_\bbR \otimes_\bbR \bbC$ defined by
the action of complex conjugation on $\bbC$.  The class  $\ul\omega$ is represented by $dz$, and since
$d E_{0,1}(z) = d F_1 -  d \ol z/A = \omegas - d \ol z/A$ where $E_{0,1}(z)$ is a single valued real analytic function on $U$
(see Remark \ref{rem: zero}), the class of $\uomegas$ is represented by $d\ol z/A$.
Hence the complex conjugation $\phi_\infty$ acts on these classes as 
\begin{align}\label{eq: c c on H}
	\phi_\infty(\uomegav) &= \uomegasv/A,& \phi_\infty(\uomegasv) &= A \uomegav.
\end{align}
If we let $\gamma_1 := (\ul\omega^\vee + (\uomegasv/A))$ and $\gamma_2 := i (\ul\omega^\vee - (\uomegasv/A))$,
then $\gamma_1$ and $\gamma_2$ form a basis of $\sH_\bbR$.

For any smooth scheme $X$ over $S$, we denote by $\sH_X$ the constant variation of $\bbR$-Hodge structures on $X$.
The underlying coherent module with connection of $\sH$ is $\cH_X$.   We have a natural isomorphism
$$
	H^0_\sA(S, H^1(E^\an, \sH_E)) = H^0_\sA( S, \sH^\vee \otimes \sH) = \Hom_{\MHS_\bbR}(\sH, \sH).
$$
Hence the short exact sequence \eqref{eq: C SES} defined from the Leray spectral sequence 
gives the short exact sequence
\begin{equation*}
	0 \rightarrow H^1_\sA(S, \sH) \rightarrow H^1_\sA(E, \sH_E) \rightarrow
	\Hom_{\MHS_\bbR}( \sH, \sH) \rightarrow 0.
\end{equation*}
In addition, we have in this case
$$
	H^1_\sA(S, \sH) = \Ext^1_{\MHS_\bbR}(\bbR(0), \sH) = \sH_\bbC/( \sH_\bbR + F^0 \sH_\bbC) = 0.
$$
Hence the above exact sequence gives an isomorphism
\begin{equation}\label{eq: C deg}
	 H^1_\sA(E, \sH_E) \xrightarrow\cong \Hom_{\MHS_\bbR}( \sH, \sH).
\end{equation}

\begin{definition}\label{def: C log}
	We define the sheaf  $\sLog^{(1)}$ to be any extension of $\sH_E$ by $\bbR(0)$ in $\VMHS_\bbR(E)$,
	whose extension class in
	$$
		\Ext^1_{\VMHS_\bbR}(\bbR(0), \sH_E) \cong H^1_\sA(E, \sH_E)
	$$
	is mapped to the identity through the isomorphism of \eqref{eq: C deg}.
	We define the $N$-th logarithm sheaf $\sLog^\theN$ to be the $N$-th symmetric tensor product of $\sLog^{(1)}$.
\end{definition}

Unlike the case of the de Rham realization (see Remark \ref{rem: not canonical}), the logarithm sheaf 
for the Hodge realization is determined uniquely up to unique isomorphism.  
Denote by $i_{[0]}^*$  the pull-back by the identity $i_{[0]}: S \rightarrow E$ of the elliptic curve.
Since the extension of $\sH$ by $\bbR(0)$ is split on $S$, there exists a splitting 
$$
	\varphi: i_{[0]}^* \sLog^{(1)} \cong \bbR(0) \bigoplus \sH
$$ as mixed $\bbR$-Hodge structures.  
This splitting is unique due to weight reasons.

We will next describe  $\sLog^{(1)}$ as a variation of $\bbR$-mixed Hodge structures on $E$.  By
\eqref{eq: C compatible}, the natural map
$
	H^1_\sA( E, \sH_E) \rightarrow H^1_\dR(E, \cH_E)
$
is given in terms of extensions by the forgetful functor.  Hence by the de Rham characterization of $\sLog^{(1)}$
in Definition \ref{def: log dR}, the underlying coherent module with connection on $E$ of $\sLog^{(1)}$ is 
the locally free $\cO_E$-module with connection $\cLog^{(1)}$ of Proposition \ref{pro: log dR}.  
We will now explicitly describe a Hodge filtration, an $\bbR$-structure, and a weight filtration on this module.

We take an affine open covering $\frU = \{ U_i \}_{i \in I}$ of $E$, and cohomology classes
$(\eta_i, u_{ij})$ and $(u_i)$ as in Definition \ref{def: basis}.
We denote by $\cLog_i^{(1)}$ the module
$$
		\cLog^{(1)}_i := \cO_{U_i} \ul e_i \bigoplus \cH_{U_i}
$$
defined in Proposition \ref{pro: log dR},  with connection 
$\nabla(\ul e_i) =\uomegav \otimes \omega + \uomegasv \otimes \uomegas_i$.
We define the Hodge filtration to be the direct sum.  The Hodge filtration extends to a filtration of $\cLog^{(1)}$ on $E$ since
the difference $\ul e_j - \ul e_i = u_{ij} \uomegasv  \in F^0 \cH_{U_i \cap U_j}$ on $U_i \cap U_j$.
The $\bbR$-structure is defined as follows.  Let $\sH_\bbR$ be the $\bbR$-structure of $\sH_X$ for any $X$.  
Let $\hbs_i := \ul e_i - z \uomegav - \xi_i(z) \uomegasv$, where $\xi_i(z) = F_1(z) - u_i(z)$. 
Then, since $d \xi_i =  \omegas_i$, by the definition of the connection,  we have $\nabla(\hbs_i) = 0$.  
The horizontal sections $\hbs_i$ are compatible and define a multi-valued section $\hbs_0$
of $\cLog^{(1)}$ on $E^\an$.  
We define the $\bbR$-structure of $\sLog^{(1)}_{\bbR}$ on the universal covering space of $E^\an$ to be the structure 
given by $\hbs_0$ and $\sH_{\bbR}$.  
One may prove that the $\bbR$-structure defined by $\hbs_0$ and $\sH_\bbR$ is invariant under the action of monodromy, 
hence this structure descends to give an $\bbR$-structure $\sLog^{(1)}_\bbR$ on $E^\an$.
We define the weight filtration $W_\bullet$ on $\sLog^{(1)}_\bbR$ by $W_{-2} \sLog^{(1)}_\bbR =0$, 
$W_{-1} \sLog^{(1)}_\bbR =\sH_\bbR$ and
$W_{0} \sLog^{(1)}_\bbR = \sLog^{(1)}_\bbR$.  The Hodge filtration, $\bbR$-structure and the weight filtration  above give
$\cLog^{(1)}$ a structure of a variation of $\bbR$-mixed Hodge structures on $E$, which we
denote by $\sLog^{(1)}$.

\begin{proposition}\label{pro: C split}
	We define $\sLog^{(1)}$ to be the variation of mixed $\bbR$-Hodge structures on $E$
	given above.  Then this satisfies the property of the first logarithm sheaf in Definition \ref{def: C log}.
\end{proposition}

\begin{proof}
	By construction, the class of the underlying coherent module with connection $\cLog^{(1)}$ of $\sLog^{(1)}$ maps
	to the element in the first de Rham cohomology which corresponds to the identity.  Hence by \eqref{eq: C compatible}, 
	the isomorphism of \eqref{eq: C deg} maps the class of $\sLog^{(1)}$ to the identity.  This gives our assertion.
\end{proof}

\begin{remark}\label{rem: C split}
	By construction, the splitting $\varphi: i_{[0]}^* \cLog^{(1)} \cong \bbC \bigoplus \cH$
	 given in \eqref{eq: split log one} of the underlying $\bbC$-vector spaces gives the unique splitting 
	 of $i^*_{[0]}\sLog^{(1)}$ compatible with
	 the Hodge filtration and the $\bbR$-structure.
\end{remark}

The $N$-th logarithm sheaf is the $N$-th symmetric tensor product of $\sLog^{(1)}$. The Hodge filtration, 
the $\bbR$-structure and the weight filtration is defined naturally by taking the symmetric tensor product of each structure.  

%
\subsection{The Polylogarithm sheaf}
%

Next, we use the logarithm sheaf of the previous subsection to define the polylogarithm class 
in absolute Hodge cohomology.  Then we will explicitly describe the polylogarithm sheaf,
which is defined to be the pro-variation of mixed Hodge structures corresponding to the
polylogarithm class.

Let $D = [0]$ and $U = E \setminus [0]$.
The calculation of the cohomology of $\cLog^\theN$ given in Lemma \ref{lem: log coh} 
is compatible with the Hodge structures as follows.

\begin{lemma}\label{lem: betti log coh} 
	The localization sequence gives isomorphisms of pro $\bbR$-mixed Hodge structures
	\begin{align*}
		\varprojlim_\theN H^0(U^\an, \sLog^\theN(1))
			&		\xleftarrow\cong
		 \varprojlim_\theN H^0(E^\an, \sLog^\theN(1))  = 0 \\
		\res: \varprojlim_\theN H^1(U^\an, \sLog^\theN(1)) & \xrightarrow\cong
		 \varprojlim_\theN H^0(D^\an,  i_{[0]}^* \sLog^\theN)=  \varprojlim_\theN i_{[0]}^* \sLog^\theN.
	\end{align*}
\end{lemma}
From \eqref{eq: C SES}, we have the short exact sequence
\begin{multline*}
	0 \rightarrow H^1_\AH(S, H^0(U^\an, \sH_U^\vee \otimes \sLog^\theN)(1))
	\rightarrow H^1_\AH(U, \sH_U^\vee \otimes \sLog^\theN(1)) \\\rightarrow
	H^0_\AH(S, H^1(U^\an, \sH_U^\vee \otimes \sLog^\theN)(1)) \rightarrow 0.
\end{multline*}
Since the projective limit of $H^0(U^\an, \sLog^\theN)$ is zero, this gives a natural isomorphism
\begin{equation}\label{eq: pol isom one}
	\varprojlim_\theN
	H^1_\AH(U, \sH_U^\vee \otimes \sLog^\theN(1)) 
	  \xrightarrow\cong\varprojlim_\theN
	H^0_\AH(S, H^1(U^\an, \sH_U^\vee \otimes \sLog^\theN(1))).
\end{equation}
We have a natural isomorphism
\begin{equation}\label{eq: pol isom two}
	H^0_\AH(S, H^0(D^\an, \sH^\vee_D \otimes i_{[0]}^* \sLog^\theN))
	\xrightarrow\cong H^0_\AH(S, \sH^\vee \otimes i_{[0]}^* \sLog^\theN)
	\cong \bbR,
\end{equation}
where the last isomorphism is obtained by mapping the identity element in 
$\sH^\vee \otimes \sH = \Hom(\sH, \sH)$  to $1 \in \bbR$.  Combining \eqref{eq: pol isom one} 
and \eqref{eq: pol isom two} with the residue isomorphism, we obtain an isomorphism
\begin{equation}\label{eq: res isom two}
	\varprojlim_\theN H^1_\AH(U, \sH_U^\vee \otimes \sLog^\theN(1)) \xrightarrow\cong \bbR.
\end{equation}

\begin{definition}
	As in \cite{BL} (see also \cite{HK} Appendix A),
	we define the polylogarithm class to be a system of classes $\pol^\theN \in  
	H^1_\AH(U, \sH_U^\vee \otimes \sLog^\theN(1))$
	such that
	$$
		\pol := \varprojlim_\theN \pol^\theN \in 	\varprojlim_\theN  H^1_\AH(U, \sH_U^\vee \otimes \sLog^\theN(1))
	$$
	maps through \eqref{eq: res isom two} to $1$.
\end{definition}

\begin{definition}
	We define the elliptic polylogarithm sheaf on $U$ to be a system of variation of mixed
	$\bbR$-Hodge structures $\sP^\theN$ on $U$ given 
	as an extension
	$$
		0 \rightarrow \sLog^\theN(1) \rightarrow \sP^\theN \rightarrow \sH_U \rightarrow 0
	$$
	whose extension class corresponds to $\pol^\theN$ in 
	$$
		\Ext^1_{\VMHS_\bbR(U)}(\sH_U, \sLog^\theN(1)) \cong H^1_\AH(U, \sH_U^\vee \otimes \sLog^\theN(1)).
	$$	
\end{definition}

Since the classes $\pol^\theN$ form an inverse system, there exists a unique surjection $\sP^{\theN+1} \rightarrow \sP^\theN$
compatible with the identity on $\sH_U$ and the projection $\sLog^{\theN+1} (1) \rightarrow \sLog^\theN(1)$.

We next explicitly describe the variations of mixed $\bbR$-Hodge structures $\sP^\theN$.
By definition,  $\sP^\theN$ is given as an extension
$$
		0 \rightarrow \sLog^\theN(1) \rightarrow \sP^\theN \rightarrow \sH_U \rightarrow 0.
$$
The polylogarithm sheaf is characterized by the image of its cohomology class mapped to de Rham cohomology.
By \eqref{eq: C compatible}, this implies that the underlying coherent module with connection of $\sP^\theN$ is
the de Rham realization $\cP^\theN$ given in Corollary \ref{cor: dR P}, and the injectivity of \eqref{eq: C compatible}
implies that $\sP^\theN$ is the variation of mixed $\bbR$-structure determined uniquely up to canonical isomorphism 
which one may equip on $\cP^\theN$.  
We define the Hodge and weight filtrations of $\cP^{\theN}$ as the direct 
sum of the Hodge and weight filtrations on $\sLog^\theN(1)$ and $\sH$.  
In order to define the $\bbR$-structure $\sP^\theN_\bbR$ on $\cP^\theN$,  we first introduce
certain horizontal sections $\hbs$ and $\hbs^*$ of $\cP^\theN$ on the universal covering space of $U^\an$ 
given as follows.

\begin{lemma}
	Let
	\begin{equation}\label{eq: horizontal basis}\begin{split}
		\hbs &=\uomegav  
		- \sum_{n=0}^\theN \frac{(-F_1(z))^{n+1}}{(n+1)!} \ul\omega^{0,n}
		+ \sum_{m=1}^\theN \sum_{n=0}^{\theN-m} \PolFun^*_{m, n+1}(z) \ul\omega^{m,n} \\
	\hbs^* &= \uomegasv 
		+ \sum_{m=0}^\theN \sum_{n=0}^{\theN-m} \PolFun_{m+1, n}(z) \ul\omega^{m,n}.
	\end{split}\end{equation}
	Then $\hbs$ and $\hbs^*$ are horizontal sections of $\cP^\theN$. 
\end{lemma}

\begin{proof}
	The statement follows from the fact that $d F_1 = \omegas$ and the differential equations   
	Lemma \ref{lemma: first} satisfied by $D_{m,n}(z)$ and $D^*_{m,n}(z)$.
\end{proof}

Our choice of the real analytic elliptic polylogarithm function given in 
Proposition \ref{prop: relation to Eisenstein} gives the following.

\begin{proposition}\label{pro: real structure}
	Let $\hbs_1 := \hbs + (\hbs^*/A)$ and $\hbs_2 =  i ( \hbs -( \hbs^*/A))$.  Then
  	the $\bbR$-structure defined by $\hbs_1$, $\hbs_2$ and $\sLog^\theN_\bbR(1)$ on 
	the universal covering space of $U^\an$ descends 
	to give an $\bbR$-structure $\sP^\theN_\bbR$ of $\cP^\theN$ on $U^\an$, which fits into the exact sequence
	$$
		0 \rightarrow \sLog^\theN_\bbR(1) \rightarrow \sP^\theN_\bbR \rightarrow \sH_{\bbR} \rightarrow 0.
	$$
\end{proposition}

The above proposition may be proved by explicitly calculating the action of the complex conjugation $\phi_\infty$
on $\hbs$ and $\hbs^*$.  By the construction of the elliptic polylogarithm, 
the choice of an $\bbR$-structure on $\cP^\theN$ is unique up to isomorphism.  Hence the $\bbR$-structure 
defined above is the real structure of the elliptic polylogarithm.  This gives the main result of our Appendix.

\begin{theorem}
	The elliptic polylogarithm sheaf $\sP^\theN$ is the variation of mixed $\bbR$-Hodge structures on
	$U$ given as an extension
	$$
		0 \rightarrow \sLog^\theN(1) \rightarrow \sP^\theN \rightarrow \sH_U \rightarrow 0, 
	$$whose underlying coherent module with connection is the object $\cP^\theN$ given in Definition 
	\ref{def: de Rham polylogarithm sheaf}, whose Hodge filtration is given as the direct sum of the
	Hodge filtrations on $\sH_U$ and $\sLog^\theN(1)$, whose real structure is the structure 
	$\sP^\theN_\bbR$ given in Proposition \ref{pro: real structure}, and whose weight filtration
	is the direct sum of the weight filtrations on $\sH$ and $\sLog^\theN(1)$.
\end{theorem}

This shows that the holomorphic functions
$\PolFun_{m,n}(z)$ and $\PolFun_{m,n}^*(z)$ which we defined in the previous section are in fact periods of the 
elliptic polylogarithm sheaf. 

%
\subsection{Specialization to points}
%

We now calculate the specialization of the elliptic polylogarithm sheaf to points of the elliptic curve.
Suppose $z_0 \in\bbC$, and we denote by $i_{z_0}^* \sLog^\theN$ the restriction of the
variation of mixed Hodge structures $\sLog^\theN$ to the point $z_0$.  We let $S = \Spec\, \bbC$.

Let $\sM:= \sH^\vee \otimes i_{z_0}^* \sLog^\theN$.
We have an isomorphism
\begin{equation}\label{equation: isom five}
	H^1_\AH(S, \sH^\vee \otimes i_{z_0}^* \sLog^\theN(1)) = M_{\bbC} /(M_{\bbR}(1) + F^{1} M_{\bbC}).
\end{equation}
Denote by $i_{z_0}^* \sP^\theN$ the pull-back of the polylogarithm sheaf to the point $z_0$. 
Note that $\ul\omega \otimes \ul\omega^{0,n} \in F^{1} M_{\bbC}$.  Hence 
the explicit calculation of the period of $i_{z_0}^* \sP^\theN$ shows that the extension class
$$
	\left[i_{z_0}^* \sP^\theN\right] \in \Ext^1_{\MHS_\bbR}(\sH, i_{z_0}^* \sLog^\theN(1))
	= H^1_\AH(S, \sH^\vee \otimes i_{z_0}^* \sLog^\theN(1))
$$
corresponds through \eqref{equation: isom five} to the element
$$
	\pol^{\theN}_{z_0} := \gamma_1^\vee \otimes (v_1 - \gamma_1) + \gamma_2^\vee \otimes (v_2 - \gamma_2)
$$
in $M_{\bbC} /(M_{\bbR}(1) + F^{1} M_{\bbC})$.    Here $\gamma_1^\vee$, $\gamma_2^\vee$ is a basis of 
$\sH_\bbR^\vee$ dual to $\gamma_1$, $\gamma_2$.
Our main result is the explicit calculation of the image of the above element through the isomorphism
\begin{equation}\label{eq: ext isom}
	 M_{\bbC} /(M_{\bbR}(1) + F^{1} M_{\bbC}) \xrightarrow\cong 
	M_{\bbR} / M_{\bbR} \cap (M_{\bbR}(1) + F^{1} M_{\bbC})
\end{equation}
defined by $u \mapsto  u + \phi_\infty(u)$, where $\phi_\infty$ denotes the complex conjugation on
$M_\bbC$.
In order to state our result, we will use a basis of $i_{z_0}^* \sLog^\theN$ which 
gives the isomorphism of mixed $\bbR$-Hodge structures between $i_{z_0}^* \sLog^\theN$
and $\prod_{k=0}^\theN \Sym^k \sH$.  Let $\ul e' := \ul e - E_{0,1}(z_0) \ul\omega^{*\vee}$
be an element of $i_{z_0}^* \sLog^{(1)}$.   We let $\ul e^{m,n} := \ul e'^{a} \ul\omega^{\vee m}
\ul\omega^{*\vee n}/a!$ for $a = \theN - m - n$, which form a basis of  $i_{z_0}^* \sLog^\theN$.
The complex conjugation acts on this base by $\phi_\infty(\ul e^{m,n}) = A^{n-m} \ul e^{n,m}$.

As in the $p$-adic case given in Lemma \ref{lem: splitting}, we have the following.

\begin{lemma}[Splitting Principle]\label{lem: split hodge}
	We have an isomorphism of mixed $\mathbb{R}$-structures
	$$
		i_{z_0}^* \sLog^\theN \cong \prod_{j=0}^\theN \Sym^j \sH
	$$
	given by mapping $\ul e^{m,k}$ to $\ul\omega^{\vee m} \ul\omega^{*\vee k}$.
\end{lemma}


Applying our choice of $G_{m,k}(z_0)$ given in Proposition \ref{prop: relation to Eisenstein}
to the explicit description of $\pol^\theN_{z_0} +\phi_\infty\left(\pol^\theN_{z_0}\right)$,
we now have the following result, originally due to Beilinson and Levin \cite{BL}.
See also the calculations of Wildeshaus (\cite{W} III Theorem 4.8).

\begin{theorem}\label{thm: Beilinson and Levin}
	Let $z_0$ be a non-zero point in $E(\bbC)$, and we let
	$$
		\pol^N_{z_0} :=\left[ \sP^N_{z_0}\right] \in H^1_\AH(S, \sH^\vee \otimes i_{z_0}^* \sLog^\theN(1))
	$$
	be the pull back by $z_0$ of the cohomology class of the elliptic polylogarithm sheaf $\sP^N$.
	Then the image of this cohomology class with respect to the isomorphism
	$$
		H^1_\AH(S, \sH^\vee \otimes i_{z_0}^* \sLog^\theN(1)) 
		\xrightarrow\cong  M_\bbR/ M_{\bbR } \cap  (M_{\bbR}(1) + F^{1} M_\bbC) 
	$$
	for $\sM: = \sH^\vee \otimes i_{z_0}^* \sLog^\theN$ is given by 
	$$
		\sum_{m=1}^\theN \sum_{k=0}^{\theN-m}  \frac{(-1)^{k}}{A^{-m}} 
		e^*_{-m, k+1}(z_0) \ul\omega \otimes \ul e^{m,k}\\
		  %
		  %
		  %
		  %
		+\sum_{m=0}^\theN \sum_{k=1}^{\theN-m}  \frac{(-1)^{k-1}}{A^{-m-1}} 
		e^*_{-m-1,k}(z_0)
		  \uomegas \otimes \ul e^{m,k},
	$$
	where the $e^*_{a,b}(z_0)$ are the Eisenstein-Kronecker numbers defined in Definition \ref{def: EK}.
\end{theorem}

The $p$-adic analogue of the above result is given in Theorem \ref{thm: main}.


\end{document}